\documentclass[10pt]{article}
  \usepackage{epsfig,amsmath,amssymb,amsthm,cancel,bm,floatflt}
\usepackage{graphicx}
\usepackage{setspace} 
\usepackage{array}
\usepackage{pstricks,pstricks-add,pst-math,pst-xkey}

\newtheorem{theorem}{Theorem}

\newtheorem{corollary}[theorem]{Corollary}

\newtheorem{property}[theorem]{Property}

\newcommand{\real}{I\!\!R}

\newcommand{\btwobytwo}{\left [ \begin{array}{rr}}
\newcommand{\etwobytwo}{ \end{array} \right ]}
\newcommand{\bthrbythr}{\left [ \begin{array}{rrr}}
\newcommand{\ethrbythr}{ \end{array} \right ]}

\title{Primal and dual algorithms\\ for the minimum covering Euclidean ball \\ 
of a set of Euclidean balls in $\real^n$}

\author{P. M. Dearing,     Mark Cawood\\
School of Mathematical and Statistical Sciences\\ Clemson University}

\begin{document}
\maketitle

\begin{abstract}

\noindent  Primal and dual algorithms are developed for solving the $n$-dimensional convex optimization problem of finding the Euclidean ball of minimum radius
that covers $m$ given Euclidean balls, each with given center and radius.  Each algorithm is based on a directional search method
in which a search path may be a ray or a two-dimensional conic section in $\real^n$.
At each iteration, a search path is constructed  by  
the intersection of bisectors of pairs of points, where the bisectors are either hyperplanes or $n$-dimensional hyperboloids.
The optimal step size along each search path is determined explicitly.   
\end{abstract}

\vspace{.5in} \noindent Keywords:  minimum covering ball, one-center
location, min-max location, nonlinear programming

\vspace{.7in} \noindent Corresponding author:  P. M. Dearing,
pmdrn@clemson.edu, School of Mathematical and Statistical Sciences, Clemson
University, 520 Bentbrook Lane, Clemson, SC 29631, USA

\newpage

\section{Introduction}

Let $P=\{\mathbf{p}_1,\ldots,\mathbf{p}_m\}$ be a given set  of $m$
distinct points in $\real^n$.
For each point $\mathbf{p}_i  \in P$, 
let $r_i$ be a non-negative radius, and let  
 $[\mathbf{p}_i, r_i] = \{ \mathbf{x}: \| \mathbf{x} - \mathbf{p}_i \|  = r_i \}$
denote  the Euclidean  ball  where $\mathbf{p}_i$ is the center, $r_i$ is the  radius, and $\| \mathbf{x} - \mathbf{p}_i \|$
is the Euclidean distance between   $\mathbf{p}_i$ and $\mathbf{x}$.
The problem of finding the \textbf{ minimum covering Euclidean  ball of a set of Euclidean balls} is to determine the ball
$[\mathbf{x}^*, z^*]$, with center $\mathbf{x}^*$ and minimum radius $z^*$ 
 that  covers, or contains, $[\mathbf{p}_i, r_i]$  for all  $\mathbf{p}_i \in P$.
The  problem is denoted by $M(P)$ and is written as:
\[ \begin{array}{llll} \text{$M(P)$:} & \min &z \\
& \text{s.t.} & z \geq \; \! \| \mathbf{x} - \mathbf{p}_i \| + r_i, & \mathbf{p}_i \in P.\\
\end{array} \]
If all the radii are equal, then all the radii may be assumed to be zero, 
and  the problem is to find the minimum covering ball of the points $\mathbf{p}_i \in P$, 
called the\textbf{ minimum covering ball problem}.
Problem $M(P)$ has the following equivalent representation:
\[  \begin{array}{llll} \text{$M(P)$:} & \underset{\mathbf{x} \in \real^n}{\min}
& \underset{\mathbf{p}_i \in P}{\max} \; \!\parallel\mathbf{p}_i-\mathbf{x}\parallel + r_i.  \\
\end{array}
\]
This version of the problem is known as the \textbf{ min-max location
problem with fixed distance} and as the \textbf{one-center delivery problem}. The center
$\mathbf{x}$ is the location of a facility that minimizes the maximum service to points
 $\mathbf{p}_i \in P$, where service is measured as the travel distance from $\mathbf{x}$ to
$\mathbf{p}_i$ plus a fixed travel distance $r_i$.

This paper presents  primal and  dual algorithms for solving problem $M(P)$.
 At each iteration of either algorithm the problem structure is used to construct a search path 
 that is either a ray or a two-dimensional conic section in $\real^n$,
 and to determine the optimal  step size along the search path.
 
\section{Literature}

The problem of finding the minimum covering circle of a set of points in $\real^2$
 was first posed by Sylvester \cite{Syl1} in 1857. 
 Various geometric solutions were reported by  Sylvester \cite{Syl2},  Crystal \cite{Chrystal}, Blumenthal and Whalen \cite{Blumenthal}, 
 Rademacher and Toeplitz \cite{Rademacher}, and Elzinga and Hearn \cite{EH1}.
Voronoi diagrams \cite{AR} have also been used to solve the problem in $\real^2$.
 Megiddo \cite{M1} reported an algorithm  for $M(P)$, with $n = 3$, that is linear in $m$.

For the  minimum covering ball problem of a set of points in $\real^n$, 
Elzinga and Hearn \cite{EH2} solved the dual problem  as a convex, quadratic
programming problem. Hopp and Reeve \cite{Hop} extended the Sylvester and
Crystal algorithm to $\real^n$; however, their results relied
on a heuristic without  proof of convergence. Fischer, G\"artner
and Kutz \cite{Fischer} expanded the Hopp and Reeve algorithm, gave a proof of finite convergence, 
and were able to solve problems for $m$ up to 10,000 and $n$  up to 2,000. 
Meggido \cite{M2} extended his approach in \cite{M1} to $\real^n$, which yielded an  algorithm that is linear in $m$ but exponential in $n$.
Dyer \cite{Dyer} improved the time-complexity of Meggido's algorithm.
Dearing and Zeck \cite{Dearing} reported a dual algorithm based on search paths
constructed from bisectors of pairs of points.
Cavaleiro and Alizadeh \cite{Cavaleiro} present computational improvements to the Dearing and Zeck approach. 

The literature for the minimum covering ball of a set of balls in $\real^n$ is more limited.
Elzinga and Hearn  \cite{EH1}  presented a geometrical algorithm for the  problem in $\real^2$.
Xu, Freund and Sun \cite{xu} reported computational results for four general approaches to problem $M(P)$ in $\real^2$: a second-order cone reformulation, a sub-gradient approach, a 
quadratic programming scheme, and a randomized incremental algorithm.
Two approaches for problem $M(P)$  in $\real^n$ were reported by
Zhou, Toh and Sun \cite{Zhou}:  an unconstrained convex program whose objective function approximates the maximum objective function, 
and a reformulation of $M(P)$ as a 
second order cone programming problem.  They solve problems with $n$ up to 10,000, and $m$ up to 5,000.
Applications of problem $M(P)$ are referenced in \cite{Fischer} and \cite{M1}.  
Plastria \cite{Plastria}  presents a survey of the min-max location problems.

The approach taken here develops primal and dual algorithms that are based based on  properties and  structures of problem $M(P)$.
Both algorithms are examples of the active set method \cite{GMW}, where a set  $S \subset P$  is active if the points in $S$ are affinely independent and 
if the constraints of $M(P)$ corresponding to points in S hold at equality.
At each iteration, a search path is constructed from the intersections of bisectors of pairs of points in the active set $S$.
For the primal algorithm, primal feasibility and complementary slackness conditions are maintained for points on the search path.
For the dual algorithm, dual feasibility and complementary slackness conditions are maintained for points on the search path.
The step size is determined explicitly for each algorithm.

If the radii are equal for all the balls corresponding to points in $S$, the search path is a ray (the intersection of hyperplanes), but if some of the radii are unequal 
for the balls corresponding to points in $S$, the search path is a two-dimensional conic section (the intersection of hyperboloids).


\section{Properties of Problem $M(P)$}

The following properties of problem $M(P)$ are used to  develop  the primal and dual  algorithms.
The first property follows from the statement of  problem $M(P)$ as a min-max location problem.
\begin{property}
The  objective function $ \underset{\mathbf{p}_i \in P}{\max} \; \| \mathbf{x} - \mathbf{p}_i \| + r_i$,
is continuous and strictly convex, so  
that a minimum solution  $\mathbf{x}^*$ exists and is unique. \hfill $\Box$
\end{property}

Observe that  a  ball $[\mathbf{p}_k,  r_k]$  is contained in a ball   $[\mathbf{p}_j,  r_j]$, 
if and only if  $r_j \geq \; \| \mathbf{p}_{k} - \mathbf{p}_j \| + r_k $.
In this case, any ball that covers $[\mathbf{p}_j,  r_j]$ also covers $[\mathbf{p}_k,  r_k]$, so that 
$[\mathbf{p}_k,  r_k]$ is redundant and $\mathbf{p}_k$ can be eliminated  from the set $P$.  
Furthermore, the ball $[\mathbf{p}_j,  r_j]$ is the \textbf{trivial solution}  for $M(P)$ if
$r_j \geq  \;\|\mathbf{p}_{k} - \mathbf{p}_j \| + r_k $ for all $\mathbf{p}_k \in P$.
Subsequently, it is assumed for problem $M(P)$ that no ball is contained in any other ball, that is,
\begin{equation}
 r_j < \; \|\mathbf{p}_{k} - \mathbf{p}_j \| + r_k   \hspace{.3in} \text{for each pair} \hspace{.3in} \mathbf{p}_{j}, \mathbf{p}_{k} \in P.  \label{eqn1}
\end{equation}
Under Assumption (1), problem $M(P)$ has no redundant points and hence no trivial solution.

For each pair of points $\mathbf{p}_{j}, \mathbf{p}_{k} \in P$,
with corresponding radii $r_{j}$ and $r_{k}$,
 the \textbf{bisector} of $\mathbf{p}_{j}$ and $ \mathbf{p}_{k}$,
is the set:    
\begin{equation}
B_{j,k} = \{ \mathbf{x}: \| \mathbf{x} - \mathbf{p}_j \|+ r_j= 
 \| \mathbf{x} - \mathbf{p}_k \| + r_k  \}. 
 \label{eqn2}\end{equation}
 If $r_j = r_k$, the bisector $B_{j,k}$  is the hyperplane that is orthogonal to the line through $\mathbf{p}_{j}$ and $\mathbf{p}_{k}$, and contains the midpoint between $\mathbf{p}_{j}$ and $\mathbf{p}_{k}$.
If  $r_j \neq r_k$, and if Assumption (1) holds, the bisector $B_{j,k}$ is one sheet of a  hyperboloid  of two sheets in $\real^n$,  
whose foci are the points  $\mathbf{p}_{j}$ and $\mathbf{p}_{k}$, and whose parameters are determined by $\mathbf{p}_{j}, \mathbf{p}_{k}, r_j$ and $r_k$.

 Bisectors of pairs of points in $P$ provide a key structure for developing 
the solution algorithms presented here.
 The following  necessary condition for optimality is part of this structure.
 
 \begin{property} If Assumption (1) holds, and if  $[\mathbf{x}^*,z^*]$ is optimal to problem $M(P)$, 
 then $\mathbf{x}^*$ is on at least one bisector.
 \end{property} 
 \noindent \textbf{Proof}:  Assume to the contrary that the optimal center is not on any bisector, 
 so that $z^* = \| \mathbf{x}^* - \mathbf{p}_j \| + r_j$ for exactly one point $\mathbf{p}_j \in P$,
 and $z^* >  \| \mathbf{x}^* - \mathbf{p}_k \| + r_k$  for all  $\mathbf{p}_k \neq \mathbf{p}_j $.
 Define the ray  $X = \{ \mathbf{x}(\alpha) = \mathbf{x}^* + \alpha(\mathbf{p}_j - \mathbf{x}^*), 0 \leq \alpha \leq 1 \}$,
 and let $z_k(\alpha) = \| \mathbf{x}(\alpha) - \mathbf{p}_k \| + r_k$, for each $\mathbf{p}_k \in P$.
 For $\mathbf{p}_j$, $z_j(\alpha)$ is continuous for $0 \leq \alpha \leq 1$, with $z_j(0) = z^*$, and $z_j(1) = r_j$.
 For $\mathbf{p}_k \neq \mathbf{p}_j$, $z_k(0) = \|\mathbf{x}^* - \mathbf{p}_k \| + r_k < z^*$, and 
 $z_k(1) = \| \mathbf{p}_j - \mathbf{p}_k \| + r_k > r_j$ by Assumption 1.
 Thus for each $\mathbf{p}_k \neq \mathbf{p}_j$, there exists $\alpha_{j,k} > 0$, such that $z_k(\alpha_{j,k}) = z_j(\alpha_{j,k})$.
 Let $\alpha^* = \min_{\mathbf{p}_k \neq \mathbf{p}_j} \{ \alpha_{j,k} \}$.
 Then $z_j(\alpha^*) < z^*$ and $z_k(\alpha^*)  < z^*$, which contradicts the optimality of $[\mathbf{x}^*, z^*]$. \hfill $\Box$
 
 A subset $S$  of $P$ is an \textbf{active set} \cite{GMW} corresponding to a ball $[\mathbf{x}_S, z_S]$ if the points in $S$ are affinely independent, and if
the constraint of $M(P)$ corresponding to each point in $S$ holds at equality, that is,
$ \| \mathbf{p}_i - x \| + r_i = z$, for each $\mathbf{p}_i \in S$.
The primal and dual algorithms developed here are examples of the active-set method \cite{GMW} of constrained optimization. 
At each iteration the search path is determined by the points in $S$ and their radii.  
At each iteration the size of $S$ may increase, decrease or remain unchanged.
 
Property 3.2 and Assumption (\ref{eqn1}) lead to the following result.
\begin{property} If $S$ is an active set corresponding to the covering ball  $[\mathbf{x}_S, z_S]$ and $|S| > 1$,
then $\mathbf{x}_S \neq \mathbf{p}_i$\; for all $\mathbf{p}_i \in S$.
\end{property}
\noindent \textbf{Proof}:  
If $\mathbf{x}_S = \mathbf{p}_i$ for some $\mathbf{p}_i \in S$, 
then $r_i = z_S$, and
$r_i = \parallel \mathbf{p}_i - \mathbf{p}_j \parallel + r_j$  for each $\mathbf{p}_j \in S$, 
which violates  Assumption (\ref{eqn1}).  \hfill $\Box$

Given an active set $S$, 
Property 3.3 shows that for each $\mathbf{p}_i \in S$,  the norm  $ \| \mathbf{x} - \mathbf{p}_i \| $ 
is differentiable with respect to $\mathbf{x}$ over the domain  $\real^n \setminus S$.
Several properties of $M(P)$ follow from the Karush Kuhn Tucker (KKT)  conditions \cite{GMW} which are stated next.
\begin{property} Let $[\mathbf{x}_S, z_S]$ be a ball corresponding to the active set $S$.
Then the ball $[\mathbf{x}_S, z_S]$ is optimal to $M(P)$ if and only if 
there exists  variables $\lambda_i \geq 0$ for $\mathbf{p}_i \in S$ and  $\lambda_i = 0$ for $\mathbf{p}_i \notin S$ such that 
\begin{alignat}{2}
&z_S \geq   \| \mathbf{x}_S - \mathbf{p}_i \|  + r_i  && \mathbf{p}_i \in P  \label{kktch1}  \\
&\sum_{\mathbf{p}_i \in S} \lambda_i = 1&&   \label{kktch2}  \\
&\sum_{\mathbf{p}_i \in S}  \frac{(\mathbf{x}_S - \mathbf{p}_i) }{ \| \mathbf{x}_S - \mathbf{p}_i \|  }\lambda_i  = \mathbf{0} &&  \label{kktch3}  \\
&\lambda_i \geq 0  &&\mathbf{p}_i \in S   \label{kktch4}    \\
&(z_S -   \| \mathbf{x}_S - \mathbf{p}_i \|  - r_i )\lambda_i = 0 \qquad && \mathbf{p}_i \in P. \label{kktch5}
\end{alignat}
\end{property}

\begin{property} At optimality, statements (\ref{kktch2}), (\ref{kktch3}) and (\ref{kktch4}) of the KKT conditions are equivalent to: 
\begin{alignat}{2}
&\sum_{\mathbf{p}_i \in S} \pi_{i}= 1 & &  \label{denom1}  \\
&\sum_{\mathbf{p}_i \in S}  (\mathbf{x}_S - \mathbf{p}_i)\pi_i  = \mathbf{0} & &   \label{denom2}  \\
&\pi_{i} \geq 0 &&\mathbf{p}_i \in S.  \label{denom3} 
\end{alignat}
\end{property}
\noindent \textbf{Proof}:  
Since  $\parallel \mathbf{x}_S - \mathbf{p}_i \parallel > 0$, for  all $\mathbf{p}_i \in S$, 
the change of variables
 $\pi_i =\frac{\lambda_i/ \| \mathbf{x}_S - \mathbf{p}_i \|  }{\sum_{\mathbf{p}_i \in S} \lambda_i/ \| \mathbf{x}_S - \mathbf{p}_i \|  }$
 and $\lambda_i = \frac{  \| \mathbf{x}_S - \mathbf{p}_i \|  \pi_i}{ \sum_{\mathbf{p}_i \in S}  \| \mathbf{x}_S - \mathbf{p}_i \|  \pi_i}$
 for each $\mathbf{p}_i \in S$, shows the equivalence of conditions
 (\ref{kktch2}), (\ref{kktch3}), and (\ref{kktch4})  to conditions 
  (\ref{denom1}), (\ref{denom2}), and (\ref{denom3}). \hfill $\Box$

Properties 3.4 and 3.5 lead to additional properties of an optimal ball $[\mathbf{x}_S, z_S]$.
\begin{property} The minimum covering ball  $[\mathbf{x}_S,z_S]$ for problem $M(P)$ is determined by an active set $S$ of at most $n+1$ affinely independent points in $P$, 
and the optimal center $\mathbf{x}_S$ lies in the convex hull of $S$, denoted by conv$(S)$.  
Furthermore, $\mathbf{x}_S$ is in the relative interior of conv$(S)$, denoted by ri$(S)$ if and only if
$S$ is a minimal active set for which the KKT conditions are satisfied at $[\mathbf{x}_S, z_S]$.
\end{property}
\noindent \textbf{Proof:}
Conditions  (\ref{denom1}) and (\ref{denom2}) determine a linear system with $n+1$ linear equations. A solution is determined by at most $n+1$
linearly independent columns of the system which correspond to at most $n+1$ affinely independent points from $P$.  Conditions 
 (\ref{denom1}), (\ref{denom2})  and (\ref{denom3}) are equivalent to $\mathbf{x}_S \in$ conv$(S)$.  Finally, each 
 column of (\ref{denom1}) and (\ref{denom2}) with  $\pi_i = 0$ can be eliminated from $S$, so that $\pi_i > 0$ for all $\mathbf{p}_i \in S,$ if and only if
 $\mathbf{x}_S \in$ ri$(S),$ if and only if $S$ is a minimal active set satisfying the KKT conditions at $[\mathbf{x}_S,z_S]$. \hfill $\Box$

If $S$ is an active set of size $s$ corresponding to an optimal covering ball $[\mathbf{x}_S, z_S]$, there may be more than $s$ constraints that are active.  
Property 3.6 shows that the  constraints corresponding  to points in an active set $S$ are sufficient to determine the optimal covering ball $[\mathbf{x}_S, z_S]$.

If  $S$ is an active set corresponding to the ball  $[\mathbf{x}_S, z_S]$, then 
$\mathbf{x}_S$ is on the bisector $B_{j,k}$ for each pair of points 
$\mathbf{p}_j, \mathbf{p}_k \in S$.  
Consequently, $\mathbf{x}_S$ is on the intersection of bisectors  $B_{j,k}$ over all pairs of points 
$\mathbf{p}_j, \mathbf{p}_k \in S$, denoted by $B_S$.
That is, $\mathbf{x}_S \in B_S = \cap_{\{\mathbf{p}_j, \mathbf{p}_k\} \subseteq S} B_{j,k}$.

Each search path constructed here for the primal and dual algorithms 
is a space curve contained in $B_S,$ so that each point on the search path is 
also on each bisector for all pairs of points in $S$.
Thus, the complementary slackness conditions (\ref{kktch5}) are maintained at each point on the search path. 

Reference \cite{Dearing1} presents results on intersections of hyperplanes with $n$-dimensional conic sections, 
along with expressions for computing the parameters and vectors of the intersections.
Results from \cite{Dearing1} that are applicable to problem $M(P)$ are stated in the Appendix.
A key property is that the non-empty  intersection of two $n$-dimensional hyperboloids with a common focal point is contained in a hyperplane, and that the
intersection of the two hyperboloids is equivalent to the intersection of either hyperboloid with the hyperplane.  
Using this property, $B_S$ is equivalent  to the intersection of one bisector (a hyperboloid)  with $s-2$ hyperplanes, 
so that $B_S$ is a conic section of dimension $n-s+2$.  
The Appendix includes closed form expressions for computing the  parameters and vectors of $B_S$.
 
If all the points in $S$ have equal radii, then all the bisectors are hyperplanes, and the intersection $B_S$ is a hyperplane of dimension $n-s+1$.
However, if some pair of points in $S$ have unequal radii, $B_S$ is 
a conic section of dimension $n - s+2$.


\section{Primal Algorithm}

At each iteration of the primal algorithm there is an active set $S$ and a  corresponding ball $[ \mathbf{x}_S, z_S ]$
that  satisfy expressions \eqref{kktch1} and \eqref{kktch5} of the KKT conditions,  
but  do not satisfy expressions \eqref{kktch2}, \eqref{kktch3} and \eqref{kktch4} of the KKT conditions.
That is, $S$ and $[ \mathbf{x}_S, z_S ]$ are primal feasible but not dual feasible.

Since $[ \mathbf{x}_S, z_S ]$ is primal feasible at each iteration, $z_S$ is an upper bound on the optimal objective function value of $M(P)$.
 Assuming non-degeneracy, the radius $z_S$ will be shown to decrease at each iteration of the primal algorithm.

The primal algorithm is initialized  by choosing $\mathbf{x}_S$ to be any point in $\real^n$, computing
 $z_S =  \max_{\mathbf{p}_i \in P} \| \mathbf{x}_S - \mathbf{p}_i \| + r_i $, and choosing $S = \{ \mathbf{p}_j \}$,  where $\mathbf{p}_j$
 is some point such that  $z_S = \| \mathbf{x}_S - \mathbf{p}_j \| + r_j$.
Observe that $[ \mathbf{x}_S, z_S ]$ is primal feasible, but $\mathbf{x}_S \notin$ aff$(S)$, which implies  $\mathbf{x}_S \notin $ conv$(S)$ so that $\mathbf{x}_S$ is not dual feasible.  

\noindent \begin{bfseries} Primal Search Phase \end{bfseries}
 
Given an active set $S$ and a primal feasible ball $[ \mathbf{x}_S, z_S ]$, a search path  $X_S = \{ \mathbf{x}(\alpha) : \alpha \geq 0 \}$ is constructed
 so that  $X_S \subset B_S$, and  $\mathbf{x}_S = \mathbf{x}(0)$.
Thus,  $S$ is an active set for  the ball  $[\mathbf{x}(\alpha), z(\alpha)]$, for $\alpha \geq 0$, 
where $z(\alpha) = \|\mathbf{x}(\alpha)-\mathbf{p}_i \| +r_i$, for $\mathbf{p}_i \in S$.
Also, the complementary slackness conditions \eqref{kktch5} are satisfied by $S$ and $[\mathbf{x}(\alpha), z(\alpha)]$, for $\alpha \geq 0$ .  
If all the points in $S$ have equal radii, the search path will be a ray, but if some points in $S$ have unequal radii, the search path will be a two-dimensional conic section in $\real^n$.


 \noindent 
\begin{bfseries}Case 1: All points in $S$ have equal radii  \end{bfseries}

\noindent \textbf{The Search Path:} The points in $S$ are denoted by  $S = \{ \mathbf{p}_{i_1},  \ldots, \mathbf{p}_{i_{s}} \}$, 
where by assumption $r_{i_j} = r_{i_1}$ for $\mathbf{p}_{i_j} \in S$.
In this case, the search path is the ray  
\begin{equation}
X_S = \{ \mathbf{x}(\alpha) = \mathbf{x}_S + \alpha \mathbf{d}_S,  \alpha \geq 0 \}, \text{\;\;where\;\;} 
\mathbf{d}_S = (\mathbf{p}_{i_1} - \mathbf{x}_S) - \text{Proj}_{\text{sub}(S)} (\mathbf{p}_{i_1} - \mathbf{x}_S), \label{pds}
\end{equation}  
and sub$(S) = $ span$[ (\mathbf{p}_{i_1} - \mathbf{p}_{i_2}), \ldots, (\mathbf{p}_{i_1} - \mathbf{p}_{i_{s}}) ]$.  
That is, $\mathbf{d}_S$ is the component of $(\mathbf{p}_{i_1} - \mathbf{x}_S)$ that is orthogonal to the projection of 
$(\mathbf{p}_{i_1} - \mathbf{x}_S)$ onto sub$(S)$.
If  $s = 1$, as in the initial step with $S = \{ \mathbf{p}_{i_1} \}$, $\mathbf{d}_S = \mathbf{p}_{i_1} -\mathbf{x}_S$. 

The next property shows that for a positive step size along the search path,
the set $S$ remains active and the radius of the ball is decreased.
\begin{property}
Suppose the set $S = \{ \mathbf{p}_{i_1},  \ldots, \mathbf{p}_{i_{s}} \}$  is an active set for 
 the primal feasible ball $[\mathbf{x}_S, z_S ]$, with $r_{i_1} = r_{i_j}$ for each $\mathbf{p}_{i_j} \in S$.
Let $X_S = \{ \mathbf{x}(\alpha) = \mathbf{x}_S + \alpha \mathbf{d}_S, \alpha \geq 0 \}$ where
$\mathbf{d}_S = (\mathbf{p}_{i_1} - \mathbf{x}_S) -$ Proj$_{sub(S)} (\mathbf{p}_{i_1} - \mathbf{x}_S)$.
Then $S$ is an active set for the ball $[\mathbf{x}(\alpha), z(\alpha) ]$, for $ \alpha \geq 0$, where 
 $z(\alpha) = \|\mathbf{x}(\alpha) - \mathbf{p}_{i_1} \| + r_{i_1}$.
 Furthermore,  $z(\alpha)$ is decreasing for 
$0 \leq \alpha \leq \hat{\alpha}_S$, 
where $\hat{\alpha}_S = \frac{(  \mathbf{p}_{i_1} - \mathbf{x}_S )  \mathbf{d}_S}{\| \mathbf{d}_S \|^2} > 0$.
\end{property}
 \noindent \textbf{Proof:}  For each pair of points $\mathbf{p}_{i_1}$ and $\mathbf{p}_{i_j}$ in $S$, with $r_{i_1} = r_{i_j}$, the bisector  $B_{i_1,i_j} = \{ \mathbf{x}: (\mathbf{p}_{i_1} - \mathbf{p}_{i_j})\mathbf{x} = (\mathbf{p}_{i_1} - \mathbf{p}_{i_j})\mathbf{x}_S \}$ since $\mathbf{x}_S \in B_{i_1,i_j}$. Then
 $(\mathbf{p}_{i_1} - \mathbf{p}_{i_j})\mathbf{x}(\alpha) =  (\mathbf{p}_{i_1} - \mathbf{p}_{i_j})\mathbf{x}_S + \alpha  (\mathbf{p}_{i_1} - \mathbf{p}_{i_j})\mathbf{d}_S
 = (\mathbf{p}_{i_1} - \mathbf{p}_{i_j})\mathbf{x}_S$, since $\mathbf{d}_S$ is orthogonal to $(\mathbf{p}_{i_1} - \mathbf{p}_{i_j})$ by construction.
 Thus, for all $\mathbf{p}_{i_j} \in S$, and for $\alpha \geq 0$, 
 $\mathbf{x}(\alpha) \in B_{i_1,i_j}$, that is,
$z(\alpha) = \|\mathbf{x}(\alpha) - \mathbf{p}_{i_1} \| + r_{i_1} = \|\mathbf{x}(\alpha) - \mathbf{p}_{i_j} \| + r_{i_j}$,
 so that $S$ is an active set for $[\mathbf{x}(\alpha), z(\alpha) ]$.
 
By the construction of $\mathbf{d}_S$,  $(  \mathbf{p}_{i_1} -  \mathbf{x}_S) \mathbf{d}_S > 0$.
 Property 3.3 implies that $\mathbf{x}(\alpha) \neq \mathbf{p}_{i_1}$ 
for all $\alpha$, so that $z(\alpha)$ is differentiable with respect  to $\alpha$.
Furthermore, $z(\alpha)$ is convex and has a unique minimum at  $\hat{\alpha}_S = \frac{(  \mathbf{p}_{i_1} -  \mathbf{x}_S) \mathbf{d}_S}{\| \mathbf{d}_S \|^2} > 0$,
so that $z(\alpha)$ is decreasing along $X_S$ for $0 \leq \alpha \leq \hat{\alpha}_S$.   \hfill $\Box$
 
 \begin{corollary}  Given  $S = \{ \mathbf{p}_{i_1} \}$, let  $X_S = \{ \mathbf{x}(\alpha) = \mathbf{x}_S + \alpha \mathbf{d}_S, \alpha \geq 0 \}$ 
 where $\mathbf{d}_S = (\mathbf{p}_{i_1} - \mathbf{x}_S)$.  
 Then  $z(\alpha) = \| \mathbf{x}(\alpha) - \mathbf{p}_{i_1} \| + r_{i_1}$ is decreasing for $0 \leq \alpha \leq \hat{\alpha}_S$,  where $\hat{\alpha}_S = \frac{(  \mathbf{p}_{i_1} - \mathbf{x}_S )  \mathbf{d}_S}{\| \mathbf{d}_S \|^2} > 0$.
 \end{corollary}

 \noindent \begin{bfseries}The Step Size: \end{bfseries}
 For each $\mathbf{p}_k \in P \setminus S$,  the step size $\alpha_k \geq 0$ is determined, if it exists,  so that $X_S$ intersects the bisector $B_{i_1,k}$ at $\mathbf{x}(\alpha_k)$.   
 If $\mathbf{x}(\alpha_k) \in X_S \cap B_{i_1,k}$, then  $\mathbf{x}(\alpha_k) \in X_S \cap B_{i_j,k}$ for all $\mathbf{p}_j \in S$.
  That is, $X_S$  simultaneously intersects the bisectors $B_{i_j,k}$  at $\mathbf{x}(\alpha_k)$, for all $\mathbf{p}_{i_j} \in S$.
  Thus, it suffices  to consider the intersection of $X_S$ with only $B_{i_1,k}$.

 At the point $\mathbf{x}(\alpha_k) \in X_S \cap B_{i_1,k}$, 
   $z(\alpha_k) =  \| \mathbf{x}(\alpha_k) - \mathbf{p}_k \| + r_k = \| \mathbf{x}(\alpha_k) - \mathbf{p}_{i_1} \| + r_{i_1} $,
   so that the constraint corresponding to the point $\mathbf{p}_k$ is active.
   Geometrically, the search moves the center $\mathbf{x}(\alpha)$ of the ball along the path $X_S$, while decreasing the radius $z(\alpha)$.

 There are two sub-cases to consider for computing $\alpha_k$ depending on whether the radius $r_{i_1}$  equals the radius $r_{k}$.\\
  \noindent Sub-case 1: $r_{i_1} = r_{k}$, so that $B_{i_1,k}$ is a hyperplane.  The intersection of $X_S$ and $B_{i_1,k}$ is determined by solving for $\alpha$ 
using the equation
$(\mathbf{p}_{i_1} - \mathbf{p}_k)\mathbf{x}(\alpha) = (\mathbf{p}_{i_1} - \mathbf{p}_k)(\mathbf{p}_{i_1} + \mathbf{p}_k)/2$.
If $(\mathbf{p}_{i_1} - \mathbf{p}_{k}) \mathbf{d}_S = 0$, set $\alpha_k = \infty$, else set
 \begin{alignat}{2}
\alpha' =  \frac{(\mathbf{p}_{i_1} - \mathbf{p}_{k}) (\mathbf{p}_{i_1} + \mathbf{p}_{k})/2 - (\mathbf{p}_{i_1} - \mathbf{p}_{k})  \mathbf{x}_S}
 {(\mathbf{p}_{i_1} - \mathbf{p}_{k})  \mathbf{d}_S}. \label{PEQR1}
 \end{alignat}
If $\alpha' <0$, set $\alpha_k = \infty$, else set $\alpha_k = \alpha'$.\\
 \noindent Sub-case 2:   $r_{i_1} \neq r_{k}$,  so that $B_{i_1,k}$ is a hyperboloid. 
 The intersection of $X_S$ and $B_{i_1,k}$ is determined by substituting the expression for  
$ \mathbf{x}(\alpha)$ into the quadratic form \eqref{Q} for the hyperboloid $B_{i_1,k}$, with center  $ \mathbf{c}_{i_1,k}$, axis vector $ \mathbf{v}_{i_1,k}$, 
 eccentricity $\epsilon_{i_1,k}$, and parmeters $ a_{i_1,k}$ and $c_{i_1,k}$.  
This  gives the quadratic equation 
\begin{equation}
A_{k}\alpha^2 + B_{k}\alpha + C_{k} = 0, \label{quadp}
\end{equation}
where  $ A_{k} = (\mathbf{d}_S)^2 - \epsilon_{i_1,k}^2(\mathbf{d}_S \mathbf{v}_{i_1,k})^2$,
$B_{k} = 2(\mathbf{x}_S - \mathbf{c}_{i_1,k})  \mathbf{d}_S - 2\epsilon_{i_1,k}^2[(\mathbf{x}_S - \mathbf{c}_{i_1,k})  \mathbf{v}_{i_1,k}][ \mathbf{d}_S  \mathbf{v}_{i_1,k}]$, 
and $C_{k} = (\mathbf{x}_S - \mathbf{c}_{i_1,k})^2 - \epsilon_{i_1,k}^2[(\mathbf{x}_S - \mathbf{c}_{i_1,k})  \mathbf{v}_{i_1,k}]^2
 + a_{i_1,k}^2 - c_{i_1,k}^2.$
If  $A_k = 0$, there is one real solution, $\alpha_k = -C_k/B_k$.
If there are no real solutions, $\alpha_k = \infty$.
Otherwise, let $\alpha'$ and $\alpha''$ be the real solutions.
If $A_k < 0$ and $r_{i_1} > r_k$, then $\alpha_k = \max \{ \alpha', \alpha'' \}$, or if $r_{i_1} < r_k$, then $\alpha_k = \min \{ \alpha', \alpha'' \}$.
If $A_k > 0$ and  $r_{i_1} > r_k$, then $\alpha_k = \min \{ \alpha', \alpha'' \}$, or if $r_{i_1} < r_k$, then $\alpha_k = \infty$.
 
 From the Proof of Property 4.1,  $\mathbf{x}(\hat{\alpha}_S)  \in  X_S\; \cap$ aff$(S)$, where 
\begin{equation}
\hat{\alpha}_S = \frac{(\mathbf{p}_{i_1} - \mathbf{x}_S)\mathbf{d}_S}{\| \mathbf{d}_S \|^2} \label{alphats}
\end{equation}
The step size in the Primal Search phase is given by
\begin{equation}
\alpha_S = \min \{ \hat{\alpha}_S,  \alpha_k : \mathbf{p}_k \in P \setminus S  \}. \label{palphas}
\end{equation}  
If $\alpha_S = \hat{\alpha}_S$,  then $ \mathbf{x}(\alpha_S) \in $ aff$(S)$.  In this case, the Update Phase is entered
to check the optimality of the active set $S$ and the ball $ [ \mathbf{x}(\alpha_S), z(\alpha_S) ]$.

Otherwise,  let $\mathbf{p}_e$ be a point  in $P \setminus S$  so that $\alpha_S = \alpha_{e}$. 
Then $\alpha_S$  is the smallest step size so that  the ball $[\mathbf{x}(\alpha_S), z(\alpha_S) ]$ remains primal feasible and 
 contains   $[ \mathbf{p}_{e}, r_e ]$. 
 The set $S \cup  \{ \mathbf{p}_{e} \}$, and the ball $ [ \mathbf{x}(\alpha_S), z(\alpha_S)]$
 are checked for optimality in the Update Phase.

There may be a tie for the entering point $\mathbf{p}_{e}$, in which case the next iteration of the Search Phase may result in a degenerate step with zero step size.  
Cycling can be avoided by applying an adaptation of Bland's rule that chooses the point with the 
 smallest index among all points that are candidates for entering.

 \noindent 
\begin{bfseries}Case 2: At least two points  in $S$ have unequal radii \end{bfseries}

\noindent \textbf{The Search Path:}  The points in $S$ are ordered by non-increasing radii and denoted by 
$S = \{ \mathbf{p}_{i_1},  \ldots,\mathbf{p}_{i_{s}} \}$, with $r_{i_1}  \geq \ldots \geq r_{i_{s}}$.
By the assumption of unequal radii, $r_{i_1} > r_{i_{s}}$, so that $B_{i_1,i_{s}}$ is one sheet of a hyperboloid.  
A search path $X_S$ is constructed to be a two-dimensional conic section such that $X_S \subset B_S  = \cap_{\{ \mathbf{p}_{i_j},\mathbf{p}_{i_k} \} \subset S} B_{i_j,i_k}$.
Properties A.7 and A.8 show that $B_S$ is the intersection of $B_{i_1,i_s}$ with $s-2$ hyperplanes, resulting in a conic section of dimension $n-s+2$..
The vectors and parameters of $B_S$, namely
the center $\mathbf{c}_{S}$, axis vector $\mathbf{v}_{S}$,   vertex $\mathbf{a}_{S}$,   eccentricity  $e_{S}$,  and parameters  $a_{S}$ and  $b_{S}$,
are computed using expressions  (\ref{HT12})  - (\ref{pnc}) in the Appendix.
If $\epsilon_S >1$, then $B_S$ is one sheet of a hyperboloid; if $\epsilon_S <1$, $B_S$ is an ellipsoid; if $\epsilon_S = 1$, $B_S$ is a paraboloid.

Property A.9 in the Appendix shows that for any vector $\mathbf{u}_{S}$ orthogonal to 
$\mathbf{v}_{S}$, there is a two-dimensional conic section $Y_S \subset B_S$ that has the same vectors and parameters as $B_{S}$. 
Furthermore, $Y_S \subset$ aff$(\mathbf{c}_{S}, \mathbf{v}_{S}, \mathbf{u}_{S})$.
Thus, $Y_S$ is either one sheet of a hyperbola, an ellipse, or a parabola if $\epsilon_S >1, \epsilon_S <1$, or $\epsilon_S =1$, respectively.
The search path $X_S$ is constructed to be a sub-path of $Y_S$ determined by a particular choice  of $\mathbf{u}_S$.

\noindent \textbf{Case 2a: $X_S$ is a hyperbola:}
If  $\epsilon_S > 1$, then $B_S$ is one sheet of a hyperboloid, and 
\begin{equation}
Y_S =  \{ \mathbf{y}(\beta) =  \mathbf{c}_S + a_S \sec(\beta) \mathbf{v}_S  + b_S \tan(\beta) \mathbf{u}_S : -\pi/2 < \beta < \pi/2 \} \label{hyp2d}
\end{equation}
is one sheet of a two-dimensional hyperbola  in  aff$(\mathbf{c}_S, \mathbf{v}_S, \mathbf{u}_S)$.
The vector  $\mathbf{u}_S$ is  chosen to be the normalized component of 
$(\mathbf{c}_S - \mathbf{x}_S)$ that is orthogonal to the projection of $(\mathbf{c}_S - \mathbf{x}_S)$  onto $\mathbf{v}_S$. That is,
\begin{equation}
\mathbf{u}_S = [(\mathbf{c}_S - \mathbf{x}_S) -  ((\mathbf{c}_S - \mathbf{x}_S)\mathbf{v}_S)\mathbf{v}_S)] / \|  (\mathbf{c}_S - \mathbf{x}_S) -  (( \mathbf{c}_S -  \mathbf{x}_S )\mathbf{v}_S)\mathbf{v}_S) \|.  \label{pus}
\end{equation}
Property A.8 shows that $Y_S$ is a subset of $B_S$, and that
$\mathbf{y}(0) = \mathbf{a}_S$ is the vertex of the hyperbola $Y_S$ and of the hyperboloid $B_S$.
 Furthermore,  the point   $\mathbf{x}_S = \mathbf{y}(\beta_S) \in Y_S$, where 
 $\beta_S = -\tan^{-1} \{\sqrt{(\mathbf{c}_S - \mathbf{x}_S)^2-a_S^2} / ( \epsilon_Sa_S)\} < 0$.
 As $\beta$ increases from $\beta_S$ to $0$, the point $\mathbf{y}(\beta)$ moves from $\mathbf{x}_S$ toward $\mathbf{a}_S$ 
  and dual feasibility.
 
 The search path $X_S \subset Y_S$ is defined by  
 \begin{equation}
 X_S = \{ \mathbf{x}(\alpha) =  \mathbf{y}(\alpha + \beta_S )   \text{  for   } 0 \leq \alpha \leq - \beta_S \}. \label{xpthhyp}
 \end{equation}
 Observe that  $\mathbf{x}(0) =  \mathbf{y}(\beta_S ) = \mathbf{x}_S$, and $\mathbf{x}(-\beta_S) = \mathbf{y}(0) = \mathbf{a}_S  \in$ aff$(S)$. 
 For $ 0 \leq \alpha \leq  -\beta_S $, the points $\mathbf{x}(\alpha)$ are on $B_S$ and move along the path $X_S$ from $\mathbf{x}_S$ to $\mathbf{a}_S$, 
that is, toward aff$(\mathbf{c}_S, \mathbf{v}_S, \mathbf{u}_S)$ and dual feasibility.

 \begin{property}
Suppose the set  $S =  \{ \mathbf{p}_{i_1}, \ldots,  \mathbf{p}_{i_{s}} \}$ is an active set for the  primal feasible ball $[\mathbf{x}_S, z_S]$,
ordered so that $r_{i_1} \geq \ldots \geq r_{i_{s}}$, with $r_{i_1} > r_{i_{s}}$.  
Let $X_S = \{ \mathbf{x}(\alpha) = \mathbf{y}(\alpha+\beta_S), 0 \leq \alpha \leq -\beta_S \}$  
be the search path along a hyperbola $Y_S$ determined by (\ref{hyp2d}) and (\ref{xpthhyp}) .
Then $S$ is an active set for the ball $[\mathbf{x}(\alpha), z(\alpha) ]$, for $\alpha \geq 0$, where
$z(\alpha) = \| \mathbf{x}(\alpha) - \mathbf{p}_{i_1} \| + r_{i_1} $.
Furthermore, $z(\alpha)$ is decreasing for $0 \leq \alpha \leq - \beta_S$.
 \end{property}
 \noindent \textbf{Proof:} By construction,
 $\mathbf{x}(\alpha) \in B_{i_1,i_j}$ for each $\mathbf{p}_{i_j} \in S$ and $0 \leq \alpha \leq -\beta_S$, so that
 $z(\alpha) = \|\mathbf{x}(\alpha) - \mathbf{p}_{i_1} \| + r_{i_1}  =  \| \mathbf{x}(\alpha) - \mathbf{p}_{i_j} \| + r_{i_j}$.
 Thus, $S$ is an active set for $[\mathbf{x}(\alpha), z(\alpha) ]$, for $\alpha \geq 0$.
 
 Direct computation shows that $z(\alpha)$ has a minimum value at $\alpha = -\beta_S$, and since $z(\alpha)$ is convex with respect to $\alpha$,
 $z(\alpha)$ is decreasing from $0$ to $-\beta_S$.  \hfill  $\Box$

 \noindent \begin{bfseries}The Step Size: \end{bfseries}
 For each $\mathbf{p}_k \in P \setminus S$, the  step size $\alpha_k \geq 0$  is determined, if it exists, so that $X_S$ intersects the bisector $B_{i_1,k}$ at $\mathbf{x}(\alpha_k)$. 
  If  $\mathbf{x}(\alpha_k) \in X_S \cap B_{i_1,k}$, then  $\mathbf{x}(\alpha_k) \in X_S \cap B_{i_j,k}$, for all $\mathbf{p}_{i_j} \in S$. That is,
  $X_S$  simultaneously intersects the bisectors $B_{i_j,k}$  at $\mathbf{x}(\alpha_k)$, for all $\mathbf{p}_{i_j} \in S$.  
  Thus, it suffices  to consider the intersection of $X_S$ with only $B_{i_1,k}$.

    At the point $\mathbf{x}(\alpha_k) \in X_S \cap B_{i_1,k}$, 
   $z(\alpha_k) =  \| \mathbf{x}(\alpha_k) - \mathbf{p}_k \| + r_k = \| \mathbf{x}(\alpha_k) - \mathbf{p}_{i_1} \| + r_{i_1}$,
   so that the constraint corresponding to the point $\mathbf{p}_k$ is active.
Geometrically, the search moves the center $\mathbf{x}_S$ of the ball $[\mathbf{x}_S, z_S]$ along the path $X_S$, while decreasing the radius $z_S$.
  
 The hyperplane    $H_k = \{ \mathbf{x}:  \mathbf{h}_k \mathbf{x} = \mathbf{h}_k \mathbf{d}_k \}$ is constructed such that
 $X_S \cap B_{i_1,k} = X_S \cap H_k$.
Order the set of three points 
$\{ \mathbf{p}_{i_1}, \mathbf{p}_{i_{s}}, \mathbf{p}_k  \}$   by non-increasing radii, and denote the ordered set
by $\{ \mathbf{p}_{j_1}, \mathbf{p}_{j_{2}}, \mathbf{p}_{j_3}  \}$, so that $r_{j_1} \geq r_{j_2} \geq r_{j_3}$.
By assumption, $r_{j_1} > r_{j_3}$.
The vectors  $\mathbf{h}_k$ and $\mathbf{d}_k$ 
of the hyperplane $H_k$ are 
computed as follows. 
\begin{alignat}{4}
& \text{if\;\;} r_{j_1} = r_{j_2} > r_{j_3}\;\; &  \mathbf{h}_{k} &\;=\;  (\mathbf{p}_{j_1} -  \mathbf{p}_{j_2})/\|\mathbf{p}_{j_1} -  \mathbf{p}_{j_2}\| , \notag
\;\;\;\; \mathbf{d}_k = .5(\mathbf{p}_{j_1} + \mathbf{p}_{j_2})\\ \notag
& \text{if\;\;} r_{j_1} > r_{j_2} = r_{j_3}\;\; &  \mathbf{h}_{k} &\;=\;  (\mathbf{p}_{j_2} -  \mathbf{p}_{j_3})/\|\mathbf{p}_{j_2} -  \mathbf{p}_{j_3}\| \notag \notag,
\;\;\;\; \mathbf{d}_k = .5(\mathbf{p}_{j_2} + \mathbf{p}_{j_3}) \\
& \text{if\;\;} r_{j_1} > r_{j_2} > r_{j_3}\;\;  & \mathbf{h}_{k} &\;=\;  (\epsilon_{j_1,j_2}\mathbf{v}_{j_1,j_2} -  \epsilon_{j_1,j_3}\mathbf{v}_{j_1,j_3})/
\|\epsilon_{j_1,j_2}\mathbf{v}_{j_1,j_2} -  \epsilon_{j_1,j_3}\mathbf{v}_{j_1,j_3}\|  \label{hrr} \\
&&\mathbf{w}_{k} & \;=\; (\mathbf{h}_{k} - (\mathbf{h}_{k} \mathbf{v}_{j_1,j_3})\mathbf{v}_{j_1,j_3}) / 
\| \mathbf{h}_{k} - (\mathbf{h}_{k} \mathbf{v}_{j_1,j_3})\mathbf{v}_{j_1,j_3})\|  \notag \\
&&\mathbf{d}_{k} & \;=\;  \mathbf{d}_{j_1,j_3} + \frac{\mathbf{v}_{j_1,j_2} (\mathbf{d}_{j_1,j_2} - \mathbf{d}_{j_1,j_3})}{\mathbf{v}_{j_1,j_2} \mathbf{w}_{k}} \mathbf{w}_{k}.   \notag 
\end{alignat}
The intersection $X_S \cap H_k$ is equivalent to the intersection of $Y_S \cap H_k$, which is determined 
by solving the equation $\mathbf{h}_k\mathbf{y}(\beta_k) = \mathbf{h}_k\mathbf{d}_k$ for $\beta_k$, which is equivalent to solving the equation
\begin{equation}
A \sec(\beta_k) + B \tan(\beta_k) = C,  \label{AAAu} 
\end{equation}
where \;
$A = a_S \mathbf{h}_k \mathbf{v}_S, $ 
\; $B = b_S\mathbf{h}_k \mathbf{u}_S, $ 
\;\;and \; $C =  \mathbf{h}_k  ( \mathbf{d}_k -  \mathbf{c}_S). $

Let $D_h = \sqrt{C^2-(A^2-B^2)}$. 
The value $C^2$ is the squared distance between  $H_k$ and the center $\mathbf{c}_S$. The value   $A^2-B^2$ is the
minimum squared distance between $\mathbf{c}_S$ and $H_k$ such that $H_k \cap Y_S \neq \emptyset$.
Thus, $H_k \cap Y_S \neq \emptyset$ if and only if $C^2 - (A^2 - B^2) \geq 0$, 

If $C^2 - (A^2 - B^2) > 0$, there are two  solutions, $\beta_1$ and $\beta_2$, determined by
\begin{alignat}{2}
& \tan(\beta_1) = \frac {-AB -CD_h} {AC - BD_h}, \;\;\;\text{and}\;\;\; \sec(\beta_1) = \frac {B^2 + C^2} {AC - BD_h} \notag \\
&\text{or} \label{hysol} \\
&  \tan(\beta_2) = \frac {-AB +CD_h} {AC + BD_h}, \;\;\;\text{and}\;\;\; \sec(\beta_2) = \frac {B^2 + C^2} {AC + BD_h}. \notag 
\end{alignat}
If $\sec(\beta_1) > 0$, then $\mathbf{y}(\beta_1) \in Y_S$, and  $\beta_1 := \tan^{-1}\left(  \frac {-AB -CD_h} {AC - BD_h} \right)$.
If  $\sec(\beta_1) < 0$, then $\mathbf{y}(\beta_1) \notin Y_S$, and $\beta_1 := \infty$.
Similarly, if $\sec(\beta_2) > 0$, then $\beta_2 := \tan^{-1}\left(  \frac {-AB + CD_h} {AC + BD_h} \right)$,
 and if  $\sec(\beta_2) < 0$, $\beta_2 := \infty$.  Then, $\beta_k = \min \{ \beta_1, \beta_2 \}$.

If $C^2 - (A^2 - B^2) = 0$, there is one solution, $\beta_k := \tan^{-1}(-B/C)$.
If $C^2 - (A^2 - B^2) < 0$, set $\beta_k:= \infty$.



The step size $\alpha_k$ along the path $X_S$ from $\mathbf{x}_S$ is given by $\alpha_k = \beta_k - \beta_S$. 
If  $0 \leq \alpha_k \leq -\beta_S$, then $\mathbf{x}(\alpha_k)$ is the  point of intersection of $X_S$ and $B_{i_1,k}$. 
 If $\alpha_k < 0$, then set  $\alpha_k = \infty$, since this is a step in the opposite direction. 
If $\alpha_k > -\beta_S$, then set $\alpha_k = \beta_S$, since $\mathbf{x}(\alpha_k)$ is on the opposite side of aff$(S)$ from $\mathbf{x}_S$.

\noindent \textbf{Case 2b: $X_S$ is an ellipse:}
 If $\epsilon_S < 1$, then $B_S$ is an ellipsoid, and
\begin{equation}
Y_S =  \{ \mathbf{y}(\beta) =  \mathbf{c}_S + a_S \cos(\beta) \mathbf{v}_S  + b_S \sin(\beta) \mathbf{u}_S : 0 < \beta < 2\pi \}. \label{elp2d}
\end{equation}  
is a two-dimensional ellipse in aff$(\mathbf{c}_S, \mathbf{v}_S, \mathbf{u}_S)$.
The axis vector $\mathbf{v}_S$ is oriented as follows so that the objective function value is decreasing along the search path.  
If $(\mathbf{c}_S - \mathbf{x}_S)\mathbf{v}_S > 0$, reset $\mathbf{v}_S := - \mathbf{v}_S$. 
The vector  $\mathbf{u}_S$ is chosen to be  the normalized component of  $(\mathbf{c}_S  -  \mathbf{x}_S)$ that is orthogonal to 
the projection of $(\mathbf{c}_S  -  \mathbf{x}_S)$ onto $\mathbf{v}_S$.  That is,
\begin{equation}
\mathbf{u}_S = [(\mathbf{c}_S  -  \mathbf{x}_S)   -  ((\mathbf{c}_S  -  \mathbf{x}_S)\mathbf{v}_S)\mathbf{v}_S)] / \|  (\mathbf{c}_S  -  \mathbf{x}_S) -  (( \mathbf{c}_S  -  \mathbf{x}_S )\mathbf{v}_S)\mathbf{v}_S) \|. \label{use}
\end{equation}
Property A.9 shows that $Y_S$ is a subset of $B_S$, and that
 the points $\mathbf{y}(0)$ and $\mathbf{y}(\pi)$ are the vertices of the ellipse $Y_S$ and of the ellipsoid $B_S$.
 Furthermore, the point   $\mathbf{x}_S = \mathbf{y}(\beta_S) \in Y_S$, where $\beta_S$ is computed as follows:
  $\beta_S = -\cos^{-1} \{\sqrt{(\mathbf{c}_S - \mathbf{x}_S)^2-b_S^2} / ( \epsilon_Sa_S)\} < 0$.
As $\beta$ increases from $\beta_S$ to $0$, the point $\mathbf{y}(\beta)$ moves from $\mathbf{x}_S$ to the vertex $\mathbf{y}(0)$, 
 that is,  toward aff$(\mathbf{c}_S, \mathbf{v}_S, \mathbf{u}_S)$ and dual feasibility.
 
 The search path $X_S$ is defined in terms of $Y_S$ by  
 \begin{equation}
 X_S = \{ \mathbf{x}(\alpha) =  \mathbf{y}(\alpha + \beta_S )   \text{  for   } 0 \leq \alpha \leq   - \beta_S \}. \label{xpthelp}
 \end{equation}
 Observe that  $\mathbf{x}(0) =  \mathbf{y}(\beta_S ) = \mathbf{x}_S$, and $\mathbf{x}(  -\beta_S) = \mathbf{y}(0) = \mathbf{a}_S$.
 For $ 0 \leq \alpha \leq  -\beta_S $, the points $\mathbf{x}(\alpha)$ are on $B_S$ and move along the path $X_S$ from $\mathbf{x}_S$ toward the vertex $\mathbf{a}_S$, 
that is, toward aff$(\mathbf{c}_S, \mathbf{v}_S, \mathbf{u}_S)$ and dual feasibility. 
  
 For the case when $X_S$ is an ellipse, the following property and its proof are analogous to Property 4.2 for the case when $X_S$ is a hyperbola.
 \begin{property}
Given the active set  $S =  \{ \mathbf{p}_{i_1}, \ldots,  \mathbf{p}_{i_{s}} \}$  for the  primal feasible ball $[\mathbf{x}_S, z_S]$,
ordered so that $r_{i_1} \geq \ldots \geq r_{i_{s}}$, with $r_{i_1} > r_{i_{s}}$.  
Let $X_S = \{ \mathbf{x}(\alpha) = \mathbf{y}(\alpha+\beta_S), 0 \leq \alpha \leq -\beta_S \}$  
be the search path along the ellipse  $\mathbf{y}(\alpha)$ determined by (\ref{elp2d}) and (\ref{xpthelp}) .
Then $S$ is an active set for the ball $[\mathbf{x}(\alpha), z(\alpha) ]$, for $\alpha \geq 0$, where
$z(\alpha) = \| \mathbf{x}(\alpha) - \mathbf{p}_{i_1} \| + r_{i_1} $.
Furthermore, $z(\alpha)$ is decreasing for $0 \leq \alpha \leq  - \beta_S$.
 \end{property}

 \noindent \begin{bfseries}The Step Size: \end{bfseries}
 For each $\mathbf{p}_k \in P \setminus S$, the  step size $\alpha_k \geq 0$  is determined, if it exists, so that $X_S$ intersects the bisector $B_{i_1,k}$ at $\mathbf{x}(\alpha_k)$.  
 If  $\mathbf{x}(\alpha_k) \in X_S \cap B_{i_1,k}$, then  $\mathbf{x}(\alpha_k) \in X_S \cap B_{i_j,k}$, for all $\mathbf{p}_{i_j} \in S$. That is,
  $X_S$  simultaneously intersects the bisectors $B_{i_j,k}$  at $\mathbf{x}(\alpha_k)$, for all $\mathbf{p}_{i_j} \in S$.  
  Thus, it suffices  to consider the intersection of $X_S$ with only $B_{i_1,k}$.
 
  At the point $\mathbf{x}(\alpha_k) \in X_S \cap B_{i_1,k}$, 
   $ \| \mathbf{x}(\alpha_k) - \mathbf{p}_k \| + r_k = \| \mathbf{x}(\alpha_k) - \mathbf{p}_{i_1} \| + r_{i_1}= z(\alpha_k)$,
   so that the constraint corresponding to the point $\mathbf{p}_k$ is active.
 Geometrically, the search moves the center $\mathbf{x}_S$ of the ball $[\mathbf{x}_S, z_S]$ along the path $X_S$, while decreasing the radius $z_S$.

  The hyperplane $H_k = \{ \mathbf{x}: \mathbf{h}_k\mathbf{x} = \mathbf{h}_k\mathbf{d}_k \}$
  is constructed so that   $X_S \cap H_k = B_{i_1,k} \cap H_k$.
   Order the set of three points 
$\{ \mathbf{p}_{i_1}, \mathbf{p}_{i_{s}}, \mathbf{p}_k  \}$   by non-increasing radii, and denote the ordered set
by $\{ \mathbf{p}_{j_1}, \mathbf{p}_{j_{2}}, \mathbf{p}_{j_3}  \}$, so that $r_{j_1} \geq r_{j_2} \geq r_{j_3}$.
By assumption, $r_{j_1} > r_{j_3}$.
The vectors  $\mathbf{h}_k$ and $\mathbf{d}_k$ 
of the hyperplane $H_k = \{ \mathbf{x}:  \mathbf{h}_k \mathbf{x} = \mathbf{h}_k \mathbf{d}_k \}$ are 
computed by  (\ref{hrr}).

The intersection  $X_S \cap B_{i_1,i_j}$ is equivalent to the intersection of $Y_S \cap H_k$, which is 
determined by solving the equation $\mathbf{h}_k\mathbf{y}(\beta_k) =  \mathbf{h}_k\mathbf{d}_k$ for $\beta_k$, which is equivalent to solving the equation
\begin{equation}
A \cos(\beta_k) + B \sin(\beta_k) = C,  \label{AAAv} 
\end{equation}
where \;
$A = a_S \mathbf{h}_k \mathbf{v}_S, $ 
\; $B = b_S \mathbf{h}_k \mathbf{u}_S, $ 
\;\;and \; $C =  \mathbf{h}_k  ( \mathbf{d}_k-  \mathbf{c}_S). $

Let $D_e = \sqrt{A^2 +B^2-C^2}$. 
The value $C^2$ is the squared distance between  $H_k$ and the center $\mathbf{c}_S$. The expression   $A^2+B^2$ is the
maximum squared distance between $\mathbf{c}_S$ and $H_k$ such that $H_k \cap Y_S \neq \emptyset$.
Thus, $H_k \cap Y_S \neq \emptyset$ if and only if $A^2 + B^2 -C^2 \geq 0$.

If $A^2 + B^2 -C^2 \geq 0$, there are two possible solutions, $\beta_1$ and $\beta_2$,  determined by
\begin{alignat}{2}
& \cos(\beta_1) = \frac {AC + BD_e} {A^2 + B^2}, \;\;\;\text{and}\;\;\; \sin(\beta_1) = \frac {BC -AD_e} {A^2 + B^2} \notag \\
&\text{or} \label{elsol} \\
&  \cos(\beta_2) = \frac {AC - BD_e} {A^2 + B^2}, \;\;\;\text{and}\;\;\; \sin(\beta_2) = \frac {BC + AD_e} {A^2 + B^2} \notag 
\end{alignat}
If $\sin(\beta_1) > 0$, then $\mathbf{y}(\beta_1) \in X_S$, and  $\beta_1 := \tan^{-1}\left(  \frac {BC -AD_e} {AC + BD_e} \right)$.
If  $\sin(\beta_1) < 0$, then $\mathbf{y}(\beta_1) \notin X_S$, and $\beta_1 := \infty$.
Similarly, if $\sin(\beta_2) > 0$, then $\beta_2 := \tan^{-1}\left(  \frac {BC + AD_e} {AC - BD_e} \right)$,
 and if  $\sec(\beta_2) < 0$, $\beta_2 := \infty$.  Then, $\beta_k = \min \{ \beta_1, \beta_2 \}$.

If $A^2 + B^2 -C^2 = 0$, there is one solution $\beta_k := \tan^{-1}(B/A)$.
 If $A^2 + B^2 -C^2 < 0$, set $\beta_k:= \infty$.



The step size $\alpha_k$ along the path $X_S$ from $\mathbf{x}_S$ is given by $\alpha_k = \beta_k - \beta_S$. 
If  $0 \leq \alpha_k \leq -\beta_S$, then $\mathbf{x}(\alpha_k)$ is the  point of intersection of $X_S$ and $B_{i_1,k}$. 
 If $\alpha_k < 0$, then set  $\alpha_k = \infty$, since this is a step in the opposite direction. 
If $\alpha_k > -\beta_S$, then set $\alpha_k = \beta_S$, since $\mathbf{x}(\alpha_k)$ is on the opposite side of aff$(S)$ from $\mathbf{x}_S$.

 
\noindent \textbf{Case 2c: $X_S$ is a paraboloid:}
If  $\epsilon_{S} = 1$, then $B_{S} $ is a paraboloid of dimension $n-s+2$, and 
\begin{equation}
 Y_{S}  =  \{ \mathbf{y}(\beta) =  \hat{\mathbf{c}}_{S} + \tilde{c}_{S} \beta^2 \mathbf{v}_{S}  + 2 \tilde{c}_{S}\beta \mathbf{u}_{S} : -\infty < \beta < \infty \} \label{para2d}
 \end{equation}
is a two-dimensional parabola in  aff$( \hat{\mathbf{c}}_{S}, \mathbf{v}_{S}, \mathbf{u}_{S})$, where $\hat{\mathbf{c}}_S$ and  $\tilde{c}_S$ are determined by (\ref{pch}) - (\ref{pnc}).
The vector $\mathbf{u}_{S}$ is  chosen to be the normalized component of  $(\hat{\mathbf{c}}_{S} - \mathbf{x}_S)$ 
that is orthogonal to the projection of $(\hat{\mathbf{c}}_{S} - \mathbf{x}_S)$
 onto sub($S$).  That is, 
 \begin{equation}
\mathbf{u}_{S} =( (\hat{\mathbf{c}}_{S} - \mathbf{x}_S) - \text{Proj}_{\text{sub}(S)} (\hat{\mathbf{c}}_{S} - \mathbf{x}_S)) / \| (\hat{\mathbf{c}}_{S} - \mathbf{x}_S) - \text{Proj}_{\text{sub}(S)} (\hat{\mathbf{c}}_{S} - \mathbf{x}_S) \| \label{parau}.
 \end{equation}

Property A.9 shows that $Y_{S} $ is a subset of $B_{S}$, and that
 the point $\mathbf{y}(0)$ is the vertex of the parabola $Y_{S} $ and of the paraboloid $B_{S}$.
  Let $\beta_S$ be the parameter such that $\mathbf{y}(\beta_S) = \mathbf{x}_S$.  
 If $\mathbf{x}_S = \hat{\mathbf{c}}_S$, then $\beta_S = 0$, but if  $\mathbf{x}_S \neq \hat{\mathbf{c}}_S$, then
 $\beta_S = \sqrt{ -2 +  \sqrt{ 4 +  (\hat{\mathbf{c}}_S - \mathbf{x}_S )^2 /(\tilde{c}_S^2 ) }}.$ 
 If a point was deleted from $S$ in the last iteration, $\mathbf{x}_S \neq \hat{\mathbf{c}}_S$, else, $\mathbf{x}_S = \hat{\mathbf{c}}_S$.
 
  The search path $X_S$ is defined to be a subset of $Y_S$ given by 
 \begin{equation}
 X_S = \{ \mathbf{x}(\alpha) = \mathbf{y}(\alpha + \beta_S), \alpha \geq 0\}. \label{xpthpara}
 \end{equation} 
 Observe that $\mathbf{x}(0) = \mathbf{y}(\beta_S) = \mathbf{x}_S$, and for $\alpha \geq 0$, 
 the points $\mathbf{x}(\alpha)$ are on $B_S$ and move along $X_S$ from $\mathbf{x}_S$ toward $\mathbf{p}_e$ and dual feasibility.

  For the case when $X_S$ is a parabola, the following property and its proof are analogous to Properties 4.2 and 4.3 for the case when $X_S$ is based on a hyperbola or an ellipse, respectively.
 \begin{property}
Given the active set  $S =  \{ \mathbf{p}_{i_1}, \ldots,  \mathbf{p}_{i_{s}} \}$  for the  primal feasible ball $[\mathbf{x}_S, z_S]$,
ordered so that $r_{i_1} \geq \ldots \geq r_{i_{s}}$, with $r_{i_1} > r_{i_{s}}$.  
Let $X_S = \{ \mathbf{x}(\alpha) = \mathbf{y}(\alpha+\beta_S), 0 \leq \alpha \leq -\beta_S \}$  
be the search path along the parabola  $\mathbf{x}(\alpha)$ determined by (\ref{para2d}) and (\ref{xpthpara}) .
Then $S$ is an active set for the ball $[\mathbf{x}(\alpha), z(\alpha) ]$, for $\alpha \geq 0$, where
$z(\alpha) = \| \mathbf{x}(\alpha) - \mathbf{p}_{i_1} \| + r_{i_1} $.
Furthermore, $z(\alpha)$ is increasing for $0 \leq \alpha \leq -\beta_S  \pi/2$.
 \end{property}
  
 \noindent \begin{bfseries}The Step Size: \end{bfseries}
 For each point $\mathbf{p}_k \in P$, the step size $\alpha_k \geq 0$ is determined, if it exists, so that $X_S$ intersects the bisector $B_{i_1,k}$ at $\mathbf{x}(\alpha_k)$.
 If $\mathbf{x}(\alpha_k) \in X_S \cap B_{i_1,k}$, then $\mathbf{x}(\alpha_k) \in X_S \cap B_{i_j,k}$, for all  $\mathbf{p}_{i_j} \in S$, that is, 
$X_S$ simultaneously intersects the bisectors $B_{i_j,k}$ at $\mathbf{x}(\alpha_k)$, for all $\mathbf{p}_{i_j} \in S$.
Thus, it suffices to determine the intersection of $X_S$ with only $B_{i_1,k}$.
 
 At the  point $\mathbf{x}(\alpha_k) \in X_S \cap B_{i_1,k}$, 
 $ z(\alpha_k) = \| \mathbf{x}(\alpha_k) - \mathbf{p}_k \| + r_k =  \| \mathbf{x}(\alpha_k) - \mathbf{p}_{i_1} \| + r_{i_1}$,
 so that the constraint corresponding to the point $\mathbf{p}_k$ is active. 
  Geometrically, the search moves the center  $\mathbf{x}(\alpha)$ along the path $X_{S}$ while increasing  
 $z(\alpha)$.
 
Theorem A.1 shows that $X_S \cap B_{i_1,k} = X_S \cap H_k$, where 
the hyperplane $H_k =  \{ \mathbf{x} : \mathbf{h}_k\mathbf{x} = \mathbf{h}_k\mathbf{d}_k \}$ is constructed as follows.
Order the set of three points $\{\mathbf{p}_{i_1}, \mathbf{p}_{i_{s}}, \mathbf{p}_{k} \}$, 
by non-increasing radii, and denote the ordered set by $\{ \mathbf{p}_{j_1},  \mathbf{p}_{j_2},   \mathbf{p}_{j_3} \}$, so that $r_{j_1} \geq r_{j_2} \geq r_{j_3}$.
By assumption, $r_{j_1} > r_{j_3}$.
The vectors $\mathbf{h}_k$ and $\mathbf{d}_k$ of the hyperplane
$H_k = \{ \mathbf{x} : \mathbf{h}_k\mathbf{x} = \mathbf{h}_k\mathbf{d}_k \}$, are computed by the expressions (\ref{hrr}).
 
The intersection of $X_{S}$ and $H_k$ is equivalent to the intersection of $Y_{S}$ and $H_k$, and is determined by  
by solving the equation $\mathbf{h}_k\mathbf{y}(\beta_k) = \mathbf{h}_k\mathbf{d}_k$, for $\beta_k$,
which yields the quadratic equation  
\begin{equation}
A \beta_k^2+ B \beta_k = C,  \label{AAAP} 
\end{equation}
where \;
$A = \tilde{c}_{S} \mathbf{h}_{k}  \mathbf{v}_{S}, $ 
\; $B = \tilde{c}_{S} \mathbf{h}_{k}  \mathbf{u}_{S}, $ 
\;\;and \; $C =  \mathbf{h}_{k}   ( \mathbf{d}_{k} -  \hat{\mathbf{c}}_{S}). $
If $Y_{S} \cap H_{e} \neq \emptyset$, there are possibly two intersection points
corresponding to the following two solutions to equation  (\ref{AAAP}):

If $B^2+A^2-C^2 <0$, there is no solution and no intersection point, and $\beta_k := \infty$.  
If $B^2+A^2-C^2 =0$, or if $C = 0$, then there is one solution $\beta_k$.
Else, $B^2+A^2-C^2  > 0$, and there are two solutions  $\beta_1$ and $\beta_2$.
If $\beta_1 < \beta_S$, or $\mathbf{h}_k\mathbf{y}(\beta_1) \neq \mathbf{h}_e\mathbf{d}_k$, $\beta_1 := \infty$. 
If $\beta_2 < \beta_S$, or $\mathbf{h}_k\mathbf{y}(\beta_2) \neq \mathbf{h}_e\mathbf{d}_k$, $\beta_2 := \infty$. 
Then $\beta_k := \min\{\beta_1, \beta_2 \}$.
The step size $\alpha_k$ along the path $X_S$ from $\mathbf{x}_S$ is given by $\alpha_k = \beta_k - \beta_S$.

Given the step size $\alpha_k$ for each $\mathbf{p}_k \in P\setminus S$, 
the step size in the Primal Search phase is given by
 $\alpha_S = \min \{ \beta_S,  \alpha_k : \mathbf{p}_k \in P \setminus S  \}$.  
If $\alpha_S = \beta_S$,  then $ \mathbf{x}(\alpha_S) \in $ aff$(S)$.  In this case, $[ \mathbf{x}_S, z_S ]: = [ \mathbf{x}(\alpha_S), z(\alpha_S) ]$ and the set $S$
are checked for optimality in the Update Phase.

Otherwise,  let $\mathbf{p}_e$ be a point  in $P \setminus S$  so that $\alpha_S = \alpha_{e}$. 
Then $\alpha_S$  is the smallest step size from $\mathbf{x}_S$ so that  the ball $[\mathbf{x}(\alpha_S), z(\alpha_S) ]$ remains primal feasible and 
contains the ball $[\mathbf{p}_e,r_e]$.    
The set  $S: = S \cup \{\mathbf{p}_e \}$ and the solution $[ \mathbf{x}_S, z_S ]: = [ \mathbf{x}(\alpha_S), z(\alpha_S) ]$ are checked for optimality in the Update Phase.

There may be a tie for the entering point $\mathbf{p}_{e}$, in which case the next iteration of the Search Phase may result in a degenerate step with zero step size.  
It can be shown that cycling will be avoided by an adaptation of Bland's rule that chooses the point with the 
 smallest index among all points that are candidates for entering.


\noindent \begin{bfseries} Primal Update Phase \end{bfseries}
 
In the Update Phase the  set $S$ and the ball $[ \mathbf{x}_S, z_S ]$
 are checked for optimality  using
equations (\ref{kktch2}) and (\ref{kktch3}) of the KKT conditions.
If a solution exists with $\lambda_{i_j} \geq 0$ for all $\mathbf{p}_{i_j} \in S$, then the ball $[ \mathbf{x}_S, z_S ]$ is dual feasible, 
and hence optimal to $M(P)$.

If a solution exists, but some $\lambda_{i_l} < 0$, then $ \mathbf{x}_S \in $ aff$(S)$ 
 but $ \mathbf{x}_S \notin $ conv$(S)$, so that 
$[ \mathbf{x}_S, z_S ]$ is not optimal.
In this case, a point $\mathbf{p}_{i_l} \in S$, such that $\lambda_{i_l} < 0$, is chosen to leave $S$.
The  algorithm returns to the Update Phase with the set 
$S:= S  \setminus \{\mathbf{p}_{i_l} \}$
and the ball $[ \mathbf{x}_S, z_S ]$.

If no solution exists, then $\mathbf{x}_S \notin$ aff$(S)$, and  $[ \mathbf{x}_S, z_S ]$ is not optimal. In this case,  the points in  $S$ are affinely independent, and
the algorithm returns to the Search Phase with the set $S$ 
and the ball $[ \mathbf{x}_S, z_S ]$.

\begin{property} If the points in $S$ are affinely dependent, and if the equations (\ref{kktch2}) and (\ref{kktch3}) of the KKT conditions, 
with $\mathbf{x}^* = \mathbf{x}(\alpha_S)$ and $S$ have a solution with  $\lambda_{i_l} < 0$, then the 
points in $S \setminus \{\mathbf{p}_{i_l} \}$ are affinely independent.
\end{property}
\noindent \textbf{Proof:}
The points in $S \setminus \{ \mathbf{p}_e \}$ are affinely independent if and only if the columns of the linear system  determined by (\ref{kktch2}) and (\ref{kktch3}) 
over the set $S  \setminus \{ \mathbf{p}_e \} $  are linearly independent.  
The linear system determined by  (\ref{kktch2}) and (\ref{kktch3}) 
over the set $S \setminus \{\mathbf{p}_{i_l} \}$ has the same columns as the linear system over the set $S  \setminus \{ \mathbf{p}_e \}$ 
except that the column corresponding to $\mathbf{p}_{i_l}$ is replaced by the column corresponding to $\mathbf{p}_e$. 
 Since the multiplier $\lambda_{i_l} \neq 0$, the columns over the set $S \setminus \{\mathbf{p}_{i_l} \}$ are linearly independent,
which implies that the points in set 
$S \setminus \{\mathbf{p}_{i_l} \}$ are affinely independent. \hfill $\Box$

The next property shows that if the point $\mathbf{p}_{i_l}$ leaves $S$, because $\lambda_{i_l} < 0$, 
then the point $\mathbf{p}_{i_l}$ is covered by the ball $[ \mathbf{x}(\alpha),z(\alpha) ]$ during the next search phase. 

\begin{property} In the Update Phase, suppose the point $\mathbf{p}_{i_l}$ is chosen to leave the set $S$ 
because $\lambda_{i_l} < 0$  in the solution to equations (\ref{kktch2}) and (\ref{kktch3}) of the KKT conditions over the set $S$.  
Let $X_{S  \setminus \{\mathbf{p}_{i_l} \}}$  denote the search path determined by $S \setminus \{\mathbf{p}_{i_l} \}$,  
which may be a ray, a hyperbola determined by (\ref{hyp2d}), an ellipse determined by (\ref{elp2d}), or a parabola.
In either case, the  ball $[\mathbf{x}(\alpha), z(\alpha)]$  remains feasible with  respect to the leaving point $\mathbf{p}_{i_l}$, 
that is, $[\mathbf{p}_{i_l}, r_{i_l}] \subset [\mathbf{x}(\alpha), z(\alpha)]$ for $0 \leq \alpha \leq -\beta_{S  \setminus \{\mathbf{p}_{i_l} \}}$.
\end{property}
\noindent \textbf{Proof:}
 Let $z_{i_j}(\alpha) = \|\mathbf{x}(\alpha) - \mathbf{p}_{i_j} \| +r_{i_j}$, for each $\mathbf{p}_{i_j} \in S$,
 and let $z'_{i_j}(\alpha)$ denote the directional derivative 
 of $z_{i_j}(\alpha)$ at the point $\mathbf{x}(\alpha) \neq \mathbf{p}_{i_j}$ on the path $X_{S  \setminus \{\mathbf{p}_{i_l} \}}$.  That is, 
$z'_{i_j}(\alpha) = \nabla z_{i_j}(\alpha)\mathbf{x}'(\alpha) = \frac{(\mathbf{x}(\alpha) -  \mathbf{p}_{i_j} )}{\|\mathbf{x}(\alpha) - \mathbf{p}_{i_j} \| } \mathbf{x}'(\alpha)$,
where $\mathbf{x}'(\alpha)$ is the normalized tangent vector to the search path at $\mathbf{x}(\alpha)$.
By construction of the search path $X_{S  \setminus \{\mathbf{p}_{i_l} \}}$, $z(\alpha) = z_{i_j}(\alpha)$,  for each  $\mathbf{p}_{i_j} \in S  \setminus \{\mathbf{p}_{i_l} \}$, 
and for $0 \leq \alpha \leq -\beta_{S \setminus \{\mathbf{p}_{i_l} \}}$.  This implies the directional derivatives are equal,
that is, $z'_{i_1}(\alpha) =z'_{i_j}(\alpha)$,  for each  $\mathbf{p}_{i_j} \in S  \setminus \{\mathbf{p}_{i_l} \}$, and for
 $0 \leq \alpha \leq -\beta_{S \setminus \{\mathbf{p}_{i_l} \}}$.
Property 4.2 shows that $z_{i_j}(\alpha)$ is decreasing along the path $X_{S  \setminus \{\mathbf{p}_{i_l} \}}$, that is, 
$\nabla z'_{i_j}(\alpha) < 0,$ for each  $\mathbf{p}_{i_j} \in S  \setminus \{\mathbf{p}_{i_l} \}$, 
and $0 \leq \alpha \leq -\beta_{S \setminus \{\mathbf{p}_{i_l} \}}$.

Equation (\ref{kktch3}) is written as   $\sum_{\mathbf{p}_{i_j} \in S \setminus \{\mathbf{p}{i_l} \}} \frac{ (\mathbf{x}_S - \mathbf{p}_{i_j} )}{ \| \mathbf{x}_S - \mathbf{p}_{i_j} \| } \lambda_{i_j} +\frac{ (\mathbf{x}_S - \mathbf{p}_{i_l} )}{ \| \mathbf{x}_S - \mathbf{p}_{i_l} \| }(\lambda_{i_l} ) = 0.$
Substituting $\mathbf{x}(0) = \mathbf{x}_S$, multiplying each summand by the tangent vector $\mathbf{x}'(0)$,
substituting $z'_{i_1}(0) =z'_{i_j}(0)$,  for each  $\mathbf{p}_{i_j} \in S  \setminus \{\mathbf{p}_{i_l} \}$, and applying equation (\ref{kktch2}),
equation (\ref{kktch3})  may be written as $z'_{i_1}(0) (1-\lambda_{i_l}) = z_{i_l}'(0)( -\lambda_{i_l} )$.

Since $z'_{i_1}(0) < 0$ and $\lambda_{i_l} < 0$, then $z'_{i_l}(0) < 0$, so that $z_{i_l}(\alpha)$ is decreasing at $\alpha = 0$. 
Furthermore  $z'_{i_1}(0) = z'_{i_l}(0)(-\lambda_{i_l})/(1- \lambda_{i_l}) > z'_{i_l}(0)$, so that $z_{i_l}(\alpha)$
  is decreasing at a faster rate than $z_{i_j}(\alpha)$ for   $\mathbf{p}_{i_j} \in S$. This shows that 
   $[\mathbf{p}_{i_l}, r_{i_l}] \subset [\mathbf{x}(\alpha), z(\alpha)]$ for $0 \leq \alpha \leq -\beta_{S \setminus \{\mathbf{p}_{i_l} \}}$. \hfill $\Box$ 

\section{Pseudocode for Primal Algorithm}

\begin{enumerate}
\item Input: a set of distinct points $P=\{{\mathbf{p}_1,...,\mathbf{p}_m}\} \subset \real^n$,\\
 a radius $r_i \geq 0$, a  ball $[\mathbf{p}_i, r_i] = \{ \mathbf{x} : \parallel \mathbf{x} - \mathbf{p}_i \parallel \leq r_i \}$ for each $\mathbf{p}_i \in P$. \\
Output: a unique  ball $[\mathbf{x}^*, z^*]$, with minimum radius $z^*$, containing  $[\mathbf{p}_i, r_i]$ for $\mathbf{p}_i$ in $P$.\\
Assume condition (1) for each pair $ \mathbf{p}_{j}, \mathbf{p}_{k} \in P.$

\item Initialize: Choose any point $ \mathbf{x}_S \in \real^n$,  compute $z_S = \max_{\mathbf{p}_i \in P}\|\mathbf{x}_S - \mathbf{p}_i \| + r_i$,\\
Choose $S = \{ \mathbf{p}_j \} \ni z_S = \|\mathbf{x}_S - \mathbf{p}_j \| + r_j $
 
\item Search Path: Denote the active set  $S = \{ \mathbf{p}_{i_1}, \ldots, \mathbf{p}_{i_{s}} \}$ 
ordered so that  $r_{i_1} \geq \ldots \geq r_{i_{s}}$.
	 \begin{enumerate}  
	\item[3 a:]  $r_{i_1} = r_{i_s}$ so that all the points in $S$ have equal radii. \\ 
	 The search path is $X_S$ is given by  (\ref{pds}).  \\
	 Find the step size $\alpha_k$ so that  $\mathbf{x}(\alpha_k) \in X_S \cap B_{i_1,k}$, for $\mathbf{p}_k \in P \setminus S$. 
		 \begin{enumerate}  
		 \item[] If $r_{i_1} = r_k$, then $B_{i_1,k}$ is a hyperplane.
			 Determine the step size $\alpha_k$ using  (\ref{PEQR1}).
		 \item[] If $r_{i_1} \neq r_k$, then $B_{i_1,k}$ is a hyperboloid. Determine $\alpha_k$ by solving  (\ref{quadp}).
		 \end{enumerate}
	\item[] Determine $\hat{\alpha}_S$ using (\ref{alphats}).
	\item[] Choose $\alpha_S = \min \{ \hat{\alpha}_S, \alpha_{k} : \mathbf{p}_{k} \in P \setminus S \}$, 
		\begin{enumerate}
		\item[] If $\alpha_S = \hat{\alpha}_S$, go to Step 4 with active set $S$ and $[\mathbf{x}_S, z_S ] = [\mathbf{x}(\alpha_S), z(\alpha_S) ]$.
		\item[] Else, chose $\mathbf{p}_e \in P \setminus S$, so that $\alpha_S = \alpha_e$.  
		\item[] Reset $S: = S \cup \{ \mathbf{p}_e \}$, and  $[\mathbf{x}_S, z_S ] := [\mathbf{x}(\alpha_S), z(\alpha_S) ]$. Go to Step 4.
		\end{enumerate}
	
	\item[ 3 b:] $r_{i_1} > r_{i_{s}}$, so that the points in $S$ have unequal radii.\\
	Determine  $\mathbf{v}_S$, $\mathbf{c}_S$, $a_S$,  $b_S$, and $\epsilon_S$ of $B_S = \cap_{j=2}^s B_{i_1.i_{j}}$ 
	using (\ref{HT12})  - (\ref{pnc}).
		\begin{enumerate}
		\item[3 b (1):] If $\epsilon_S >1$, $B_S$ is a hyperboloid of dimension $n-s+2$.
			\begin{enumerate}
			\item[] If $\mathbf{x}_S = \mathbf{a}_S$, then  $\beta_S = 0$; else compute $\beta_S$.\\
			Construct the  hyperbola $Y_S$ using  (\ref{hyp2d})
			where $\mathbf{u}_S$ is determined by (\ref{pus}).\\
			Construct the search path $X_S$  usiing (\ref{xpthhyp}) .\\
			Find the step size $\alpha_k$ so that  $\mathbf{x}(\alpha_k) \in X_S \cap B_{i_1,k}$, for $\mathbf{p}_k \in P \setminus S$.\\
			Find the step size $\beta_k$ so that $\mathbf{h}_k\mathbf{y}(\beta_k) =\mathbf{h}_k\mathbf{d}_k$, where $\mathbf{h}_k$ is given by (\ref{hrr}), \\
			and $\beta_k$ is determined by (\ref{AAAu}) and (\ref{hysol}).  \\
			Then $\alpha_k = \beta_k - \beta_S$.
			\end{enumerate}
		\item[3 b (2):] If $\epsilon_S <1$, $B_S$ is an ellipsoid of dimension $n-s+2$.
			\begin{enumerate}
			\item[] If $\mathbf{x}_S = \mathbf{a}_S$, then  $\beta_S = 0$; else compute $\beta_S$.\\
			Construct $\mathbf{v}_S$, $\mathbf{u}_S$ and the  ellipse $Y_S$ using  (\ref{use}) and (\ref{elp2d})\\
			Construct the search path $X_S$, using  (\ref {xpthelp}).\\
			Find the step size $\alpha_k$ so that  $\mathbf{x}(\alpha_k) \in X_S \cap B_{i_1,k}$, for $\mathbf{p}_k \in P \setminus S$.\\
			Find the step size $\beta_k$ so that $\mathbf{h}_k\mathbf{y}(\beta_k) =\mathbf{h}_k\mathbf{d}_S$, where $\mathbf{h}_k$ is given by (\ref{hrr}), \\
			and $\beta_k$ is determined by (\ref{AAAv}) and (\ref{elsol}).  \\
			Then $\alpha_k = \beta_k - \beta_S$.
			\end{enumerate}
		\item[3 b (3):] If $\epsilon_S =1$, $B_S$ is a paraboloid of dimension $n-s+2$.
			\begin{enumerate}
			\item[] If $\mathbf{x}_S = \mathbf{a}_S$, then  $\beta_S = 0$; else compute $\beta_S$.\\
			Construct $\mathbf{v}_S$, $\mathbf{u}_S$ and the parabola $Y_S$ using  (\ref{para2d}) and (\ref{parau}). \\
			Construct the search path $X_S$  using (\ref{xpthpara} ).\\
			Find the step size $\alpha_k$ so that  $\mathbf{x}(\alpha_k) \in X_S \cap B_{i_1,k}$, for $\mathbf{p}_k \in P \setminus S$.\\
			Find the step size $\beta_k$ so that $\mathbf{h}_k\mathbf{y}(\beta_k) =\mathbf{h}_k\mathbf{d}_S$, where $\mathbf{h}_k$ is given by (\ref{hrr}), \\
			and $\beta_k$ is determined by (\ref{AAAP}).  \\
			Then $\alpha_k = \beta_k - \beta_S$.
			\end{enumerate}
		 \end{enumerate}
	\item[] Determine $\hat{\alpha}_S$ using (\ref{alphats}).
	\item[] Choose $\alpha_S = \min \{ \hat{\alpha}_S, \alpha_{k} : \mathbf{p}_{k} \in P \setminus S \}$, 
		\begin{enumerate}
		\item[] If $\alpha_S = \hat{\alpha}_S$, go to Step 4 with active set $S$ and $[\mathbf{x}_S, z_S ] = [\mathbf{x}(\alpha_S), z(\alpha_S) ]$.
		\item[] Else, choose $\mathbf{p}_e \in P \setminus S$, so that $\alpha_S = \alpha_e$.  
		\item[] Reset $S: = S \cup \{ \mathbf{p}_e \}$, and  $[\mathbf{x}_S, z_S ] = [\mathbf{x}(\alpha_S), z(\alpha_S) ]$. Go to Step 4.
		\end{enumerate}
	\end{enumerate}
\item Update $S$: Given $[\mathbf{x}_S, z_S]$ and $S$, 
	\begin{enumerate}
	\item[] Compute solution of (\ref{kktch2}) and (\ref{kktch3}) with $S$ and  $\mathbf{x}_S$.
	\item[] $\hspace{.2in}$ If there is a solution with $\lambda_{i_j} \geq 0$, $[{\bf x}_S,  z_S]$ is optimal
	\item[] $\hspace{.2in}$ If there is a solution but $\lambda_{i_l} < 0$ for some $\mathbf{p}_{i_l} \in S$, then $\mathbf{p}_{i_l}$ leaves.  
	\item[] \hspace{.4in} Return to Search Phase with $S:= S  \setminus \{\mathbf{p}_{i_l} \}$ and $[{\bf x}_S,  z_S]$.
	\item[] $\hspace{.2in}$ If no solution,  return to Search Phase with $S: = S$,  and $[{\bf x}_S, z_S]$.
	\end{enumerate}
\end{enumerate}


\section{Dual Algorithm}

At each iteration of the dual algorithm there is an active set $S$ and a corresponding ball $[ \mathbf{x}_S, z_S ]$  that satisfy 
expressions (\ref{kktch2}), (\ref{kktch3}), (\ref{kktch4}), and (\ref{kktch5}) of the KKT conditions, but do not satisfy 
expressions (3) of the KKT condition.  That is, $S$ and $[ \mathbf{x}_S, z_S ]$
are dual feasible but not primal feasible.

Since  $[\mathbf{x}_S, z_S ]$ is a feasible covering of the balls $[ \mathbf{p}_i, r_i ]$ for 
$\mathbf{p}_i \in S$, then $[\mathbf{x}_S, z_S ]$ is an optimal solution to 
 the problem $M(S)$, a relaxation of $M(P)$.  Therefore, 
$z_S$ is a lower bound on the optimal objective function value for $M(P)$.
Assuming non-degeneracy,  the  radius $z_S$ will be shown to increase at each iteration of the dual algorithm.

 The dual algorithm is initialized by choosing  $S = \{ \mathbf{p}_{j}, \mathbf{p}_{k} \}$, for any two points 
$\mathbf{p}_{j}, \mathbf{p}_{k} \in P$, with  $r_{j} \geq r_{k}$.
The initial solution is given by
\begin{equation}
\mathbf{x}_S = \frac{(\mathbf{p}_{j} + \mathbf{p}_{k})}{2} + \frac{(r_{j}-r_{k})}{2} \frac{(\mathbf{p}_{j} - \mathbf{p}_{k})}{\| \mathbf{p}_{j} - \mathbf{p}_{k} \|}, \text{  and   }
z_S =\|\mathbf{x}_S -\mathbf{p}_{j}\!\| + r_{j} =\| \mathbf{x}_S-\mathbf{p}_{k}\!\| + r_{k}.  \label{dinit}
\end{equation}
Observe that $S$ is an active set for the ball  $[\mathbf{x}_S, z_S]$,   $\mathbf{x}_S \in$ ri$(S)$, and  $\mathbf{x}_S \in$ conv$(S)$,
so that  $[\mathbf{x}_S, z_S ]$ is dual feasible to $M(P)$.
 If $r_{j} = r_{k}$, then  $\mathbf{x}_S$ is the center point between $\mathbf{p}_j$ and $\mathbf{p}_k$, and 
 $B_{j,k}$ is the hyperplane $\{ \mathbf{x}: (\mathbf{p}_{j} - \mathbf{p}_{k}) \mathbf{x} = (\mathbf{p}_{j} - \mathbf{p}_{k})\mathbf{x}_S \}$.
 If $r_{j} > r_{k}$, then $\mathbf{x}_S$ is the vertex $\mathbf{a}_{j,k}$ of the sheet $B_{j,k}$ of the hyperboloid $H\!B_{j,k}$.

\noindent \begin{bfseries} Dual Update Phase \end{bfseries}

A dual feasible ball $[\mathbf{x}_S, z_S ]$ is optimal to $M(P)$ if it is primal feasible to $M(P)$, that is, if 
 $z_S \geq \| \mathbf{p}_i-\mathbf{x}_S \| + r_i$,
for all $\mathbf{p}_i \in P$.  
If not optimal, then an entering point $\mathbf{p}_{e} \in P \setminus S$ is chosen  
such that  $z_S < \| \mathbf{p}_e-\mathbf{x}_S \| + r_e$.

If the points in $S \cup \{ \mathbf{p}_e \}$ are affinely independent, then  $|S \cup \{\mathbf{p}_e \}| \leq n+1$, 
and $\mathbf{x}_S \in$ conv$(S \cup \{ \mathbf{p}_e\})$, so that $[\mathbf{x}_S, z_S]$ remains dual feasible with respect to $S \cup \{\mathbf{p}_e \}$.
The algorithm enters the Search Phase with the set  $S$,  the ball $[\mathbf{x}_S, z_S ]$, and the entering point $\mathbf{p}_e$.
The Search Phase searches for a new solution that is active for the set $S \cup  \{\mathbf{p}_e \}$.
 
If the points in $S \cup \{ \mathbf{p}_e \}$ are affinely dependent, then a point $\mathbf{p}_{l} \in S$ is chosen 
to leave $S$.
Since $[\mathbf{x}_S, z_S ]$ is dual feasible to $M(P)$, there exists 
a non-negative solution $(\pi_{1},\ldots,\pi_{s})$  to the 
linear system (\ref{denom1}), (\ref{denom2}),   over the set $S$.
Since the points in $S \cup \{\mathbf{p}_{e}\}$ are affinely dependent, the linear system (\ref{linsys1}) and (\ref{linsys2})
has a solution with  $\lambda_{j}<0$, for some $\mathbf{p}_j \in S$.
\begin{alignat}{3}
&\sum_{\mathbf{p}_j \in S} \lambda_{j}= -1 & & \label{linsys1} \\
&\sum_{\mathbf{p}_j \in S}  (\mathbf{x}_S - \mathbf{p}_{j})\lambda_{j}  =-( \mathbf{x}_S - \mathbf{p}_e).  & & \label{linsys2}   
\end{alignat} 

\noindent A leaving point $\mathbf{p}_l$ is chosen by the `minimum ratio rule,' that is, 
\begin{equation}
\frac{\pi_{l}}{-\lambda_{l}} = \underset{\mathbf{p}_j \in S}{\min} \left\{ \frac{\pi_{j}}{-\lambda_{j}}:
\lambda_{j}<0 \right\}.
\label{minratio1}\end{equation}
The next Theorem shows that  $\mathbf{x}_S$ remains dual feasible for the set $S \cup \{\mathbf{p}_{e}\} \setminus\{\mathbf{p}_{l}\}$.

\begin{theorem}
Suppose that  $[\mathbf{x}_S, z_S]$, and the corresponding active set $S$, is a dual feasible, but not an optimal, solution to $M(P)$.
Suppose that $\mathbf{p}_{e}$ is the point chosen to enter the set $S$ and that 
the set $S \cup \{\mathbf{p}_{e} \}$ is affinely dependent.
If  the leaving point $\mathbf{p}_{l}$ is chosen by  (\ref{minratio1}),
then the points in $S \cup \{\mathbf{p}_{e}\} \setminus\{\mathbf{p}_{l}\}$
are affinely independent, and $\mathbf{x}_S \in $ conv$(S \cup \{\mathbf{p}_{e}\} \setminus\{\mathbf{p}_{l}\})$.
\end{theorem}

\noindent \textbf{Proof:}
A non-negative solution to the linear system  (\ref{denom1}), (\ref{denom2}) over the set $S$
corresponds to a basic feasible solution if and only if the points $\mathbf{p}_j$ in $S$ are affinely independent.
The entering point
$\mathbf{p}_{e}$ corresponds to the right hand side  column of the linear system (\ref{linsys1}), (\ref{linsys2}). 
The minimum ratio rule guarantees that if the column corresponding to $\mathbf{p}_l$ is replaced by the column corresponding to $\mathbf{p}_e$,
then the columns of the new basis, which correspond to points in  $S \cup \{\mathbf{p}_{e}\} \setminus\{\mathbf{p}_{l}\}$, are
linearly independent, and that the new basis yields a non-negative solution for the system  (\ref{denom1}), (\ref{denom2}),
with $S$ replaced by $S \cup \{\mathbf{p}_{e}\} \setminus\{\mathbf{p}_{l}\}$. 
Thus, the points in $S \cup \{\mathbf{p}_{e}\} \setminus\{\mathbf{p}_{l}\}$
are affinely independent, $|S \cup \{\mathbf{p}_{e}\} \setminus\{\mathbf{p}_{l}\}| \leq n+1$, 
and $\mathbf{x}_S \in $ conv$(S \cup \{\mathbf{p}_{e}\} \setminus\{\mathbf{p}_{l}\})$.  \hfill $\Box$ 

If the point $\mathbf{p}_l$ leaves $S$, then the set  $S$ is reset to $S : = S \setminus \{ \mathbf{p}_{l} \}$.
Thus  $|S \cup \{\mathbf{p}_e \}| \leq n+1$, and $\mathbf{x}_S \in$ conv$(S \cup \{ \mathbf{p}_e\})$, 
so that $[\mathbf{x}_S, z_S]$ remains dual feasible with respect to $S \cup \{\mathbf{p}_e \}$.
The updated set $S$ remains active for $[\mathbf{x}_S, z_S ]$.
The algorithm enters the Search Phase with the set  $S$,  the ball $[\mathbf{x}_S, z_S ]$, 
and the entering point $\mathbf{p}_e$.
The Search Phase searches for a new solution that is active for the set $S \cup  \{\mathbf{p}_e \}$.

\noindent \begin{bfseries} Dual Search Phase \end{bfseries}

Given a dual feasible ball $[ \mathbf{x}_S, z_S ]$, its active set $S$, and an entering point  $\mathbf{p}_e$ as determined in the Update Phase, 
a search path $X_{S} = \{ \mathbf{x}(\alpha): \alpha \geq 0 \}$  is constructed so that $X_S \subset B_S$,
$\mathbf{x}_S = \mathbf{x}(0)$, and  $S$ is an active set for the ball  $[\mathbf{x}(\alpha), z(\alpha)]$, for $\alpha \geq 0$, where $z(\alpha) = \| \mathbf{x}(\alpha) -\mathbf{p}_i \| + r_i$ for 
$\mathbf{p}_i \in S$.  
Also, the complementary slackness conditions \eqref{kktch5}   are satisfied by $S$ and  $[\mathbf{x}(\alpha), z(\alpha)]$, for $\alpha \geq 0$.

A step size $\alpha'$ is determined so that $[\mathbf{x}(\alpha'), z(\alpha')]$ is active for $S \cup \{\mathbf{p}_e \}$.
If $\mathbf{x}(\alpha') \in$ conv$(S \cup \{\mathbf{p}_e\} )$, then $[\mathbf{x}(\alpha'), z(\alpha')]$ is also dual feasible for $S \cup \{\mathbf{p}_e \}$.
In this case, reset $S: = S \cup \{\mathbf{p}_e \}$, and set $[\mathbf{x}_S, z_S]:= [\mathbf{x}(\alpha'), z(\alpha')]$.
The Update Phase is entered with $S$ and $[\mathbf{x}_S, z_S]$.

 If $\mathbf{x}(\alpha')  \notin$ conv$(S \cup \{\mathbf{p}_e\} )$, some point $\mathbf{p}_{i_l}$ is chosen to leave $S$, and a step size $\alpha^*$ is determined 
 so that the ball  $[\mathbf{x}(\alpha^*), z(\alpha^*)]$ is dual feasible 
with respect to $S \cup \{\mathbf{p}_e\} \setminus \{\mathbf{p}_{i_l}\} $.
 The set $S$ is reset to $S:= S \setminus \{ \mathbf{p}_{i_l} \}$,
and the ball is reset to  $[\mathbf{x}_S, z_S] :=   [\mathbf{x}(\alpha^*), z(\alpha^*)]$. 
The Search Phase is re-entered with the active set $S$, its corresponding dual feasible ball $[\mathbf{x}_S, z_S]$, and the entering point $\mathbf{p}_e$.

If all the points in $S$ have equal radii, the search path will be a ray, but if some points in $S$ have unequal radii, the search path will be a two dimensional conic section.


 \noindent 
\begin{bfseries}Case 1: All points in $S$ have equal radii \end{bfseries}

The points in $S$ are denoted by $S = \{ \mathbf{p}_{i_1}, \ldots,  \mathbf{p}_{i_{s}} \}$, where $r_{i_j} = r_{i_1}$ for $\mathbf{p}_{i_j} \in S$.
In this case, the search path is the ray:
\begin{alignat}{1}
&X_S = \{ \mathbf{x}(\alpha) = \mathbf{x}_{S} + \alpha \mathbf{d}_{S}, \alpha \geq 0 \}, \text{\;\;\;where\;\;\;}
\mathbf{d}_{S} = (\mathbf{p}_{e} - \mathbf{p}_{i_1}) - \text{Proj}_{\text{sub}(S)} (\mathbf{p}_{e} - \mathbf{p}_{i_1}), \label{dird}
\end{alignat} 
and sub$(S) = $ span$[ (\mathbf{p}_{i_1} - \mathbf{p}_{i_2}),  \ldots, (\mathbf{p}_{i_1} - \mathbf{p}_{i_{s}}) ]$.  
Thus, $\mathbf{d}_{S}$ is the component of $(\mathbf{p}_{e} - \mathbf{p}_{i_1})$  orthogonal to the projection of 
$(\mathbf{p}_{e} - \mathbf{p}_{i_1})$ onto sub$(S)$.
 If $s = 1$, with $S = \{ \mathbf{p}_{i_1} \}$, $\mathbf{d}_{S} = \mathbf{p}_{e} - \mathbf{p}_{i_1}$.  

\begin{property}
Suppose  the active set $S = \{ \mathbf{p}_{i_1}, \ldots, \mathbf{p}_{i_{s}} \}$  
corresponds to the dual feasible ball $[\mathbf{x}_{S}, z_{S} ]$, and the radius of each point $\mathbf{p}_{i_j} \in S$ is $r_{i_1}$.
Let $X_{S} = \{ \mathbf{x}(\alpha) =  \mathbf{x}_{S} + \alpha \mathbf{d}_{S}, \alpha \geq 0 \}$ where
$\mathbf{d}_{S} = (\mathbf{p}_{e} - \mathbf{p}_{i_1}) -  \text{Proj}_{\text{sub}(S)} (\mathbf{p}_{e} - \mathbf{p}_{i_1})$. 
Let $z(\alpha) = \|\mathbf{x}_{S} - \mathbf{p}_{i_1} \| + r_{i_1}$.
Then $S$ is an active set corresponding to the ball $[\mathbf{x}(\alpha), z(\alpha) ]$, for $\alpha \geq 0$, 
 Furthermore,  $z(\alpha) $ is increasing for $ 0 \leq \alpha \leq \hat{\alpha}_S$, 
 where $\hat{\alpha}_S = \frac{(  \mathbf{p}_{e} - \mathbf{x}_{S} )  \mathbf{d}_{S}}{\| \mathbf{d}_{S} \|^2} > 0$.
 \end{property}
 \noindent \textbf{Proof:}  For each pair of points $\mathbf{p}_{i_j}$ and $\mathbf{p}_{i_k}$ in $S$, with $r_{i_j} = r_{i_k}$, the bisector $B_{i_j,i_k}$ may be written as $B_{i_j,i_k} = \{ \mathbf{x}: (\mathbf{p}_{i_j} - \mathbf{p}_{i_k})\mathbf{x} = (\mathbf{p}_{i_j} - \mathbf{p}_{i_k})\mathbf{x}_{S} \}$ since $\mathbf{x}_{S} \in B_{i_j,i_k}$. Then
 $(\mathbf{p}_{i_j} - \mathbf{p}_{i_k})\mathbf{x}(\alpha) =  (\mathbf{p}_{i_j} - \mathbf{p}_{i_k})\mathbf{x}_{S} + \alpha  (\mathbf{p}_{i_j} - \mathbf{p}_{i_k})\mathbf{d}_{S}
 = (\mathbf{p}_{i_j} - \mathbf{p}_{i_k})\mathbf{x}_{S}$, since $\mathbf{d}_{S}$ is orthogonal to $(\mathbf{p}_{i_j} - \mathbf{p}_{i_k})$ by construction.
 Thus for all $\mathbf{p}_{i_j}, \mathbf{p}_{i_k} \in S$, and  for $\alpha \geq 0$,
 $\mathbf{x}(\alpha) \in B_{i_j,i_k}$,  so that $S$ is an active set for $[\mathbf{x}(\alpha), z(\alpha) ]$.
 
 For each $ \mathbf{p}_{i_j} \in S$, $(  \mathbf{p}_{e} -  \mathbf{p}_{i_j}) \mathbf{d}_{S} > 0$ by the construction of $\mathbf{d}_{S}$.
 Property 3.3 implies that $\mathbf{x}(\alpha) \neq \mathbf{p}_{i_j}$ 
for all $\alpha$, so that $z(\alpha)$ is differentiable with respect  to $\alpha$.
 Thus  for each $\mathbf{p}_{i_j} \in S$, $\hat{\alpha}_S = \frac{(  \mathbf{p}_{e} -  \mathbf{x}_{S}) \mathbf{d}_{S}}{\| \mathbf{d}_{S} \|^2} > 0$ 
 is the  unique maximum of  $z(\alpha)$, so $z(\alpha)$ is increasing along $X_{S}$ for $0 \leq \alpha \leq \hat{\alpha}_S$.   \hfill $\Box$
 
 \begin{corollary}  Given  $S = \{ \mathbf{p}_{i_1} \}$, let  $X_{S} = \{ \mathbf{x}(\alpha) = \mathbf{x}_{S} + \alpha \mathbf{d}_{S}, \alpha \geq 0 \}$ 
 where $\mathbf{d}_{S} = (\mathbf{p}_{i_1} - \mathbf{x}_{S})$.  
 Then  $z(\alpha) = \| \mathbf{x}(\alpha) - \mathbf{p}_{i_1} \| + r_{i_1}$ is increasing for $0 \leq \alpha \leq \hat{\alpha}_S$,  where $\hat{\alpha}_S = \frac{(  \mathbf{p}_{i_1} - \mathbf{x}_{S} )  \mathbf{d}_{S}}{\| \mathbf{d}_{S} \|^2} > 0$.
 \end{corollary}

\noindent \textbf{The Step Size:} The step size $\alpha' \geq 0$ is determined, if it exists, so that $X_{S}$ intersects the bisector $B_{i_1,e}$ at $\mathbf{x}(\alpha')$.
If $\mathbf{x}(\alpha') \in X_{S} \cap B_{i_1,e}$, then  $\mathbf{x}(\alpha') \in X_{S} \cap B_{i_j,e}$ for each $\mathbf{p}_{i_j} \in S$.  
That is, $X_{S}$ simultaneously intersects the bisectors $B_{i_j,e}$ at 
$\mathbf{x}(\alpha')$, for all $\mathbf{p}_{i_j} \in S$.  Thus it sufficies to consider the intersection of $X_{S}$ with only $B_{i_1,e}$.

At the point $\mathbf{x}(\alpha') \in X_{S} \cap B_{i_1,e}$,  $ z(\alpha') = \| \mathbf{x}(\alpha')  - \mathbf{p}_e \| +r_e = \| \mathbf{x}(\alpha')  - \mathbf{p}_{i_1} \| + r_{i_1}$,
that is, $S \cup \{ \mathbf{p}_e \}$  is active at $[ \mathbf{x}(\alpha'), z(\alpha') ]$.
Geometrically, the search moves the center $\mathbf{x}(\alpha)$ of the ball $[\mathbf{x}(\alpha), z(\alpha)]$ along the path $X_{S}$, while the radius $z(\alpha)$ increases.

There are two sub-cases to consider for computing $\alpha'$ depending on whether the radius $r_{i_1}$ equals the radius $r_e$.\\
Sub-case 1: $r_{i_1} = r_e$, so that  $B_{i_1,e}$ is a hyperplane. The intersection of $X_S$ and $B_{i_1,e}$ is determined by 
substituting the expression for the ray $\mathbf{x}(\alpha)$ into the equation for the hyperplane $B_{i_1,e}$ which yields
 $(\mathbf{p}_{i_1} - \mathbf{p}_e)\mathbf{x}(\alpha) = (\mathbf{p}_{i_1} - \mathbf{p}_e)(\mathbf{p}_{i_1} + \mathbf{p}_e)/2$, and
 \begin{alignat}{2}
\alpha =  \frac{(\mathbf{p}_{i_1} - \mathbf{p}_{e}) (\mathbf{p}_{i_1} + \mathbf{p}_{e})/2 - (\mathbf{p}_{i_1} - \mathbf{p}_{e})  \mathbf{x}_S}
 {(\mathbf{p}_{i_1} - \mathbf{p}_{e})  \mathbf{d}_S}. \label{PEQR}
 \end{alignat}
 If $(\mathbf{p}_{i_1} - \mathbf{p}_{e}) \mathbf{d}_S = 0$, or if $\alpha < 0$, set $\alpha' := \infty$.  Otherwise,
 set $\alpha' = \alpha$.

\noindent Sub-case 2: $r_{i_1} \neq r_e$, so that  $B_{i_1,e}$ is a hyperboloid.  The intersection of $X_S$ and $B_{i_1,e}$ is determined by substituting the expression for 
the ray $\mathbf{x}(\alpha)$ into 
 the quadratic form  \eqref{Q} for the hyperboloid $B_{i_1,e}$, with center  $ \mathbf{c}_{i_1,e}$, axis vector $ \mathbf{v}_{i_1,e}$, 
 eccentricity $\epsilon_{i_1,e}$, and parameters $ a_{i_1,e}$ and $c_{i_1,e}$.  
This  gives the quadratic equation 
\begin{alignat}{3}
&A_{e}\alpha^2 + B_{e}\alpha + C_{e} = 0, \label{QUADe} 
\end{alignat} 
 where $A_{e} = (\mathbf{d}_S)^2 - \epsilon_{i_1,e}^2(\mathbf{d}_S \mathbf{v}_{i_1,e})^2,$   
$B_{e} = 2(\mathbf{x}_S - \mathbf{c}_{i_1,e})  \mathbf{d}_S - 2\epsilon_{i_1,e}^2 [(\mathbf{x}_S - \mathbf{c}_{i_1,e})  \mathbf{v}_{i_1,e} ] 
[\mathbf{d}_S  \mathbf{v}_{i_1,e}] ,$
and $C_{e} = (\mathbf{x}_S - \mathbf{c}_{i_1,e})^2 - \epsilon_{i_1,e}^2[(\mathbf{x}_S - \mathbf{c}_{i_1,e})  \mathbf{v}_{i_1,e}]^2
 + a_{i_1,e}^2 - c_{i_1,e}^2$.
If there are no real solutions, $\alpha' = \infty$.
If there is one real solution $\alpha$,  then $\alpha' = \alpha$.
Otherwise, let $\alpha_1$ and $\alpha_2$ be the real solutions.
If $A_e < 0$ and $r_{i_1} < r_e$, then $\alpha' = \max \{ \alpha_1, \alpha_2 \}$, or if $r_{i_1} > r_e$, then $\alpha' = \min \{ \alpha_1, \alpha_2 \}$.
If $A_e > 0$ and $r_{i_1} < r_e$, then $\alpha' = \min \{ \alpha_1, \alpha_2 \}$, or if $r_{i_1} > r_e$, then $\alpha' = \infty$.

If $\mathbf{x}(\alpha') \in $ conv$(S \cup \{ \mathbf{p}_e \})$, then $\mathbf{x}(\alpha')$ is dual feasible. 
To determine if $\mathbf{x}(\alpha') \in$ conv$(S \cup \{ \mathbf{p}_e \})$, 
substitute  $\mathbf{x}(\alpha) = \mathbf{x}_{S} + \alpha \mathbf{d}_{S}$ 
for $\mathbf{x}_S$ in equation (\ref{denom2}),  expand and simplify to obtain the linear system: 
\begin{alignat}{2}
&T\boldsymbol{\pi}(\alpha) =  \boldsymbol{\varepsilon}_1 - \alpha\hat{\mathbf{d}}_{S}  \label{Teqr}
\end{alignat}
\begin{alignat}{3}
\text{where\;\;\;}T = \left[
\begin{array}{*{1}{cccccc}}
1&  \ldots & 1 &1 \\
\mathbf{x}_{S} - \mathbf{p}_{i_1}  & \ldots & \mathbf{x}_{S} - \mathbf{p}_{i_{s}}  & \mathbf{x}_{S} - \mathbf{p}_{e} \\
\end{array}
\right]\!, \;\; \boldsymbol{\varepsilon}_1 =  \left[
\begin{array}{*{1}{c}}
1\\
\mathbf{0} \\
\end{array}
\right]\!, \;\;  \hat{\mathbf{d}}_{S}  =  \left[
\begin{array}{*{1}{c}}
0\\
\mathbf{d}_{S} \\
\end{array}
\right]\!, \notag
\end{alignat}
and each variable is a function of $\alpha$, that is, $\boldsymbol{\pi}(\alpha) = [\pi_{i_1}(\alpha), \ldots, \pi_{i_s}(\alpha), \pi_{e}(\alpha) ]$.
The matrix $T$ is $n+1 \times s+1$ with rank $s+1$.  
The linear system (\ref{Teqr}) has a solution  since $\mathbf{x}(\alpha) \in $ aff$(S \cup \mathbf{p}_e)$ for $\alpha \geq 0$.

If $\boldsymbol{\pi}(\alpha') \geq 0$, then $\mathbf{x}(\alpha') \in $ conv$(S \cup \{ \mathbf{p}_e \})$.
In this case, reset $S:= S \cup \{\mathbf{p}_e\}$, and  $[\mathbf{x}_S, z_S ] :=  [\mathbf{x}(\alpha'),z(\alpha')]$.
Then the ball $[\mathbf{x}_S, z_S ]$ and its active set $S$ are checked for optimality by the Update Phase.

Otherwise, some component of $\boldsymbol{\pi}(\alpha')$ is negative, and $\mathbf{x}(\alpha') \notin$ conv$(S \cup \{ \mathbf{p}_e \})$. 
In this case, a new step size $\alpha^*$  is determined so that   $0 \leq \alpha^* \leq \alpha'$ and $\mathbf{x}(\alpha^*)  \in$ conv$(S \cup \{ \mathbf{p}_e \})$. 
Since the $s+1$ points in $S \cup \{ \mathbf{p}_e\}$ are affinely independent,
conv$(S  \cup \{ \mathbf{p}_e\})$ is a simplex with $s+1$ vertices and $s+1$ facets.  The $s+1$ vertices are the points in $S \cup \{ \mathbf{p}_e\}$.  
The $s+1$  facets  are denoted by $F_{i_j} =$ conv$(S \cup \{ \mathbf{p}_e\} \setminus \{\mathbf{p}_{i_j} \})$, 
 for each $\mathbf{p}_{i_j} \in S$, and by $F_e = $ conv$(S)$, for $\mathbf{p}_e$.
Each  point $\mathbf{p}_{i_j} \in S$ corresponds to the 
component $\pi_{i_j}(\alpha)$ of $\boldsymbol{\pi}(\alpha)$ and  to the facet $F_{i_j}$.
The  point $\mathbf{p}_e$ corresponds to the component $\pi_e(\alpha)$ of $\boldsymbol{\pi}(\alpha)$ and  to the facet $F_e$. 

Since $\mathbf{x}(0) \in$ conv$(S \cup \{ \mathbf{p}_e \})$,  and $\mathbf{x}(\alpha') \notin$ conv$(S \cup \{ \mathbf{p}_e \})$, 
$X_{S}$ must intersect some facet of conv$(S \cup \{ \mathbf{p}_e \})$
between $\mathbf{x}(0)$ and $\mathbf{x}(\alpha')$.
 For each $\mathbf{p}_{i_j} \in S$, the  step-size $\alpha_{i_j} \geq 0$ is computed, if it exists, 
 so that $\mathbf{x}(\alpha_{i_j}) \in X_S \cap F_{i_j}$,
 which is equivalent to  $\pi_{i_j}(\alpha_{i_j}) = 0$ in the solution to 
$ T\boldsymbol{\pi}(\alpha_{i_j}) =  \boldsymbol{\varepsilon}_1 - \alpha_{i_j}\hat{\mathbf{d}}_{S}$.
 
 To determine $\alpha_{i_j}$ for each $\mathbf{p}_{i_j} \in S$, so that $\pi_{i_j}(\alpha_{i_j}) = 0$,
 in the solution $\boldsymbol{\pi}(\alpha_{i_j})$ to 
 $ T\boldsymbol{\pi}(\alpha_{i_j}) =  \boldsymbol{\varepsilon}_1 - \alpha_{i_j}\hat{\mathbf{d}}_{S}$, 
solve $T \boldsymbol{\gamma}  = \boldsymbol{\varepsilon}_1$ for $\boldsymbol{\gamma}$, and 
 solve $T \boldsymbol{\delta} = \hat{\mathbf{d}}_S$ for $\boldsymbol{\delta}$.
 Then, $T \boldsymbol{\pi}(\alpha_{i_j})  = \boldsymbol{\varepsilon}_1 - \alpha_{i_j}\hat{\mathbf{d}}_{S} = T(\boldsymbol{\gamma} - \alpha_{i_j} \boldsymbol{\delta})$.
If $ \gamma_{i_j} - \alpha_{i_j} \delta_{i_j} = 0$, then $\pi_{i_j}(\alpha_{i_j}) = 0$.
For each $\mathbf{p}_{i_j} \in S$, set $\alpha_{i_j} = \gamma_{i_j}/\delta_{i_j}$ if $\delta_{i_j} > 0$, and  set $\alpha_{i_j} = \infty$ if $\delta_{i_j} \leq 0$.  
 For the point $\mathbf{p}_e$, $\mathbf{x}(\alpha_e) \in$ conv$(S)$ implies $\alpha_e \leq 0$.  
 That is,  if $ \mathbf{x}_S \in$ conv$(S)$, then $\mathbf{x}(0) = \mathbf{x}_S$, so that $\alpha_e = 0$.
 If  $\mathbf{x}_S \notin $ conv$(S)$, then $\mathbf{x}(\alpha_e) \in$ aff$(S)$, where $\alpha_e = \frac{(\mathbf{p}_{i_1} - \mathbf{x}_S)\mathbf{d}_S}{\mathbf{d}_S^2} < 0$.
 
 The intersection of $X_S$ and the facet of conv$(S \cup \{ \mathbf{p}_e \})$ first encountered along $X_S$ occurs 
at the step size $\alpha^*$, where
\begin{equation}
 \alpha^*  = \min_{\mathbf{p}_{i_j} \in S } \{ \alpha_{i_j}   \}. \label{alstr}
\end{equation}
Cavaleiro and Alizadeh \cite{Cavaleiro}  present an equivalent procedure for finding $\alpha^*$.

The leaving point $\mathbf{p}_{i_l}  \in S$ is chosen so that $\alpha^* = \alpha_{i_l}$.
Then the solution to the linear system (\ref{Teqr}) with $\mathbf{x}(\alpha^*)$ substituted for $\mathbf{x}_S$, yields
$\pi_{i_l} (\alpha^*)= 0$, and $\pi_{i_j}(\alpha^*) \geq 0$, for $\mathbf{p}_{i_j} \in S \cup \{ \mathbf{p}_e \} \setminus \{ \mathbf{p}_{i_l} \}$, 
so that $\mathbf{x}(\alpha^*) \in F_{i_l}$.
This construction shows that  $S  \setminus \{ \mathbf{p}_{i_l} \}$ is an active set for $[\mathbf{x}(\alpha^*), z(\alpha^*)]$
and that  $[\mathbf{x}(\alpha^*), z(\alpha^*)]$ is dual feasible with respect to $S \cup \{ \mathbf{p}_e \} \setminus \{ \mathbf{p}_{i_l} \}$.
However, the constraint corresponding to $\mathbf{p}_e$ is not active at $[\mathbf{x}(\alpha^*), z(\alpha^*)]$, and $\mathbf{p}_e$ is not added to $S$.

 Delete the point $\mathbf{p}_{i_l}$ from $S$.  Reset $S: =  S \setminus \{ \mathbf{p}_{i_l} \}$, 
 and  $[\mathbf{x}_S, z_S] := [ \mathbf{x}(\alpha^*), z(\alpha^*) ]$.
The set $S$ is active with respect to the ball $[\mathbf{x}_S, z_{S}]$, and $[\mathbf{x}_S, z_{S}]$ is dual feasible with respect to 
$S \cup \{ \mathbf{p}_e \} $.  
The Search Phase is re-entered with the active set $S$, the entering point $\mathbf{p}_e$, and the ball $[\mathbf{x}_S, z_{S}]$.

If there is more than one  point
$\mathbf{p}_{i_j} \in S \cup \{\mathbf{p}_{e}\}$ such that $\alpha^* = \alpha_{i_l}$, then  $X_S$
intersects the corresponding facets $F_{i_j}$ simultaneously, and there is a tie for the leaving point. 
In this case choose any point $\mathbf{p}_{i_j}$
 such that $\alpha_{i_j} = \alpha^*$, to be deleted from $S$. 
Reset  $S:=S \setminus \{ \mathbf{p}_{i_j} \}$, and $[\mathbf{x}_S,z_S] := [\mathbf{x}(\alpha^*), z(\alpha^*)]$. 
The Search Phase is re-entered with the active set $S$,
the entering point $\mathbf{p}_e$, and the ball $[\mathbf{x}_S,z_S]$.
This case may lead to a degenerate iteration at  the next step
with $\alpha^* = \alpha_{i_l} =0$ for some $\mathbf{p}_{i_l}$.  
Cycling will not occur since at each degenerate iteration one point is deleted from the finite set $S$.  
After a finite number of points are deleted, 
 $S$ is reduced to two points, and the step size $\alpha$ will be positive at the next iteration.

\begin{property}
If all the points in $S$ have equal radii, the radius of the covering ball,  given by 
$z_{S}(\alpha) = \parallel \mathbf{x}(\alpha) - \mathbf{p}_{i_j} \parallel + r_{i_j}$ 
for any $\mathbf{p}_{i_j} \in S \setminus \{ \mathbf{p}_e \}$,
is strictly increasing for $\alpha > 0 $, 
and the distance  $\parallel \mathbf{p}_e - \mathbf{x}(\alpha) \parallel$ is strictly decreasing for $\alpha > 0 $.
\end{property}
\noindent Proof: By construction of $\mathbf{d}_e$, 
  $(\mathbf{p}_{i_j} - \mathbf{x}_S) \mathbf{d}_e= 0$,
  for each $\mathbf{p}_{i_j} \in S\setminus \{\mathbf{p}_{e} \}$.
  and $(\mathbf{p}_{e} - \mathbf{x}_S)  \mathbf{d}_e > 0$.
  The convex function $g_j(\alpha) =  \parallel \mathbf{p}_{i_j} - \mathbf{x}(\alpha) \parallel^2$ 
has a unique minimum point at  $\alpha = \frac{(\mathbf{p}_{i_j} - \mathbf{x}_S) \mathbf{d}_e}{
 \parallel \mathbf{d} \parallel^2} = 0$, 
so that $g_j(\alpha)$ is increasing for $\alpha > 0$.   Thus $z_{S}(\alpha)$ is increasing for $\alpha > 0$.
The convex function  $f(\alpha) = \parallel \mathbf{p}_{e} - \mathbf{x}(\alpha) \parallel^2$
has a unique minimum at 
 $\alpha = \frac{(\mathbf{p}_{e} - \mathbf{x}_S) \mathbf{d}_S}{\parallel \mathbf{d} \parallel^2} > 0$ ,
 so that $f(\alpha)$ is decreasing for $0 < \alpha \leq  \frac{(\mathbf{p}_{e} - \mathbf{x}_S) \mathbf{d}}{\parallel \mathbf{d} \parallel^2}$.
\hfill $\Box$

Thus, for each iteration with $\alpha > 0$, $z_{S}(\alpha) > z_S$.
Since  $\mathbf{x}_S \in$ ri$(S)$, a positive step size must occur at least at the first iteration.

 \noindent  \begin{bfseries}Case 2:  At least two points in $S$ have unequal radii \end{bfseries}
 
\noindent \textbf{The Search Path:}  The points in $S$ are ordered by non-increasing radii and denoted by 
$S = \{ \mathbf{p}_{i_1},  \ldots,\mathbf{p}_{i_{s}} \}$, with $r_{i_1}  \geq \ldots \geq r_{i_{s}}$.
By the assumption of unequal radii, $r_{i_1} > r_{i_{s}}$, so that $B_{i_1,i_{s}}$ is one sheet of a hyperboloid.  
The intersection of bisectors $B_S  = \cap_{\{ \mathbf{p}_{i_j},\mathbf{p}_{i_k} \} \subset S} B_{i_j,i_k}$
is shown, by Properties A.7 and A.8, to be the intersection of  the sheet $B_{i_1,i_s}$ with $s-2$ hyperplanes resulting in a conic section of dimension $n-s+2$.
The vectors and parameters of $B_S$, namely
the center $\mathbf{c}_{S}$, axis vector $\mathbf{v}_{S}$,   vertex $\mathbf{a}_{S}$,   eccentricity  $e_{S}$,  and parameters  $a_{S}$ and  $b_{S}$
are computed using expressions  (\ref{HT12})  - (\ref{pnc}) in the Appendix.
If $\epsilon_S > 1$, then $B_S$ is one sheet of a hyperboloid. If $\epsilon_S <1$, then $B_S$ is an ellipsoid.  If $\epsilon_S = 1$, then $B_S$ is a paraboloid.

Property A.9 in the Appendix shows that for any vector $\mathbf{u}_{S}$ orthogonal to 
$\mathbf{v}_{S}$, there is a two dimensional conic section $Y_S \subset B_S$ that has the same vectors and parameters as $B_{S}$. 
Furthermore, $Y_S \subset$ aff$(\mathbf{c}_{S}, \mathbf{v}_{S}, \mathbf{u}_{S})$.
Thus $Y_S$ is either one sheet of a hyperbola, an ellipse, or a parabola if $\epsilon_S <1, \epsilon_S <1,$ or $\epsilon_S = 1$, respectively.  
The search path $X_S$ is chosen to be a sub-path of $Y_S$ determined by a particular choice  of $\mathbf{u}_S$.

\noindent \textbf{Case 2a: $X_S$ is a hyperbola:}

If  $\epsilon_{S} > 1$, then $B_{S} $ is one sheet of a hyperboloid of dimension $n-s+2$, and 
\begin{equation}
 Y_{S}  =  \{ \mathbf{y}(\beta) =  \mathbf{c}_{S} + a_{S} \sec(\beta) \mathbf{v}_{S}  + b_{S} \tan(\beta) \mathbf{u}_{S} : -\pi/2 < \beta < \pi/2 \} \label{ysdualhy}
 \end{equation}
is one sheet of a two-dimensional hyperbola in  aff$(\mathbf{c}_{S}, \mathbf{v}_{S}, \mathbf{u}_{S})$.  
The vector $\mathbf{u}_{S}$ is  chosen to be the normalized component of  $(\mathbf{p}_e - \mathbf{c}_{S})$ 
that is orthogonal to the projection of  $(\mathbf{p}_e - \mathbf{c}_{S})$
 onto sub($S$), that is, 
 \begin{equation}
 \mathbf{u}_{S} = ((\mathbf{p}_{e} - \mathbf{c}_{S}) - \text{Proj}_{\text{sub}(S)} (\mathbf{p}_{e} - \mathbf{c}_{S})) / \| ((\mathbf{p}_{e} - \mathbf{c}_{S}) - \text{Proj}_{\text{sub}(S)} (\mathbf{p}_{e} - \mathbf{c}_{S})) \label{usdhyp}\|. 
 \end{equation}
Property A.8 shows that $Y_{S} $ is a subset of $B_{S}$, and that
 the point $\mathbf{y}(0) = \mathbf{a}_{S}$ is the vertex of the hyperbola $Y_{S} $ and of the hyperboloid $B_{S}$.
 Let $\beta_S$ be the parameter such that $\mathbf{y}(\beta_S) = \mathbf{x}_S$.  
 If $\mathbf{x}_S = \mathbf{a}_S$, then $\beta_S = 0$, but if $\mathbf{x}_S \neq \mathbf{a}_S$, then
 $\beta_S = -\tan^{-1} \{\sqrt{(\mathbf{c}_S - \mathbf{x}_S)^2-a_S^2} / ( \epsilon_Sa_S)\} $.
 If a point was deleted from $S$ in the last iteration, $\mathbf{x}_S \neq \mathbf{a}_S$, else, $\mathbf{x}_S = \mathbf{a}_S$.
 
 The search path $X_S \subset Y_S$ is defined to be a subset of $Y_S$ given by 
 \begin{equation}
 X_S = \{ \mathbf{x}(\alpha) = \mathbf{y}(\alpha + \beta_S), 0 \leq  \alpha \leq -\beta_S \}. \label{xsdualhy}
 \end{equation} 
 
 Observe that $\mathbf{x}(0) = \mathbf{y}(\beta_S) = \mathbf{x}_S$, and for $\alpha \geq 0$, 
 the points $\mathbf{x}(\alpha)$ are on $B_S$ and move along $X_S$ from $\mathbf{x}_S$ toward $\mathbf{p}_e$ and primal feasibility.

  \begin{property}
Suppose the set  $S =  \{ \mathbf{p}_{i_1}, \ldots,  \mathbf{p}_{i_{s}} \}$ is an active set for the  dual feasible ball $[\mathbf{x}_S, z_S]$,
ordered so that $r_{i_1} \geq \ldots \geq r_{i_{s}}$, with $r_{i_1} > r_{i_{s}}$.  
Let $X_S = \{ \mathbf{x}(\alpha) = \mathbf{y}(\alpha + \beta_S),  0 \leq \alpha \leq -\beta_S \}$  
be the search path along a hyperbola $Y_S$ determined by (\ref{ysdualhy}) and (\ref{xsdualhy}) .
Then $S$ is an active set for the ball $[\mathbf{x}(\alpha), z(\alpha) ]$, for $\alpha \geq 0$, where
$z(\alpha) = \| \mathbf{x}(\alpha) - \mathbf{p}_{i_1} \| + r_{i_1} $.
Furthermore, $z(\alpha)$ is increasing for $0 \leq \alpha \leq -\beta$.
 \end{property}
 \noindent \textbf{Proof:} By construction,
 $\mathbf{x}(\alpha) \in B_{i_1,i_j}$ for each $\mathbf{p}_{i_j} \in S$ and $\alpha \geq 0$, so that
 $z(\alpha) = \|\mathbf{x}(\alpha) - \mathbf{p}_{i_1} \| + r_{i_1}  =  \| \mathbf{x}(\alpha) - \mathbf{p}_{i_j} \| + r_{i_j}$.
 Thus, $S$ is an active set for $[\mathbf{x}(\alpha), z(\alpha) ]$, for $\alpha \geq 0$.
 Direct computation shows that $z(\alpha)$ has a minimum value at $\alpha = 0$, and since $z(\alpha)$ is convex with respect to $\alpha$,
 $z(\alpha)$ is increasing for $0 \leq \alpha \leq -\beta_S$. \hfill  $\Box$

  \noindent \begin{bfseries}The Step Size: \end{bfseries}
 The step size $\alpha'$ is determined, if it exists, so that $X_{S}$ intersects the bisector  $B_{i_1,e}$
 at $\mathbf{x}(\alpha')$.
 If $\mathbf{x}(\alpha') \in X_S \cap B_{i_1,e}$, then $\mathbf{x}(\alpha') \in X_S \cap B_{i_j,e}$, for each $\mathbf{p}_{i_j} \in S$, that is, 
$X_S$ simultaneously intersects the bisectors $B_{i_j,e}$ at $\mathbf{x}(\alpha')$, for all $\mathbf{p}_{i_j} \in S$.
Thus, it suffices to determine the intersection of $X_S$ with only $B_{i_1,e}$.

At the point $\mathbf{x}(\alpha') \in X_S \cap B_{i_1,e}$, 
 $z(\alpha') =  \| \mathbf{x}(\alpha') - \mathbf{p}_e \| + r_e =  \| \mathbf{x}(\alpha') - \mathbf{p}_{i_1} \| + r_{i_1}$,
 that is, $S \cup \{ \mathbf{p}_e \}$ is an active set for  $[\mathbf{x}(\alpha'), z(\alpha')]$.
  Geometrically, the search moves the center  $\mathbf{x}(\alpha)$ along the path $X_{S}$ while increasing  
 $z(\alpha)$.

Theorem A.1  shows that $X_S \cap B_{i_1,e} = X_S \cap H_e$, where 
the hyperplane $H_e =  \{ \mathbf{x} : \mathbf{h}_e\mathbf{x} = \mathbf{h}_e\mathbf{d}_e \}$ is constructed as follows.
Order the set of three points $\{\mathbf{p}_{i_1}, \mathbf{p}_{i_{s}}, \mathbf{p}_{e} \}$, 
by non-increasing radii, and denote the ordered set by $\{ \mathbf{p}_{j_1},  \mathbf{p}_{j_2},   \mathbf{p}_{j_3} \}$, so that $r_{j_1} \geq r_{j_2} \geq r_{j_3}$.
By assumption, $r_{j_1} > r_{j_3}$.
The vectors $\mathbf{h}_e$ and $\mathbf{d}_e$ of the hyperplane
$H_e = \{ \mathbf{x} : \mathbf{h}_e\mathbf{x} = \mathbf{h}_e\mathbf{d}_e \}$, are computed by the following.
\begin{alignat}{4}
& \text{if\;\;} r_{j_1} = r_{j_2} > r_{j_3}\;\; &  \mathbf{h}_{e} &\;=\;  (\mathbf{p}_{j_1} -  \mathbf{p}_{j_2})/\|\mathbf{p}_{j_1} -  \mathbf{p}_{j_2}\|, \;\;\;
\mathbf{d}_e = .5(\mathbf{p}_{j_1} + \mathbf{p}_{j_2})   \notag \\
& \text{if\;\;} r_{j_1} > r_{j_2} = r_{j_3}\;\; &  \mathbf{h}_{e} &\;=\;  (\mathbf{p}_{j_2} -  \mathbf{p}_{j_3})/\|\mathbf{p}_{j_2} -  \mathbf{p}_{j_3}\|, \;\;\;
\mathbf{d}_e = .5(\mathbf{p}_{j_2} + \mathbf{p}_{j_3})  \notag \\
& \text{if\;\;} r_{j_1} > r_{j_2} > r_{j_3}\;\;  & \mathbf{h}_{e} &\;=\;  (\epsilon_{j_1,j_2}\mathbf{v}_{j_1,j_2} -  \epsilon_{j_1,j_3}\mathbf{v}_{j_1,j_3})/
\|\epsilon_{j_1,j_2}\mathbf{v}_{j_1,j_2} -  \epsilon_{j_1,j_3}\mathbf{v}_{j_1,j_3}\|  \label{he} \\
&&\mathbf{w}_{e} & \;=\; (\mathbf{h}_{e} - (\mathbf{h}_{e} \mathbf{v}_{j_1,j_3})\mathbf{v}_{j_1,j_3}) / 
\| \mathbf{h}_{e} - (\mathbf{h}_{e} \mathbf{v}_{j_1,j_3})\mathbf{v}_{j_1,j_3})\|  \notag \\
&&\mathbf{d}_{e} & \;=\;  \mathbf{d}_{j_1,j_3} + \frac{\mathbf{v}_{j_1,j_2} (\mathbf{d}_{j_1,j_2} - \mathbf{d}_{j_1,j_3})}{\mathbf{v}_{j_1,j_2} \mathbf{w}_{e}} \mathbf{w}_{e}.   \notag 
\end{alignat}
The intersection of $X_{S}$ and $H_e$ is equivalent to the intersection of $Y_{S}$ and $H_e$, and is determined by  
by solving the equation $\mathbf{h}_e\mathbf{y}(\beta') = \mathbf{h}_e\mathbf{d}_e$, for $\beta'$. 
Since  $Y_{S}$ is a hyperbola, with $\epsilon_S > 1$,  the equation $\mathbf{h}_e\mathbf{y}(\beta') = \mathbf{h}_e\mathbf{d}_e$
is equivalent to the equation
\begin{equation}
A \sec(\beta') + B \tan(\beta') = C,  \label{AAAA} 
\end{equation}
where \;
$A = a_{S} \mathbf{h}_{e}  \mathbf{v}_{S}, $ 
\; $B = b_{S} \mathbf{h}_{e}  \mathbf{u}_{S}, $ 
\;\;and \; $C =  \mathbf{h}_{e}   ( \mathbf{d}_{e} -  \mathbf{c}_{S}). $

Let $D_h = \sqrt{C^2-(A^2-B^2)}$. 
The value $C^2$ is the squared distance between  $H_e$ and the center $\mathbf{c}_S$. The value   $A^2-B^2$ is the
minimum squared distance between $\mathbf{c}_S$ and $H_e$ such that $H_e \cap Y_S \neq \emptyset$.
Thus, $H_e \cap Y_S \neq \emptyset$ if and only if $C^2 - (A^2 - B^2) \geq 0$, 

If $C^2 - (A^2 - B^2) > 0$, there are two  solutions, $\beta_1$ and $\beta_2$, determined by
\begin{alignat}{2}
& \tan(\beta_1) = \frac {-AB -CD_h} {AC - BD_h}, \;\;\;\text{and}\;\;\; \sec(\beta_1) = \frac {B^2 + C^2} {AC - BD_h} \notag \\
&\text{or} \label{hysold} \\
&  \tan(\beta_2) = \frac {-AB +CD_h} {AC + BD_h}, \;\;\;\text{and}\;\;\; \sec(\beta_2) = \frac {B^2 + C^2} {AC + BD_h}. \notag 
\end{alignat}

If $\sec(\beta_1) > 0$, then $\mathbf{y}(\beta_1) \in Y_S$, and  $\beta_1 := \tan^{-1}\left(  \frac {-AB -CD_h} {AC - BD_h} \right)$.
If  $\sec(\beta_1) < 0$, then $\mathbf{y}(\beta_1) \notin Y_S$, and $\beta_1 := \infty$.
Similarly, if $\sec(\beta_2) > 0$, then $\beta_2 := \tan^{-1}\left(  \frac {-AB + CD_h} {AC + BD_h} \right)$,
 and if  $\sec(\beta_2) < 0$, $\beta_2 := \infty$.  Then $\beta' = \min\{\beta_1, \beta_2 \}$.

If $C^2 - (A^2 - B^2) = 0$, there is one solution, $\beta' := \tan^{-1}(-B/C)$.
If $C^2 - (A^2 - B^2) < 0$, set $\beta':= \infty$.



The step size $\alpha'$ along the path $X_S$ from $\mathbf{x}_S$ is given by $\alpha' = \beta' - \beta_S$.

If $\mathbf{x}(\alpha') \in $ conv$(S \cup \{ \mathbf{p}_e \})$, then $\mathbf{x}(\alpha')$ is dual feasible. 
To determine if $\mathbf{x}(\alpha') \in$ conv$(S \cup \{ \mathbf{p}_e \})$, 
substitute $\mathbf{x}(\alpha) = \mathbf{c}_{S} + a_{S}\sec(\alpha)\mathbf{v}_{S} + b_{S}\tan(\alpha)\mathbf{u}_{S}$ 
for $\mathbf{x}_S$ in equation (\ref{denom2}),  expand and simplify to obtain the following linear system: 
\begin{alignat}{3}
T \boldsymbol{\pi}(\alpha) =  \boldsymbol{\varepsilon}_1 - a_S \sec(\alpha)\hat{\mathbf{v}}_{S}  - b_S \tan(\alpha)\hat{\mathbf{u}}_{S}, \label{Tuneqh}
\end{alignat}
where
\begin{alignat}{3}
T = \left[\!
\begin{array}{*{1}{cccccc}}
1&  \ldots & 1 &1 \\
\mathbf{c}_{S} - \mathbf{p}_{i_1}  & \ldots & \mathbf{c}_{S} - \mathbf{p}_{i_{s}}  & \mathbf{c}_{S} - \mathbf{p}_{e} \\
\end{array} \!
\right]\!, \; \boldsymbol{\varepsilon}_1 =  \left[ \!
\begin{array}{*{1}{c}}
1\\
\mathbf{0} \\
\end{array}
\right]\!, \;  \hat{\mathbf{v}}_{S}  =  \left[
\begin{array}{*{1}{c}}
0\\
\mathbf{v}_{S} \\
\end{array} \!
\right], \;  \hat{\mathbf{u}}_{S}  =  \left[ \!
\begin{array}{*{1}{c}}
0\\
\mathbf{u}_{S} \\
\end{array} \!
\right]\!, \label{Tevu}
\end{alignat}
and each variable is a function of  $\alpha$, that is, $\boldsymbol{\pi}(\alpha) = [\pi_{i_1}(\alpha), \ldots, \pi_{i_s}(\alpha), \pi_{e}(\alpha) ]$.
The matrix $T$ is $n+1$ by $s+1$ with rank $s+1$.  
The linear system (\ref{Tuneqh}) has a solution  since $\mathbf{x}(\alpha) \in $ aff$(S \cup \mathbf{p}_e)$ for $\alpha \geq 0$.

If $\boldsymbol{\pi}(\alpha') \geq 0$, then $\mathbf{x}(\alpha') \in $ conv$(S \cup \{ \mathbf{p}_e \})$.
In this case, reset $S:= S \cup \{\mathbf{p}_e\}$, and  $[\mathbf{x}_S, z_S ] :=  [\mathbf{x}(\alpha'),z(\alpha')]$.
Then the ball $[\mathbf{x}_S, z_S ]$ and its active set $S$ are checked for optimality by the Update Phase.

Otherwise, some component of $\boldsymbol{\pi}(\alpha')$ is negative, and $\mathbf{x}(\alpha') \notin$ conv$(S \cup \{ \mathbf{p}_e \})$. 
In this case, a new step size $\alpha^*$  is determined so that   $0 \leq  \alpha^* < \alpha'$, and $\mathbf{x}(\alpha^*)  \in$ conv$(S \cup \{ \mathbf{p}_e \})$. 
Since the $s+1$ points in $S \cup \{ \mathbf{p}_e\}$ are affinely independent,
conv$(S  \cup \{ \mathbf{p}_e\})$ is a simplex with $s+1$ vertices and $s+1$ facets.  The $s+1$ vertices are the points in $S \cup \{ \mathbf{p}_e\}$.  
The $s+1$  facets  are denoted by $F_{i_j} =$ conv$(S \cup \{ \mathbf{p}_e\} \setminus \{\mathbf{p}_{i_j} \})$, 
 for each $\mathbf{p}_{i_j} \in S$, and by $F_e = $ conv$(S)$, for $\mathbf{p}_e$.
Each  point $\mathbf{p}_{i_j} \in S$ corresponds to the 
component $\pi_{i_j}(\alpha)$ of $\boldsymbol{\pi}(\alpha)$ and  to the facet $F_{i_j}$.
The  point $\mathbf{p}_e$ corresponds to the component $\pi_e(\alpha)$ of $\boldsymbol{\pi}(\alpha)$ and  to the facet $F_e$. 

Since $\mathbf{x}(0) \in$ conv$(S \cup \{ \mathbf{p}_e \})$,  and $\mathbf{x}(\alpha') \notin$ conv$(S \cup \{ \mathbf{p}_e \})$, 
$X_{S}$ must intersect some facet of conv$(S \cup \{ \mathbf{p}_e \})$
between $\mathbf{x}(0)$ and $\mathbf{x}(\alpha')$.
 For each $\mathbf{p}_{i_j} \in S$, the  step-size $\alpha_{i_j} \geq 0$ is computed, if it exists, 
 so that $\mathbf{x}(\alpha_{i_j}) \in X_S \cap F_{i_j}$,
 which is equivalent to  $\pi_{i_j}(\alpha_{i_j}) = 0$ in the solution to 
$ T \boldsymbol{\pi}(\alpha_{i_j}) =  \boldsymbol{\varepsilon}_1 - a_S \sec(\alpha_{i_j})\hat{\mathbf{v}}_{S}  - b_S \tan(\alpha_{i_j})\hat{\mathbf{u}}_{S}$.
 
 To determine $\alpha_{i_j}$ for each $\mathbf{p}_{i_j} \in S$, so that $\pi_{i_j}(\alpha_{i_j}) = 0$,
 in the solution /$\boldsymbol{\pi}(\alpha_{i_j})$ to 
 $ T \boldsymbol{\pi}(\alpha_{i_j}) =  \boldsymbol{\varepsilon}_1 - a_S \sec(\alpha_{i_j})\hat{\mathbf{v}}_{S}  - b_S \tan(\alpha_{i_j})\hat{\mathbf{u}}_{S}$, 
solve $T \boldsymbol{\gamma}  = \boldsymbol{\varepsilon}_1$ for $\boldsymbol{\gamma}$, 
 solve $T \boldsymbol{\delta} = \hat{\mathbf{v}}_S$ for $\boldsymbol{\delta}$, and
  and  solve $T \boldsymbol{\xi} = \hat{\mathbf{u}}_S$ for $\boldsymbol{\xi}$.
 Then, $ T \boldsymbol{\pi}(\alpha_{i_j}) =  \boldsymbol{\varepsilon}_1 - a_S \sec(\alpha_{i_j})\hat{\mathbf{v}}_{S}  - b_S \tan(\alpha_{i_j})\hat{\mathbf{u}}_{S}
  = T( \boldsymbol{\gamma} - a_S \sec(\alpha_{i_j})\boldsymbol{\delta}   - b_S \tan(\alpha_{i_j})\boldsymbol{\xi})$.
If $ \gamma_{i_j} - a_S \sec(\alpha_{i_j})\delta_{i_j}   - b_S \tan(\alpha_{i_j})\xi_{i_j} = 0$, then $\pi_{i_j}(\alpha_{i_j}) = 0$.
Thus, for each $\mathbf{p}_{i_j} \in S$, solve $A\sec(\alpha_{i_j}) + B\tan(\alpha_{i_j}) = C$, for $\alpha_{i_j}$, using (\ref{BBB1}), where
$A = a_S\delta_{i_j}$, $B = b_S\xi_{i_j}$, and $C = \gamma_{i_j}$.
 If $\alpha_{i_j} < 0$ or $\alpha_{i_j} > \alpha'$, set $\alpha_{i_j} = \infty$.  
 If $C^2-(A^2-B^2) < 0$, set $\alpha_{i_j} = \infty$.  If $C^2-(A^2-B^2) = 0$, set $\alpha_{i_j} = \tan^{-1}(B/A)$.
 If $C^2-(A^2-B^2) > 0$, there are two solutions $\beta_1$ and  $\beta_2$ determined by (\ref{hysold}) so that
  $\alpha_{i_j} =\min \{ \beta_1, \beta_2 \}$.

 For the point $\mathbf{p}_e$, $\mathbf{x}(\alpha_e) \in$ conv$(S)$ implies $\alpha_e \leq 0$.  
 That is,  if $ \mathbf{x}_S \in$ conv$(S)$, then $\mathbf{x}(0) = \mathbf{x}_S$, so that $\alpha_e = 0$.
 If  $\mathbf{x}_S \notin $ conv$(S)$, then $\mathbf{x}(\alpha_e) \in$ aff$(S)$, where $\alpha_e = \frac{(\mathbf{p}_{i_1} - \mathbf{x}_S)\mathbf{d}_S}{\mathbf{d}_S^2} < 0$.
 
Thus the intersection of $X_S$ and the facet of conv$(S \cup \{ \mathbf{p}_e \})$ first encountered along $X_S$ occurs 
at the step size $\alpha^*$, where
\begin{equation}
 \alpha^*  = \min_{\mathbf{p}_{i_j} \in S } \{ \alpha_{i_j}   \}. \label{alst2}
\end{equation}
The leaving point $\mathbf{p}_{i_l}  \in S$  is chosen so that $\alpha^* = \alpha_{i_l}$.
Then the solution to the linear system (\ref{Tuneqh}) with $\alpha = \alpha^*$  yields
$\pi_{i_l}(\alpha^*) = 0$, and $\pi_{i_j} (\alpha^*) \geq 0$, for $\mathbf{p}_{i_j} \in S \cup \{ \mathbf{p}_e \} \setminus \{ \mathbf{p}_{i_l} \}$ 
so that $\mathbf{x}(\alpha^*) \in F_{i_l}$.
This construction shows that  $S  \setminus \{ \mathbf{p}_{i_l} \}$ is an active set for $[\mathbf{x}(\alpha^*), z(\alpha^*)]$
and that  $[\mathbf{x}(\alpha^*), z(\alpha^*)]$ is dual feasible with respect to $S \setminus \{ \mathbf{p}_{i_l} \} \cup \{ \mathbf{p}_e \}$.
However, the constraint corresponding to $\mathbf{p}_e$ is not active at $[\mathbf{x}(\alpha^*), z(\alpha^*)]$.
 
 Delete the point $\mathbf{p}_{i_l}$ from $S$.  Reset $S: =  S \setminus \{ \mathbf{p}_{i_l} \}$, 
 and  $[\mathbf{x}_S, z_S] := [ \mathbf{x}(\alpha^*), z(\alpha^*) ]$.
The set $S$ is active with respect to the ball $[\mathbf{x}_S, z_{S}]$, and $[\mathbf{x}_S, z_{S}]$ is dual feasible with respect to 
$S \cup \{ \mathbf{p}_e \} $.  
The Search Phase is re-entered with the active set $S$, the entering point $\mathbf{p}_e$, and the ball $[\mathbf{x}_S, z_{S}]$.

\noindent \textbf{Case 2b: $X_S$ is an ellipse:}
If  $\epsilon_{S} < 1$, then $B_{S} $ is an ellipsoid of dimension $n-s+2$, and 
\begin{equation}
 Y_{S}  =  \{ \mathbf{y}(\beta) =  \mathbf{c}_{S} + a_{S} \cos(\beta) \mathbf{v}_{S}  + b_{S} \sin(\beta) \mathbf{u}_{S} : 0 \leq \beta \leq 2\pi \} \label{ysduale}
 \end{equation}
is a two-dimensional ellipse in  aff$(\mathbf{c}_{S}, \mathbf{v}_{S}, \mathbf{u}_{S})$.  
The vector $\mathbf{v}_S$ is oriented as follows so that the objective function value is increasing along the search path.  
If $(\mathbf{x}_S - \mathbf{c}_S)\mathbf{v}_S < 0$, reset $\mathbf{v}_S := - \mathbf{v}_S$.
The vector $\mathbf{u}_{S}$ is  chosen to be the normalized component of  $(\mathbf{p}_e - \mathbf{c}_{S})$ 
that is orthogonal to the projection of  $(\mathbf{p}_e - \mathbf{c}_{S})$
 onto sub($S$), that is, 
 \begin{equation}
 \mathbf{u}_{S} = ((\mathbf{p}_{e} - \mathbf{c}_{S}) - \text{Proj}_{\text{sub}(S)} (\mathbf{p}_{e} - \mathbf{c}_{S})) / \| ((\mathbf{p}_{e} - \mathbf{c}_{S}) - \text{Proj}_{\text{sub}(S)} (\mathbf{p}_{e} - \mathbf{c}_{S})) \| \label{uella}.
 \end{equation}

Property A.9 shows that $Y_{S} $ is a subset of $B_{S}$, and that
 the points $\mathbf{y}(0)$ and $\mathbf{y}(\pi)$ are vertices of the ellipse $Y_{S} $ and of the ellipsoid $B_{S}$.
  Let $\beta_S$ be the parameter such that $\mathbf{y}(\beta_S) = \mathbf{x}_S$.  
 If $\mathbf{x}_S = \mathbf{a}_S$, then $\beta_S = 0$, but if  $\mathbf{x}_S \neq \mathbf{a}_S$, then
 $\beta_S =  \cos^{-1} \{\sqrt{(\mathbf{x}_S - \mathbf{c}_S )^2 - b_S^2}/(\epsilon_Sa_S) \}.$
 If a point was deleted from $S$ in the last iteration, $\mathbf{x}_S \neq \mathbf{a}_S$, else, $\mathbf{x}_S = \mathbf{a}_S$.
 
 The search path $X_S$ is defined to be a subset of $Y_S$ given by 
 \begin{equation}
 X_S = \{ \mathbf{x}(\alpha) = \mathbf{y}(\alpha + \beta_S), \alpha \geq 0\}. \label{xsduale}
 \end{equation} 
 Observe that $\mathbf{x}(0) = \mathbf{y}(\beta_S) = \mathbf{x}_S$, and for $\alpha \geq 0$, 
 the points $\mathbf{x}(\alpha)$ are on $B_S$ and move along $X_S$ from $\mathbf{x}_S$ toward $\mathbf{p}_e$ and primal feasibility.

  For the case when $X_S$ is an ellipse, the following property and its proof are analogous to Property 6.3 for the case when $X_S$ is based on a hyperbola.
 \begin{property}
Given the active set  $S =  \{ \mathbf{p}_{i_1}, \ldots,  \mathbf{p}_{i_{s}} \}$  for the  primal feasible ball $[\mathbf{x}_S, z_S]$,
ordered so that $r_{i_1} \geq \ldots \geq r_{i_{s}}$, with $r_{i_1} > r_{i_{s}}$.  
Let $X_S = \{ \mathbf{x}(\alpha) = \mathbf{y}(\alpha+\beta_S), 0 \leq \alpha \leq -\beta_S \}$  
be the search path along the ellipse  $\mathbf{x}(\alpha)$ determined by (\ref{ysduale}) and (\ref{xsduale}) .
Then $S$ is an active set for the ball $[\mathbf{x}(\alpha), z(\alpha) ]$, for $\alpha \geq 0$, where
$z(\alpha) = \| \mathbf{x}(\alpha) - \mathbf{p}_{i_1} \| + r_{i_1} $.
Furthermore, $z(\alpha)$ is increasing for $0 \leq \alpha \leq  -\beta_S$.
 \end{property}
  
 \noindent \begin{bfseries}The Step Size: \end{bfseries}
 The step size $\alpha'$ is determined, if it exists, so that $X_{S}$ intersects the bisector  $B_{i_1,e}$ at $\mathbf{x}(\alpha')$.
 If $\mathbf{x}(\alpha') \in X_S \cap B_{i_1,e}$, then $\mathbf{x}(\alpha') \in X_S \cap B_{i_j,e}$, for all  $\mathbf{p}_{i_j} \in S$, that is, 
$X_S$ simultaneously intersects the bisectors $B_{i_j,e}$ at $\mathbf{x}(\alpha')$, for all $\mathbf{p}_{i_j} \in S$.
Thus, it suffices to determine the intersection of $X_S$ with only $B_{i_1,e}$.
 
 At the  point $\mathbf{x}(\alpha') \in X_S \cap B_{i_1,e}$, 
 $ z(\alpha') = \| \mathbf{x}(\alpha') - \mathbf{p}_e \| + r_e =  \| \mathbf{x}(\alpha') - \mathbf{p}_{i_1} \| + r_{i_1}$,
 that is, $S \cup \{ \mathbf{p}_e \}$ is an active set for the ball $[\mathbf{x}(\alpha'),z(\alpha')]$.
  Geometrically, the search moves the center  $\mathbf{x}(\alpha)$ along the path $X_{S}$ while increasing  
 $z(\alpha)$.
 
Theorem A.1 shows that $X_S \cap B_{i_1,e} = X_S \cap H_e$, where 
the hyperplane $H_e =  \{ \mathbf{x} : \mathbf{h}_e\mathbf{x} = \mathbf{h}_e\mathbf{d}_e \}$ is constructed as follows.
Order the set of three points $\{\mathbf{p}_{i_1}, \mathbf{p}_{i_{s}}, \mathbf{p}_{e} \}$
by non-increasing radii, and denote the ordered set by $\{ \mathbf{p}_{j_1},  \mathbf{p}_{j_2},   \mathbf{p}_{j_3} \}$, so that $r_{j_1} \geq r_{j_2} \geq r_{j_3}$.
By assumption, $r_{j_1} > r_{j_3}$.
The vectors $\mathbf{h}_e$ and $\mathbf{d}_e$ of the hyperplane
$H_e = \{ \mathbf{x} : \mathbf{h}_e\mathbf{x} = \mathbf{h}_e\mathbf{d}_e \}$ are computed by the expressions (\ref{he}).
 
The intersection of $X_{S}$ and $H_e$ is equivalent to the intersection of $Y_{S}$ and $H_e$, and is determined 
by solving the equation $\mathbf{h}_e\mathbf{y}(\beta') = \mathbf{h}_e\mathbf{d}_e$ for $\beta'$. 
Since  $Y_{S}$ is an ellipse, with $\epsilon_S < 1$,  the equation $\mathbf{h}_e\mathbf{y}(\beta') = \mathbf{h}_e\mathbf{d}_e$
is equivalent to the equation
\begin{equation}
A \cos(\beta') + B \sin(\beta') = C,  \label{AAAE} 
\end{equation}
where \;
$A = a_{S} \mathbf{h}_{e}  \mathbf{v}_{S}, $ 
\; $B = b_{S} \mathbf{h}_{e}  \mathbf{u}_{S}, $ 
\;\;and \; $C =  \mathbf{h}_{e}   ( \mathbf{d}_{e} -  \mathbf{c}_{S}). $

Let $D_e = \sqrt{A^2 +B^2-C^2}$. 
The value $C^2$ is the squared distance between  $H_e$ and the center $\mathbf{c}_S$. The expression   $A^2+B^2$ is the
maximum squared distance between $\mathbf{c}_S$ and $H_e$ such that $H_e \cap X_S \neq \emptyset$.
Thus, $H_e \cap X_S \neq \emptyset$ if and only if $A^2 + B^2 -C^2 \geq 0$,

If $A^2 + B^2 -C^2 \geq 0$, there are two possible solutions, $\beta_1$ and $\beta_2$,  determined by
\begin{alignat}{2}
& \cos(\beta_1) = \frac {AC + BD_e} {A^2 + B^2}, \;\;\;\text{and}\;\;\; \sin(\beta_1) = \frac {BC -AD_e} {A^2 + B^2} \notag \\
&\text{or} \label{elsola} \\
&  \cos(\beta_2) = \frac {AC - BD_e} {A^2 + B^2}, \;\;\;\text{and}\;\;\; \sin(\beta_2) = \frac {BC + AD_e} {A^2 + B^2} \notag 
\end{alignat}

If $\sin(\beta_1) > 0$, then $\mathbf{y}(\beta_1) \in X_S$, and  $\beta_1 := \tan^{-1}\left(  \frac {BC -AD_e} {AC + BD_e} \right)$.
If  $\sin(\beta_1) < 0$, then $\mathbf{y}(\beta_1) \notin X_S$, and $\beta_1 := \infty$.
Similarly, if $\sin(\beta_2) > 0$, then $\beta_2 := \tan^{-1}\left(  \frac {BC + AD_e} {AC - BD_e} \right)$,
 and if  $\sec(\beta_2) < 0$, $\beta_2 := \infty$.
 Then, $\beta' = \min\{ \beta_1, \beta_2 \}$.

If $A^2 + B^2 -C^2 = 0$, there is one solution $\beta' := \tan^{-1}(B/A)$.
 If $A^2 + B^2 -C^2 < 0$, set $\beta' := \infty$.


The step size $\alpha'$ along the path $X_S$ from $\mathbf{x}_S$ is given by $\alpha' = \beta' - \beta_S$.

If $\mathbf{x}(\alpha') \in $ conv$(S \cup \{ \mathbf{p}_e \})$, then $\mathbf{x}(\alpha')$ is dual feasible. 
To determine if $\mathbf{x}(\alpha') \in$ conv$(S \cup \{ \mathbf{p}_e \})$, 
substitute $\mathbf{x}(\alpha) = \mathbf{c}_{S} + a_{S}\cos(\alpha)\mathbf{v}_{S} + b_{S}\sin(\alpha)\mathbf{u}_{S}$ 
for $\mathbf{x}_S$ in equation (\ref{denom2}),  expand and simplify to obtain the following linear system: 
\begin{alignat}{3}
T \boldsymbol{\pi}(\alpha) =  \boldsymbol{\varepsilon}_1 - a_S \cos(\alpha)\hat{\mathbf{v}}_{S}  - b_S \sin(\alpha)\hat{\mathbf{u}}_{S} \label{Teqhe}
\end{alignat}
where $T$, $\boldsymbol{\varepsilon}_1$, $\hat{\mathbf{v}}_{S}$, and $\hat{\mathbf{u}}_{S}$ are defined by (\ref{Tevu}),
and each variable is a function of  $\alpha$, that is, $\boldsymbol{\pi}(\alpha) = [\pi_{i_1}(\alpha), \ldots, \pi_{i_s}(\alpha), \pi_{e}(\alpha) ]$.
The matrix $T$ is $s+1 \times n+1$ with rank $s+1$.  
The linear system (\ref{Teqhe}) has a solution  since $\mathbf{x}(\alpha) \in $ aff$(S \cup \mathbf{p}_e)$ for $\alpha \geq 0$.

If $\boldsymbol{\pi}(\alpha') \geq 0$, then $\mathbf{x}(\alpha') \in $ conv$(S \cup \{ \mathbf{p}_e \})$.
In this case, reset $S:= S \cup \{\mathbf{p}_e\}$, and  $[\mathbf{x}_S, z_S ] :=  [\mathbf{x}(\alpha'),z(\alpha')]$.
Then the ball $[\mathbf{x}_S, z_S ]$ and its active set $S$ are checked for optimality by the Update Phase.

Otherwise, some component of $\boldsymbol{\pi}(\alpha')$ is negative, and $\mathbf{x}(\alpha') \notin$ conv$(S \cup \{ \mathbf{p}_e \})$. 
In this case, a new step size $\alpha^*$  is determined so that   $0 \leq  \alpha^* <  \alpha'$, and $\mathbf{x}(\alpha^*)  \in$ conv$(S \cup \{ \mathbf{p}_e \})$. 
Since the $s+1$ points in $S \cup \{ \mathbf{p}_e\}$ are affinely independent,
conv$(S  \cup \{ \mathbf{p}_e\})$ is a simplex with $s+1$ vertices and $s+1$ facets.  The $s+1$ vertices are the points in $S \cup \{ \mathbf{p}_e\}$.  
The $s+1$  facets  are denoted by $F_{i_j} =$ conv$(S \cup \{ \mathbf{p}_e\} \setminus \{\mathbf{p}_{i_j} \})$, 
 for each $\mathbf{p}_{i_j} \in S$, and by $F_e = $ conv$(S)$, for $\mathbf{p}_e$.
Each  point $\mathbf{p}_{i_j} \in S$ corresponds to the 
component $\pi_{i_j}(\alpha)$ of $\boldsymbol{\pi}(\alpha)$ and  to the facet $F_{i_j}$.
The  point $\mathbf{p}_e$ corresponds to the component $\pi_e(\alpha)$ of $\boldsymbol{\pi}(\alpha)$ and  to the facet $F_e$. 

Since $\mathbf{x}(0) \in$ conv$(S \cup \{ \mathbf{p}_e \})$,  and $\mathbf{x}(\alpha') \notin$ conv$(S \cup \{ \mathbf{p}_e \})$, 
$X_{S}$ must intersect some facet of conv$(S \cup \{ \mathbf{p}_e \})$
between $\mathbf{x}(0)$ and $\mathbf{x}(\alpha')$.
 For each $\mathbf{p}_{i_j} \in S$, the  step-size $\alpha_{i_j} \geq 0$ is computed, if it exists, 
 so that $\mathbf{x}(\alpha_{i_j}) \in X_S \cap F_{i_j}$,
 which is equivalent to  $\pi_{i_j}(\alpha_{i_j}) = 0$ in the solution to 
$ T \boldsymbol{\pi}(\alpha_{i_j}) =  \boldsymbol{\varepsilon}_1 - a_S \cos(\alpha_{i_j})\hat{\mathbf{v}}_{S}  - b_S \sin(\alpha_{i_j})\hat{\mathbf{u}}_{S}$.

To determine $\alpha_{i_j}$ for each $\mathbf{p}_{i_j} \in S$, so that $\pi_{i_j}(\alpha_{i_j}) = 0$,
 in the solution $\boldsymbol{\pi}(\alpha_{i_j})$ to 
 $ T \boldsymbol{\pi}(\alpha) =  \boldsymbol{\varepsilon}_1 - a_S \cos(\alpha)\hat{\mathbf{v}}_{S}  - b_S \sin(\alpha)\hat{\mathbf{u}}_{S}$, 
solve $T \boldsymbol{\gamma}  = \boldsymbol{\varepsilon}_1$ for $\boldsymbol{\gamma}$, 
 solve $T \boldsymbol{\delta} = \hat{\mathbf{v}}_S$ for $\boldsymbol{\delta}$, 
  and  solve $T \boldsymbol{\xi} = \hat{\mathbf{u}}_S$ for $\boldsymbol{\xi}$.
 Then, $ T \boldsymbol{\pi}(\alpha) =  \boldsymbol{\varepsilon}_1 - a_S \cos(\alpha)\hat{\mathbf{v}}_{S}  - b_S \sin(\alpha)\hat{\mathbf{u}}_{S}
  = T( \boldsymbol{\gamma} - a_S \cos(\alpha_{i_j})\boldsymbol{\delta}   - b_S \sin(\alpha_{i_j})\boldsymbol{\xi})$.
If $ \gamma_{i_j} - a_S \cos(\alpha_{i_j})\delta_{i_j}   - b_S \sin(\alpha_{i_j})\xi_{i_j} = 0$, then $\pi_{i_j}(\alpha_{i_j}) = 0$.
Thus, for each $\mathbf{p}_{i_j} \in S$, solve $A\cos(\alpha_{i_j}) + B\sin(\alpha_{i_j}) = C$, for $\alpha_{i_j}$, using (\ref{BBB2}), where
$A = a_S\delta_{i_j}$, $B = b_S\xi_{i_j}$, and $C = \gamma_{i_j}$.
 If $\alpha_{i_j} < 0$ or $\alpha_{i_j} > \alpha'$, set $\alpha_{i_j}: = \infty$.  
 If $A^2+B^2-C^2 < 0$, set $\alpha_{i_j} = \infty$.  If $A^2+B^2-C^2 = 0$, set $\alpha_{i_j} = \tan^{-1}(B/A)$.
 If $A^2+B^2-C^2 > 0$, there are two solutions $\beta_1$ and  $\beta_2$ determined by (\ref{elsola}), so that 
  $\alpha_{i_j} = \min\{ \beta_1,\beta_2 \}$.

 For the point $\mathbf{p}_e$, $\mathbf{x}(\alpha_e) \in$ conv$(S)$ implies $\alpha_e \leq 0$.  
 That is,  if $ \mathbf{x}_S \in$ conv$(S)$, then $\mathbf{x}(0) = \mathbf{x}_S$, so that $\alpha_e = 0$.
 If  $\mathbf{x}_S \notin $ conv$(S)$, then $\mathbf{x}(\alpha_e) \in$ aff$(S)$, where $\alpha_e = \frac{(\mathbf{p}_{i_1} - \mathbf{x}_S)\mathbf{d}_S}{\mathbf{d}_S^2} < 0$.

Thus, the intersection of $X_S$ and the facet of conv$(S \cup \{ \mathbf{p}_e \})$ first encountered along $X_S$ occurs 
at the step size $\alpha^*$, where
\begin{equation}
 \alpha^*  = \min_{\mathbf{p}_{i_j} \in S } \{ \alpha_{i_j}   \}.  \label{alp3}
\end{equation}
Choose $\mathbf{p}_{i_l}  \in S$ so that $\alpha^* = \alpha_{i_l}$.
Then the solution to the linear system (\ref{Teqhe}) with $\alpha = \alpha^*$  yields
$\pi_{i_l}(\alpha^*) = 0$, and $\pi_{i_j} (\alpha^*) \geq 0$, for $\mathbf{p}_{i_j} \in S \cup \{ \mathbf{p}_e \} \setminus \{ \mathbf{p}_{i_l} \}$ 
so that $\mathbf{x}(\alpha^*) \in F_{i_l}$.
This construction shows that  $S  \setminus \{ \mathbf{p}_{i_l} \}$ is an active set for $[\mathbf{x}(\alpha^*), z(\alpha^*)]$
and that  $[\mathbf{x}(\alpha^*), z(\alpha^*)]$ is dual feasible with respect to $S \setminus \{ \mathbf{p}_{i_l} \} \cup \{ \mathbf{p}_e \}$.
However, the constraint corresponding to $\mathbf{p}_e$ is not active at $[\mathbf{x}(\alpha^*), z(\alpha^*)]$.
 
 Delete the point $\mathbf{p}_{i_l}$ from $S$.  Reset $S: =  S \setminus \{ \mathbf{p}_{i_l} \}$, 
 and  $[\mathbf{x}_S, z_S] := [ \mathbf{x}(\alpha^*), z(\alpha^*) ]$.
The set $S$ is active with respect to the ball $[\mathbf{x}_S, z_{S}]$, and $[\mathbf{x}_S, z_{S}]$ is dual feasible with respect to 
$S \cup \{ \mathbf{p}_e \} $.  
The Search Phase is re-entered with the active set $S$, the entering point $\mathbf{p}_e$, and the ball $[\mathbf{x}_S, z_{S}]$.

If there is more than one  point
$\mathbf{p}_{i_j} \in S\setminus \{\mathbf{p}_{e}\}$ such that $\alpha^* = \alpha_{i_l}$, then  $X_S$
intersects the corresponding facets $F_{i_j}$ simultaneously, and there is a tie for the leaving point. 
In this case choose any point $\mathbf{p}_{i_j}$
 such that $\alpha_{i_j} = \alpha^*$, to be deleted from $S$. 
The set $S$ is reset to $S:=S \cup \{ \mathbf{p}_{i_j} \|$. 
The Search Phase is re-entered with the active set $S$,
the entering point $\mathbf{p}_e$, and the ball $[\mathbf{x}(\alpha^*), z(\alpha^*) ]$,
and the search continues on the facet $F_{i_j}$.
This case may lead to a degenerate iteration at  the next step
with $\alpha^* = \alpha_{i_l} =0$ for some $\mathbf{p}_{i_l}$.  
Cycling will not occur since at each degenerate iteration one point is deleted from the finite set $S$.  
After a finite number of points are deleted, 
 $S$ is reduced to two points, and the step size $\alpha$ will be positive at the next iteration.

\noindent \textbf{Case 2c: $X_S$ is a parabola:}
If  $\epsilon_{S} = 1$, then $B_{S} $ is a paraboloid of dimension $n-s+2$, and 
\begin{equation}
 Y_{S}  =  \{ \mathbf{y}(\beta) =  \hat{\mathbf{c}}_{S} + \tilde{c}_{S} \beta^2 \mathbf{v}_{S}  + 2 \tilde{c}_{S}\beta \mathbf{u}_{S} : -\infty < \beta < \infty \} \label{ysdualp}
 \end{equation}
is a two-dimensional parabola in  aff$( \hat{\mathbf{c}}_{S}, \mathbf{v}_{S}, \mathbf{u}_{S})$, where $\hat{\mathbf{c}}_S$ and  $\tilde{c}_S$ are determined by (\ref{pch}) - (\ref{pnc}).
The vector $\mathbf{u}_{S}$ is  chosen to be the normalized component of  $(\mathbf{p}_e -\hat{\mathbf{c}}_{S})$ 
that is orthogonal to the projection of  $(\mathbf{p}_e - \hat{\mathbf{c}}_{S})$
 onto sub($S$), that is, 
 \begin{equation}
 \mathbf{u}_{S} = ((\mathbf{p}_{e} - \hat{\mathbf{c}}_{S}) - \text{Proj}_{\text{sub}(S)} (\mathbf{p}_{e} - \hat{\mathbf{c}}_{S})) / \| ((\mathbf{p}_{e} - \hat{\mathbf{c}}_{S}) - \text{Proj}_{\text{sub}(S)} (\mathbf{p}_{e} - \hat{\mathbf{c}}_{S})) \| \label{uell}.
 \end{equation}

Property A.9 shows that $Y_{S} $ is a subset of $B_{S}$, and that
 the point $\mathbf{y}(0)$ is the vertex of the parabola $Y_{S} $ and of the paraboloid $B_{S}$.
  Let $\beta_S$ be the parameter such that $\mathbf{y}(\beta_S) = \mathbf{x}_S$.  
 If $\mathbf{x}_S = \hat{\mathbf{c}}_S$, then $\beta_S = 0$, but if  $\mathbf{x}_S \neq \hat{\mathbf{c}}_S$, then
 $\beta_S = \sqrt{  -2 +  \sqrt{ 4 +(\mathbf{x}_S - \hat{\mathbf{c}}_S )^2 /(\tilde{c}_S^2 ) }}.$ 
 If a point was deleted from $S$ in the last iteration, $\mathbf{x}_S \neq \hat{\mathbf{c}}_S$, else, $\mathbf{x}_S = \hat{\mathbf{c}}_S$.
 
  The search path $X_S$ is defined to be a subset of $Y_S$ given by 
 \begin{equation}
 X_S = \{ \mathbf{x}(\alpha) = \mathbf{y}(\alpha + \beta_S), \alpha \geq 0\}. \label{xsdualp}
 \end{equation} 
 Observe that $\mathbf{x}(0) = \mathbf{y}(\beta_S) = \mathbf{x}_S$, and for $\alpha \geq 0$, 
 the points $\mathbf{x}(\alpha)$ are on $B_S$ and move along $X_S$ from $\mathbf{x}_S$ toward $\mathbf{p}_e$ and primal feasibility.

  For the case when $X_S$ is a parabola, the following property and its proof are analogous to Property 6.3 for the case when $X_S$ is based on a hyperbola.
 \begin{property}
Given the active set  $S =  \{ \mathbf{p}_{i_1}, \ldots,  \mathbf{p}_{i_{s}} \}$  for the  primal feasible ball $[\mathbf{x}_S, z_S]$,
ordered so that $r_{i_1} \geq \ldots \geq r_{i_{s}}$, with $r_{i_1} > r_{i_{s}}$.  
Let $X_S = \{ \mathbf{x}(\alpha) = \mathbf{y}(\alpha+\beta_S), 0 \leq \alpha \leq -\beta_S \}$  
be the search path along the parabola  $\mathbf{x}(\alpha)$ determined by (\ref{ysdualp}) and (\ref{xsdualp}) .
Then $S$ is an active set for the ball $[\mathbf{x}(\alpha), z(\alpha) ]$, for $\alpha \geq 0$, where
$z(\alpha) = \| \mathbf{x}(\alpha) - \mathbf{p}_{i_1} \| + r_{i_1} $.
Furthermore, $z(\alpha)$ is increasing for $0 \leq \alpha \leq  -\beta_S$.
 \end{property}
  
 \noindent \begin{bfseries}The Step Size: \end{bfseries}
 The step size $\alpha'$ is determined, if it exists, so that $X_{S}$ intersects the bisector  $B_{i_1,e}$ at $\mathbf{x}(\alpha')$.
 If $\mathbf{x}(\alpha') \in X_S \cap B_{i_1,e}$, then $\mathbf{x}(\alpha') \in X_S \cap B_{i_j,e}$, for all  $\mathbf{p}_{i_j} \in S$, that is, 
$X_S$ simultaneously intersects the bisectors $B_{i_j,e}$ at $\mathbf{x}(\alpha')$, for all $\mathbf{p}_{i_j} \in S$.
Thus, it suffices to determine the intersection of $X_S$ with only $B_{i_1,e}$.
 
 At the  point $\mathbf{x}(\alpha') \in X_S \cap B_{i_1,e}$, 
 $ z(\alpha') = \| \mathbf{x}(\alpha') - \mathbf{p}_e \| + r_e =  \| \mathbf{x}(\alpha') - \mathbf{p}_{i_1} \| + r_{i_1}$,
 that is, $S \cup \{ \mathbf{p}_e \}$ is an active set for the ball $[\mathbf{x}(\alpha'),z(\alpha')]$.
  Geometrically, the search moves the center  $\mathbf{x}(\alpha)$ along the path $X_{S}$ while increasing  
 $z(\alpha)$.
 
Theorem A.1 shows that $X_S \cap B_{i_1,e} = X_S \cap H_e$, where 
the hyperplane $H_e =  \{ \mathbf{x} : \mathbf{h}_e\mathbf{x} = \mathbf{h}_e\mathbf{d}_e \}$ is constructed as follows.
Order the set of three points $\{\mathbf{p}_{i_1}, \mathbf{p}_{i_{s}}, \mathbf{p}_{e} \}$, 
by non-increasing radii, and denote the ordered set by $\{ \mathbf{p}_{j_1},  \mathbf{p}_{j_2},   \mathbf{p}_{j_3} \}$, so that $r_{j_1} \geq r_{j_2} \geq r_{j_3}$.
By assumption, $r_{j_1} > r_{j_3}$.
The vectors $\mathbf{h}_e$ and $\mathbf{d}_e$ of the hyperplane
$H_e = \{ \mathbf{x} : \mathbf{h}_e\mathbf{x} = \mathbf{h}_e\mathbf{d}_e \}$, are computed by the expressions (\ref{he}).
 
The intersection of $X_{S}$ and $H_e$ is equivalent to the intersection of $Y_{S}$ and $H_e$, and is determined by  
solving the equation $\mathbf{h}_e\mathbf{y}(\beta) = \mathbf{h}_e\mathbf{d}_e$ for $\beta$. 
which yields the quadratic equation  
\begin{equation}
A \beta^2+ B \beta = C,  \label{AADP} 
\end{equation}
where \;
$A = \tilde{c}_{S} \mathbf{h}_{e}  \mathbf{v}_{S}, $ 
\; $B = \tilde{c}_{S} \mathbf{h}_{e}  \mathbf{u}_{S}, $ 
\;\;and \; $C =  \mathbf{h}_{e}   ( \mathbf{d}_{e} -  \hat{\mathbf{c}}_{S}). $
If $Y_{S} \cap H_{e} \neq \emptyset$, there are possibly two intersection points
corresponding to the following two solutions to equation  (\ref{AADP}):

If there are no real solutions, $\beta' = \infty$.  If there is one real solution, $\beta'$.  Otherwise, there are two real solutions $\beta_1$ and $\beta_2$.  

If $B^2+A^2-C^2 <0$, there is no solution and no intersection point, and $\beta' = \infty$.  
If $B^2+A^2-C^2 =0$, or if $C = 0$, then $\beta_1 = \beta_2 = \beta'$.
Else, $B^2+A^2-C^2  > 0$, and there are two solutions  $\beta_1$ and $\beta_2$.
If $\beta_1 < \beta_S$, or $\mathbf{h}_e\mathbf{y}(\beta_1) \neq \mathbf{h}_e\mathbf{d}_e$, set $\beta_1 := \infty$. 
If $\beta_2 < \beta_S$, or $\mathbf{h}_e\mathbf{y}(\beta_2) \neq \mathbf{h}_e\mathbf{d}_e$, set $\beta_2 := \infty$. 
Then $\beta' = \min\{\beta_1, \beta_2 \}$.
The step size $\alpha'$ along the path $X_S$ from $\mathbf{x}_S$ is given by $\alpha' = \beta' - \beta_S$.

If $\mathbf{x}(\alpha') \in $ conv$(S \cup \{ \mathbf{p}_e \})$, then $\mathbf{x}(\alpha')$ is dual feasible. 
To determine if $\mathbf{x}(\alpha') \in$ conv$(S \cup \{ \mathbf{p}_e \})$, 
substitute $\mathbf{x}(\alpha) =  \hat{\mathbf{c}}_{S} + \tilde{c}_{S} \alpha^2 \mathbf{v}_{S}  + 2 \tilde{c}_{S}\alpha \mathbf{u}_{S} $
for $\mathbf{x}_S$ in equation (\ref{denom2}),  expand and simplify to obtain the following linear system: 
\begin{alignat}{3}
T \boldsymbol{\pi}(\alpha) =  \boldsymbol{\varepsilon}_1 - \tilde{c}_S \alpha^2 \hat{\mathbf{v}}_{S}  - \tilde{c}_S \alpha \hat{\mathbf{u}}_{S} \label{Teqpa}
\end{alignat}
where $T$, $\boldsymbol{\varepsilon}_1$, $\hat{\mathbf{v}}_{S}$, and $\hat{\mathbf{u}}_{S}$ are defined by (\ref{Tevu}).
and each variable is a function of  $\alpha$, that is, $\boldsymbol{\pi}(\alpha) = [\pi_{i_1}(\alpha), \ldots, \pi_{i_s}(\alpha), \pi_{e}(\alpha) ]$.
The matrix $T$ is $s+1 \times n+1$ with rank $s+1$.  
The linear system (\ref{Teqpa}) has a solution  since $\mathbf{x}(\alpha) \in $ aff$(S \cup \mathbf{p}_e)$ for $\alpha \geq 0$.

If $\boldsymbol{\pi}(\alpha') \geq 0$, then $\mathbf{x}(\alpha') \in $ conv$(S \cup \{ \mathbf{p}_e \})$.
In this case, reset $S:= S \cup \{\mathbf{p}_e\}$, and  $[\mathbf{x}_S, z_S ] :=  [\mathbf{x}(\alpha'),z(\alpha')]$.
Then the ball $[\mathbf{x}_S, z_S ]$ and its active set $S$ are checked for optimality by the Update Phase.

Otherwise, some component of $\boldsymbol{\pi}(\alpha')$ is negative, and $\mathbf{x}(\alpha') \notin$ conv$(S \cup \{ \mathbf{p}_e \})$. 
In this case, a new step size $\alpha^*$  is determined so that   $0 \leq  \alpha^* <  \alpha'$, and $\mathbf{x}(\alpha^*)  \in$ conv$(S \cup \{ \mathbf{p}_e \})$. 
Since the $s+1$ points in $S \cup \{ \mathbf{p}_e\}$ are affinely independent,
conv$(S  \cup \{ \mathbf{p}_e\})$ is a simplex with $s+1$ vertices and $s+1$ facets.  The $s+1$ vertices are the points in $S \cup \{ \mathbf{p}_e\}$.  
The $s+1$  facets  are denoted by $F_{i_j} =$ conv$(S \cup \{ \mathbf{p}_e\} \setminus \{\mathbf{p}_{i_j} \})$, 
 for each $\mathbf{p}_{i_j} \in S$, and by $F_e = $ conv$(S)$, for $\mathbf{p}_e$.
Each  point $\mathbf{p}_{i_j} \in S$ corresponds to the 
component $\pi_{i_j}(\alpha)$ of $\boldsymbol{\pi}(\alpha)$ and  to the facet $F_{i_j}$.
The  point $\mathbf{p}_e$ corresponds to the component $\pi_e(\alpha)$ of $\boldsymbol{\pi}(\alpha)$ and  to the facet $F_e$. 

Since $\mathbf{x}(0) \in$ conv$(S \cup \{ \mathbf{p}_e \})$,  and $\mathbf{x}(\alpha') \notin$ conv$(S \cup \{ \mathbf{p}_e \})$, 
$X_{S}$ must intersect some facet of conv$(S \cup \{ \mathbf{p}_e \})$
between $\mathbf{x}(0)$ and $\mathbf{x}(\alpha')$.
 For each $\mathbf{p}_{i_j} \in S$, the  step-size $\alpha_{i_j} \geq 0$ is computed, if it exists, 
 so that $\mathbf{x}(\alpha_{i_j}) \in X_S \cap F_{i_j}$,
 which is equivalent to  $\pi_{i_j}(\alpha_{i_j}) = 0$ in the solution to 
$ T \boldsymbol{\pi}(\alpha_{i_j}) =  \boldsymbol{\varepsilon}_1 - \tilde{c}_S \alpha_{i_j}\hat{\mathbf{v}}_{S}  - \tilde{c}_S \alpha_{i_j}\hat{\mathbf{u}}_{S}$.

To determine $\alpha_{i_j}$ for each $\mathbf{p}_{i_j} \in S$, so that $\pi_{i_j}(\alpha_{i_j}) = 0$,
 in the solution  to 
$ T \boldsymbol{\pi}(\alpha_{i_j}) =  \boldsymbol{\varepsilon}_1 - \tilde{c}_S \alpha_{i_j}^2 \hat{\mathbf{v}}_{S}  - \tilde{c}_S \alpha_{i_j}\hat{\mathbf{u}}_{S}$, 
solve $T \boldsymbol{\gamma}  = \boldsymbol{\varepsilon}_1$ for $\boldsymbol{\gamma}$, 
 solve $T \boldsymbol{\delta} = \hat{\mathbf{v}}_S$ for $\boldsymbol{\delta}$, 
  and  solve $T \boldsymbol{\xi} = \hat{\mathbf{u}}_S$ for $\boldsymbol{\xi}$.
 Then, $ T \boldsymbol{\pi}(\alpha_{i_j}) =  \boldsymbol{\varepsilon}_1 - \tilde{c}_{S} \alpha_{i_j}^2 \hat{\mathbf{v}}_{S}  - 2 \tilde{c}_{S}\alpha_{i_j} \hat{\mathbf{u}}_{S} 
  = T( \boldsymbol{\gamma} - \tilde{c}_S \alpha_{i_j}^2 \boldsymbol{\delta}   - \tilde{c}_S \alpha_{i_j}\boldsymbol{\xi})$.
If $ \gamma_{i_j} - \tilde{c}_S \alpha_{i_j}^2 \delta_{i_j}   - \tilde{c}_S \alpha_{i_j}\xi_{i_j} = 0$, then $\pi_{i_j}(\alpha_{i_j}) = 0$.
Thus, for each $\mathbf{p}_{i_j} \in S$, solve $A\alpha_{i_j}^2 + B\alpha_{i_j} = C$, for $\alpha_{i_j}$, where
$A = \tilde{c}_S\delta_{i_j}$, $B = \tilde{c}_S\xi_{i_j}$, and $C = \gamma_{i_j}$.
 If $\alpha_{i_j} < 0$ or $\alpha_{i_j} > \alpha'$, set $\alpha_{i_j} := \infty$.  
 If $B^2+C^2-A^2 < 0$, set $\alpha_{i_j} := \infty$.  If $B^2+C^2-A^2 = 0$, set $\alpha_{i_j} := -B/(2A)$.
 If $B^2+C^2-A^2 > 0$, there are two solutions $\beta_1$ and  $\beta_2$.
 If  $ \gamma_{i_j} - \tilde{c}_S \beta_1^2 \delta_{i_j}   - \tilde{c}_S \beta_1\xi_{i_j} \neq 0$, or if $\beta_1 \leq \beta_S$, set $\beta_1 := \infty$.
If  $ \gamma_{i_j} - \tilde{c}_S \beta_2^2 \delta_{i_j}   - \tilde{c}_S \beta_2\xi_{i_j} \neq 0$, or if $\beta_2 \leq \beta_S$, set $\beta_2 := \infty$.
Then $\alpha_{i_j} := \min\{ \beta_1,\beta_2 \}$.

 For the point $\mathbf{p}_e$, $\mathbf{x}(\alpha_e) \in$ conv$(S)$ implies $\alpha_e \leq 0$.  
 That is,  if $ \mathbf{x}_S \in$ conv$(S)$, then $\mathbf{x}(0) = \mathbf{x}_S$, so that $\alpha_e = 0$.
 If  $\mathbf{x}_S \notin $ conv$(S)$, then $\mathbf{x}(\alpha_e) \in$ aff$(S)$, where $\alpha_e = \frac{(\mathbf{p}_{i_1} - \mathbf{x}_S)\mathbf{d}_S}{\mathbf{d}_S^2} < 0$.

Thus, the intersection of $X_S$ and the facet of conv$(S \cup \{ \mathbf{p}_e \})$ first encountered along $X_S$ occurs 
at the step size $\alpha^*$, where
\begin{equation}
 \alpha^*  = \min_{\mathbf{p}_{i_j} \in S } \{ \alpha_{i_j}   \}.  \label{alp3}
\end{equation}
Choose $\mathbf{p}_{i_l}  \in S$ so that $\alpha^* = \alpha_{i_l}$.
Then the solution to the linear system (\ref{Teqpa}) with $\alpha = \alpha^*$  yields
$\pi_{i_l}(\alpha^*) = 0$, and $\pi_{i_j} (\alpha^*) \geq 0$, for $\mathbf{p}_{i_j} \in S \cup \{ \mathbf{p}_e \} \setminus \{ \mathbf{p}_{i_l} \}$ 
so that $\mathbf{x}(\alpha^*) \in F_{i_l}$.
This construction shows that  $S  \setminus \{ \mathbf{p}_{i_l} \}$ is an active set for $[\mathbf{x}(\alpha^*), z(\alpha^*)]$
and that  $[\mathbf{x}(\alpha^*), z(\alpha^*)]$ is dual feasible with respect to $S \setminus \{ \mathbf{p}_{i_l} \} \cup \{ \mathbf{p}_e \}$.
However, the constraint corresponding to $\mathbf{p}_e$ is not active at $[\mathbf{x}(\alpha^*), z(\alpha^*)]$.
 
 Delete the point $\mathbf{p}_{i_l}$ from $S$.  Reset $S: =  S \setminus \{ \mathbf{p}_{i_l} \}$, 
 and  $[\mathbf{x}_S, z_S] := [ \mathbf{x}(\alpha^*), z(\alpha^*) ]$.
The set $S$ is active with respect to the ball $[\mathbf{x}_S, z_{S}]$, and $[\mathbf{x}_S, z_{S}]$ is dual feasible with respect to 
$S \cup \{ \mathbf{p}_e \} $.  
The Search Phase is re-entered with the active set $S$, the entering point $\mathbf{p}_e$, and the ball $[\mathbf{x}_S, z_{S}]$.

If there is more than one  point
$\mathbf{p}_{i_j} \in S\setminus \{\mathbf{p}_{e}\}$ such that $\alpha^* = \alpha_{i_l}$, then  $X_S$
intersects the corresponding facets $F_{i_j}$ simultaneously, and there is a tie for the leaving point. 
In this case choose any point $\mathbf{p}_{i_j}$
 such that $\alpha_{i_j} = \alpha^*$, to be deleted from $S$. 
The set $S$ is reset to $S:=S \cup \{ \mathbf{p}_{i_j} \}$. 
The Search Phase is re-entered with the active set $S$,
the entering point $\mathbf{p}_e$, and the ball $[\mathbf{x}(\alpha^*), z(\alpha^*) ]$,
and the search continues on the facet $F_{i_j}$.
This case may lead to a degenerate iteration at  the next step
with $\alpha^* = \alpha_{i_l} =0$ for some $\mathbf{p}_{i_l}$.  
Cycling will not occur since at each degenerate iteration one point is deleted from the finite set $S$.  
After a finite number of points are deleted, 
 $S$ is reduced to two points, and the step size $\alpha$ will be positive at the next iteration.


\section{Dual Algorithm}

\begin{enumerate}
\item Input: a set of distinct points $P=\{{\mathbf{p}_1,...,\mathbf{p}_m}\} \subset \real^n$, \\
a radius $r_1, \geq 0$, a ball $[\mathbf{p}_i, r_i] = \{ \mathbf{x} : \parallel \mathbf{x} - \mathbf{p}_i \parallel \leq r_i \}$ for each $\mathbf{p}_i \in P$. \\
Output: a unique  ball $[\mathbf{x}^*, z^*]$, with minimum radius $z^*$, containing  $[\mathbf{p}_i, r_i]$ for $\mathbf{p}_i$ in $P$.\\
Assume condition (1) for each pair $ \mathbf{p}_{j}, \mathbf{p}_{k} \in P.$

\item Initialize: Choose  $S= \{ \mathbf{p}_{i_1}, \mathbf{p}_{i_2} \}$, for any two points $ \mathbf{p}_{i_1}, \mathbf{p}_{i_2} \in P$,  ordered $\ni$ $r_{i_1} \geq r_{i_2}$.\\
Compute the dual feasible solution $[\mathbf{x}_{S}, z_S]$ using (\ref{dinit}).

\item Optimality: $[\mathbf{x}_S, z_S]$ is optimal if  $\parallel \mathbf{x}_S - \mathbf{p}_{i} \parallel + r_{i} \leq z_S$, $\forall$ $\mathbf{p}_i \in P\setminus S$.  \\
Else, choose entering point  $\mathbf{p}_{e} \in P \setminus S$ $\ni$  $\parallel \mathbf{x}_S - \mathbf{p}_{e} \parallel + r_{e} > z_S$.

\item Update: If  (\ref{linsys1}), (\ref{linsys2}) has no solution,  go to Step 5
with $[\mathbf{x}_{S}, z_S]$, $S$, and $\mathbf{p}_{e}$. \\
Else, (\ref{linsys1}), (\ref{linsys2}) has a solution $\mathbf{\lambda}$.  Solve (\ref{denom1}), (\ref{denom2}) for $\pi$. \\ 
Determine a leaving point $\mathbf{p}_{i_l} \in S$  using minimum ratio test (\ref{minratio1}).\\
Reset $S:= S \setminus \{\mathbf{p}_{i_l}\}$.  Go to Step 5
with $[\mathbf{x}_{S}, z_S]$, its active set $S$, and entering point $\mathbf{p}_{e}$.
\item Search Path: Denote $S = \{ \mathbf{p}_{i_1}, \ldots, \mathbf{p}_{i_{s}} \}$, 
ordered so that  $r_{i_1} \geq \ldots \geq r_{i_{s}}$.
	 \begin{enumerate}  
	\item[5 a.]  If $r_{i_1} = r_{i_s}$, the search path  $X_S$ is given by  (\ref{dird}).  \\
	 Find the step size $\alpha'$ so that  $\mathbf{x}(\alpha') \in X_S \cap B_{i_1,e}$. 
		 \begin{enumerate}  
		 \item[] If $r_{i_1} = r_e$, determine step size $\alpha' $ using  (\ref{PEQR}).
		 \item[] If $r_{i_1} \neq r_e$, determine step size $\alpha'$ by solving the quadratic equation  (\ref{QUADe}).
		\end{enumerate}
	Solve  (\ref{Teqr})   for $\boldsymbol{\pi}(\alpha') = (\pi_{i_1}(\alpha'), \dots \pi_{i_s}(\alpha'),\pi_{e}(\alpha'))$.  
		\begin{enumerate}
		\item[] If $\boldsymbol{\pi}(\alpha')  \geq 0$,  $\mathbf{x}(\alpha') \in$ conv$(S \cup \{\mathbf{p}_e\})$.  Go to Step 3 with $[\mathbf{x}_S, z_S]$ and $S$.
		\item[] If $\pi_{i_j}(\alpha') < 0$, for some $\mathbf{p}_{i_j} \in S$, find  $\alpha_{i_j}$ using (\ref{Teqr})  $\ni$ $\mathbf{x}(\alpha_{i_j}) \in F_{i_j}$, $\forall \mathbf{p}_{i_j} \in S$.\\
		Choose $\alpha^* = \min \{ \alpha_{i_j} : \mathbf{p}_{i_j} \in S \}$, and
		choose $\mathbf{p}_{i_l} \in S$ so that $\alpha^* = \alpha_{i_l}$. \\
		Reset $S:=S \setminus \{ \mathbf{p}_{i_l} \}$, and  $[\mathbf{x}_S, z_S] := [\mathbf{x}(\alpha^*), z(\alpha^*)]$.  Go to Step 5.
		\end{enumerate}
	\item[5 b.] If $r_{i_1} > r_{i_{s}}$, determine $\mathbf{v}_S$, $\mathbf{c}_S$, $a_S$, $b_S$,  $\epsilon_S$ of $B_S$ using  (\ref{HT12})  - (\ref{pnc}).
		\begin{enumerate}
		\item[5 b 1:]  If $\epsilon_S >1$, $B_S$ is a hyperboloid of dimension $n-s+2$.
			\begin{enumerate}
			\item[] If $\mathbf{x}_S = \mathbf{a}_S$, then  $\beta_S = 0$; else compute $\beta_S$.\\
			Construct the  hyperbola $Y_S$ using  (\ref{ysdualhy})
			where $\mathbf{u}_S$ is determined by (\ref{usdhyp}).\\
			Construct the search path $X_S$ using  (\ref{xsdualhy})\\
			Find the step size $\beta'$ using (\ref{he}), (\ref{AAAA}), (\ref{hysold}).  \\
			Then $\alpha' = \beta' - \beta_S$.
			\end{enumerate}
		\item[] Solve  (\ref{Tuneqh})   for $\boldsymbol{\pi(\alpha')} =  (\pi_{i_1}(\alpha'), \dots \pi_{i_s}(\alpha'),\pi_{e}(\alpha'))$.  
			\begin{enumerate}
			\item[] If $\boldsymbol{\pi(\alpha')} \geq 0$,  $\mathbf{x}(\alpha') \in$ conv$(S \cup \{\mathbf{p}_e\})$.  Go to Step 3 with $[\mathbf{x}_S, z_S]$ and $S$.
			\item[] If $\pi_{i_j}(\alpha') < 0$, for some $\mathbf{p}_{i_j} \in S$,  find  $\alpha_{i_j}$ using (\ref{Tevu}) $\ni$ $\mathbf{x}(\alpha_{i_j}) \in F_{i_j}$ $\forall \; \mathbf{p}_{i_j} \in S$.\\
			Choose $\alpha^* = \min \{ \alpha_{i_j} : \mathbf{p}_{i_j} \in S \}$, and
			choose $\mathbf{p}_{i_l} \in S$ so that $\alpha^* = \alpha_{i_l}$. \\
			Reset $S:=S \setminus \{ \mathbf{p}_{i_l} \}$, and  $[\mathbf{x}_S, z_S] = [\mathbf{x}(\alpha^*), z(\alpha^*)]$. Go to Step 5.
			 \end{enumerate}
		\item[5 b 2:]  If $\epsilon_S < 1$, $B_S$ is an ellipsoid of dimension $n-s+2$.
			\begin{enumerate}
			\item[] If $\mathbf{x}_S = \mathbf{a}_S$, then  $\beta_S = 0$; else compute $\beta_S$.\\
			Construct $\mathbf{v}_S$, $\mathbf{u}_S$ and the  ellipse $Y_S$ using  (\ref{ysduale}) and (\ref{uella})\\
			Construct the search path $X_S$ using (\ref{xsduale}).\\
			Find the step size $\beta'$ using  (\ref{he}), (\ref{AAAE}) (\ref{elsola}).  \\
			Then $\alpha' = \beta' - \beta_S$.
			\end{enumerate}
		\item[] Solve  (\ref{Teqhe})   for $\boldsymbol{\pi(\alpha')} =  (\pi_{i_1}(\alpha'), \dots \pi_{i_s}(\alpha'),\pi_{e}(\alpha'))$.  
			\begin{enumerate}
			\item[] If $\boldsymbol{\pi(\alpha')} \geq 0$, then $\mathbf{x}(\alpha') \in$ conv$(S \cup \{\mathbf{p}_e\})$.  Go to Step 3 with $[\mathbf{x}_S, z_S]$ and $S$ 
			\item[] If $\pi_{i_j}(\alpha') < 0$, for some $\mathbf{p}_{i_j} \in S$, find  $\alpha_{i_j}$ using (\ref{elsola}) $\ni$ $\mathbf{x}(\alpha_{i_j}) \in F_{i_j}$,$\forall \mathbf{p}_{i_j} \in S$.\\
			Choose $\alpha^* = \min \{ \alpha_{i_j} : \mathbf{p}_{i_j} \in S \}$, and
			choose $\mathbf{p}_{i_l} \in S$ so that $\alpha^* = \alpha_{i_l}$. \\
			Reset $S:=S \setminus \{ \mathbf{p}_{i_l} \}$, and  $[\mathbf{x}_S, z_S] = [\mathbf{x}(\alpha^*), z(\alpha^*)]$. Go to Step 5.
			 \end{enumerate}
		 \item[5 b 3:]  If $\epsilon_S = 1$, $B_S$ is a paraboloid of dimension $n-s+2$.
			\begin{enumerate}
			\item[] If $\mathbf{x}_S = \mathbf{a}_S$, then  $\beta_S = 0$; else compute $\beta_S$.\\
			Construct the parabola and $\mathbf{u}_S$  using  (\ref{ysdualp}) and (\ref{uell}) \\
			Construct the search path $X_S$ using  (\ref{xsdualp}).\\
			Find the step size $\beta'$ using  (\ref{AADP}).  \\
			Then $\alpha' = \beta' - \beta_S$.
			\end{enumerate}
		\item[] Solve  (\ref{Teqpa})   for $\boldsymbol{\pi(\alpha')} =  (\pi_{i_1}(\alpha'), \dots \pi_{i_s}(\alpha'),\pi_{e}(\alpha'))$.  
			\begin{enumerate}
			\item[] If $\boldsymbol{\pi(\alpha')} \geq 0$, then $\mathbf{x}(\alpha') \in$ conv$(S \cup \{\mathbf{p}_e\})$.  Go to Step 3 with $[\mathbf{x}_S, z_S]$ and $S$ 
			\item[] If $\pi_{i_j}(\alpha') < 0$, for some $\mathbf{p}_{i_j} \in S$, find  $\alpha_{i_j}$ using  (\ref{Teqpa})  $\ni$ $\mathbf{x}(\alpha_{i_j}) \in F_{i_j}, \forall  \mathbf{p}_{i_j} \in S$.\\
			Choose $\alpha^* = \min \{ \alpha_{i_j} : \mathbf{p}_{i_j} \in S \}$, and
			choose $\mathbf{p}_{i_l} \in S$ so that $\alpha^* = \alpha_{i_l}$. \\
			Reset $S:=S \setminus \{ \mathbf{p}_{i_l} \}$, and  $[\mathbf{x}_S, z_S] := [\mathbf{x}(\alpha^*), z(\alpha^*)]$.  Go to Step 5.
			Go to 5.
			 \end{enumerate}
		\end{enumerate}
	\end{enumerate}
\end{enumerate}


\section{Computational Results}

Computational results for the primal and dual algorithms are forthcoming.

\appendix
\section{Conic Sections in $\real^n$}
\renewcommand{\theequation}{\thesection.\arabic{equation}}

This Appendix presents  results from reference \cite{Dearing1} that are used to find the intersection 
of a hyperplane with each of the conic sections in $\real^n$.
These results are used to develop the search paths of the primal and dual algorithms for problem $M(P)$.
Recall from Section 2, that a pair of points $\mathbf{p}_{j}, \mathbf{p}_{k} \in P$,
with  radii $r_{j} > r_{k}$, defines the \textbf{bisector} 
\begin{equation}
B_{j,k} = \{ \mathbf{x}: \parallel \mathbf{p}_{k} - \mathbf{x} \parallel - \parallel \mathbf{p}_{j} - \mathbf{x} \parallel =  r_j  - r_k  \} \label{A1a}
\end{equation}
which is one sheet of the hyperboloid $H\!B_{j,k}$ in $\real^n$ with two sheets, where
\begin{equation}
H\!B_{j,k} = \{ \mathbf{x}: |\parallel \mathbf{p}_{k} - \mathbf{x} \parallel - \parallel \mathbf{p}_{j} - \mathbf{x} \parallel| =  r_j  - r_k  \}.\label{A1b}
\end{equation}
The other sheet of $H\!B_{j,k}$ is given by 
\begin{equation}
\hat{B}_{j,k} = \{ \mathbf{x}: \parallel \mathbf{p}_{j} - \mathbf{x} \parallel - \parallel \mathbf{p}_{k} - \mathbf{x} \parallel =  r_j  - r_k  \},\label{A1c}
\end{equation}
so that  $H\!B_{j,k} = B_{j,k} \cup \hat{B}_{j,k}$.
The hyperboloid $H\!B_{j,k}$ is specified by the \textbf{focal points} $\mathbf{p}_{j}$ and $\mathbf{p}_{k}$,
and the radii  $r_{j}$ and $ r_{k}$ with $r_{j} > r_{k}$.
The \textbf{center}   of $H\!B_{j,k}$ is  the center  point of the line segment  between   $\mathbf{p}_{j}$ and  $ \mathbf{p}_{k}$,  
given by $\mathbf{c}_{j,k} =\frac{1}{2}(\mathbf{p}_{j} + \mathbf{p}_{k})$, and the distance from $\mathbf{c}_{j,k} $ to $\mathbf{p}_{j} $, or to  $\mathbf{p}_{k} $, 
is  given by $c_{j,k} =\frac{1}{2}\!\parallel\!\!\mathbf{p}_{j} - \mathbf{p}_{k}\!\!\parallel$.
The unit vector $\mathbf{v}_{j,k} = \frac{\mathbf{p}_{j} - \mathbf{p}_{k}}{\parallel \mathbf{p}_{j} - \mathbf{p}_{k} \parallel }$ 
 from $\mathbf{p}_{k}$ to  $ \mathbf{p}_{j}$ is called the \textbf{axis vector}, and is parallel to the \textbf{major axis},
 which is the line through  $ \mathbf{p}_{j}$ and $\mathbf{p}_{k}$. 
 The center $\mathbf{c}_{j,k}$ locates the hyperboloid and the axis vector $\mathbf{v}_{j,k}$ specifies its orientation.
 
Under the assumption that $r_{j} > r_{k}$, the  point of intersection between  the major axis and the sheet $B_{j,k}$
occurs at the $ \textbf{vertex} $  $\mathbf{a}_{j,k}$, where $\mathbf{a}_{j,k} = \mathbf{c}_{j,k} + a_{j,k} \mathbf{v}_{j,k}$, with  $a_{j,k} = \frac{1}{2} (r_j - r_k)$.  
The vertex
$\mathbf{a}_{j,k}$ is the center of the smallest ball containing, or  internally tangent to,  
the two balls  $[\mathbf{p}_j, r_j]$ and $[\mathbf{p}_k, r_k]$. 
Each point $\mathbf{x}$ on $B_{j,k}$ is the center of a ball $[\mathbf{x},z_{\mathbf{x}}]$, with radius $z_{\mathbf{x}} = \parallel \mathbf{p}_j -\mathbf{x} \parallel + r_j = \parallel \mathbf{p}_k -\mathbf{x} \parallel + r_k$, that contains, and is internally tangent to, the two balls $[\mathbf{p}_j, r_j]$ and $[\mathbf{p}_k, r_k]$. 
The sheet $\hat{B}_{j,k}$ intersects the major axis at the vertex $\hat{\mathbf{a}}_{j,k} = \mathbf{c}_{j,k} - a_{j,k} \mathbf{v}_{j,k}$,
which is the center of the smallest ball externally tangent to the two balls $[\mathbf{p}_j, r_j]$ and $[\mathbf{p}_k, r_k]$. 
Problem $M(P)$ is concerned only with the sheet $B_{j,k}$.

The \textbf{eccentricity} specifies the shape of the hyperboloid $H\!B_{j,k}$ and is given by  $\epsilon_{j,k} = \frac{c_{j,k}}{a_{j,k}}$.
Assumption (1) implies that $0 < a_{j,k} < c_{j,k}$, so that $\epsilon_{j,k} >1$, as required for a hyperboloid.
The \textbf{directrix} of the sheet $B_{j,k}$ is the hyperplane with normal vector $\mathbf{v}_{j,k}$  
containing  the point $\mathbf{d}_{j,k}$, where
where $\mathbf{d}_{j,k} = \mathbf{c}_{j,k} + d_{j,k} \mathbf{v}_{j,k}$,  and 
$d_{j,k} = \frac{a_{j,k}^2}{c_{j,k}} = \frac{(r_j - r_k)^2}{2 \parallel  \mathbf{p}_{j}-\mathbf{p}_{k} \parallel}$.
The directrix is given by  $\mathbf{v}_{j,k} \mathbf{x} =  \mathbf{v}_{j,k} \mathbf{d}_{j,k} 
= \mathbf{v}_{j,k}  \mathbf{c}_{j,k} + d_{j,k}
= \frac{ \parallel \mathbf{p}_{j} \parallel^2 -  \parallel \mathbf{p}_{k} \parallel^2 +  (r_j - r_k)^2} {2  \parallel  \mathbf{p}_{j}-\mathbf{p}_{k} \parallel}$. 
The directrix of the sheet $\hat{B}_{j,k}$ is the hyperplane with normal vector $\mathbf{v}_{j,k}$ 
through the point $\hat{\mathbf{d}}_{j,k} = \mathbf{c}_{j,k} - d_{j,k} \mathbf{v}_{j,k}$.

An $n$-dimensional \textbf{ellipsoid}, $E_{j,k}$, is specified by the same vectors and parameters that specify a hyperboloid, all of which are determined by the 
focal points $\mathbf{p}_j$ and $\mathbf{p}_k$ and the positive constant $a_{j,k}$, where 
\begin{equation}
E_{j,k} = \{ \mathbf{x}: \| \mathbf{p}_j - \mathbf{x} \| +  \| \mathbf{p}_j - \mathbf{x} \| = 2a_{j,k} \}. \notag
\end{equation}
The axis vector $\mathbf{v}_{j,k}$, the center $\mathbf{c}_{j,k}$, the parameter $c_{j,k}$, vertices $\mathbf{a}_{j,k}, \hat{\mathbf{a}}_{j,k}$, eccentricity 
$\epsilon_{j,k}$, and directrix $\mathbf{d}_{j,k}, \hat{\mathbf{d}}_{j,k}$, each have the same definition for an ellipsoid as for a hyperbolod, except 
that $c_{j,k} < a_{j,k}$, so that $\epsilon_{j,k} = c_{j,k}/a_{j,k} <1$.
The next property gives a representation of a hyperboloid, or an ellipsoid, as a quadratic form.
 \begin{property} Given the focal points $\mathbf{p}_1$ and $\mathbf{p}_2$, with corresponding axis vector $\mathbf{v}$, center point $\mathbf{c}$,
 positive constants $a$ and $c$, and eccentricity $\epsilon = c/a$, let $H$ be the hyperboloid determined by these vectors and parameters if $\epsilon > 1$, 
 and let $E$ be the ellipsoid determined by these vectors and parameters if $\epsilon < 1$.  Define: 
\begin{equation}
Q = \{ \mathbf{x}: (\mathbf{x} - \mathbf{c})^T[I - \epsilon^2\mathbf{v}\mathbf{v}^T ](\mathbf{x} - \mathbf{c}) = a^2 - c^2 \}
= \{ \mathbf{x}:  [\mathbf{x} - \mathbf{c}]^2 - [\epsilon(\mathbf{x} - \mathbf{c})\mathbf{v}]^2  = a^2 - c^2 \}  \label{Q} 
\end{equation}
\begin{equation}
\text{and \;\;\;}C = \{ \mathbf{x}: (\mathbf{x} - \mathbf{c})^T[I - \epsilon^2\mathbf{v}\mathbf{v}^T ](\mathbf{x} - \mathbf{c}) = 0 \}
= \{ \mathbf{x}:  [\mathbf{x} - \mathbf{c}]^2 - [\epsilon(\mathbf{x} - \mathbf{c})\mathbf{v}]^2  = 0 \}. \notag
\end{equation}
If   $c > a$, then $Q=H$, and if $c < a$, then $Q = E$. 
If $\epsilon > 1$,  then $C$ is the  cone that is the asymptotic approximation to the hyperboloid $H$.
\end{property}

An $n$-dimensional\textbf{ paraboloid}, symmetric about its major axis, is specified by two points $\mathbf{p}_1$ and $\mathbf{p}_2$ only, 
as the set of all points $\mathbf{x} \in \real^n$ whose distance from $\mathbf{p}_1$ equals the distance from $\mathbf{x}$ to the hyperplane that is orthogonal
to the major axis and contains the point $\mathbf{p}_2$. The axis vector is $\mathbf{v} = \frac{(\mathbf{p}_1 - \mathbf{p}_2)} {\|\mathbf{p}_1 - \mathbf{p}_2\|}$. 
The paraboloid is defined as the set 
\begin{equation}
P = \{ \mathbf{x}: \| \mathbf{p}_1 - \mathbf{x} \| = \mathbf{v}(\mathbf{x} - \mathbf{p}_2) \}.
\end{equation}
The center of the paraboloid, also its vertex, is $\mathbf{c} = \frac{\mathbf{p}_1 + \mathbf{p}_2)}{2}$, and the parameter $c = \frac{\| \mathbf{p}_1 - \mathbf{p}_2 \|}{2}$.
\begin{property}  The paraboloid $P$ has the equivalent expression
\begin{equation}
P = \{ \mathbf{x} : (\mathbf{x} - \mathbf{c})^T[I - \mathbf{v}\mathbf{v}^T](\mathbf{x} - \mathbf{c}) = 4c\mathbf{v}((\mathbf{x} - \mathbf{c}) \}
\end{equation}
\end{property}

  From the classical studies of conic sections in $\real^3$, it is well known that if a  plane  and a cone  
  intersect at an appropriate angle, measured between the axis vector and the hyperplane, the intersection is a hyperbola of dimension 2.
  The following development extends these results to $\real^n$ and considers the intersection of a hyperplane with a hyperboloid, an ellipsoid, and a paraboloid.
  Conditions are given for the resulting intersection to be a hyperboloid, an ellipsoid or a paraboloid of dimension $n-1$, 
  and expressions are derived for the parameters of the resulting structure. 
  These results are then applied to the conic sections arising in problem $M(P)$.
\begin{property}  Suppose $Q$ 
 is a hyperboloid in $\real^n$,  centered at $\mathbf{c}$, 
with  axis  vector $\mathbf{v}$ of unit length,  eccentricity $\epsilon$, and 
  parameters $a < c$, and suppose 
  $H\!P = \{ \mathbf{x}: \mathbf{h} \mathbf{x}  =   \mathbf{h}(\mathbf{c} + \hat{h} \mathbf{h} \}$,
  is a hyperplane with normal vector  $\| \mathbf{h} \| =1$ through the point $\mathbf{c} + \hat{h} \mathbf{h}$.
  Let $\rho = \sqrt{1-(\mathbf{h} \mathbf{v})^2}$.  
 Then  $Q \cap H\!P$  is a non-empty hyperboloid of dimension $n-1$ iff  $\epsilon \rho > 1$,  a non-empty ellipsoid of dimension $n-1$ iff $\epsilon \rho < 1$ 
 and $\sqrt{a^2(1 - \epsilon^2\rho^2)}\leq \hat{h}$, or a paraboloid of dimension $n-1$ iff  $\epsilon \rho = 1$. 
 The expression $\sqrt{a^2(1 - \epsilon^2\rho^2)}$
 is the minimum distance between  the hyperplane $H\!P$ passing through the point $\mathbf{c}$ (with $\hat{h} = 0$), 
 and the hyperplane $H\!P$ passing through a point on $Q$. 
  \end{property}
  
  \begin{property}  Suppose $Q$ 
 is an ellipsoid in $\real^n$,  centered at $\mathbf{c}$, 
with  axis  vector $\mathbf{v}$ of unit length,  eccentricity $\epsilon$, and 
  parameters $a > c$, and suppose 
  $H\!P = \{ \mathbf{x}: \mathbf{h} \mathbf{x}  =   \mathbf{h}(\mathbf{c} + \hat{h} \mathbf{h} \}$,
  is a hyperplane with normal vector  $\| \mathbf{h} \| =1$ through the point $\mathbf{c} + \hat{h} \mathbf{h}$.
  Let $\rho = \sqrt{1-(\mathbf{h} \mathbf{v})^2}$.  
 Then  $Q \cap H\!P$  is a non-empty elipsoid of dimension $n-1$ iff  
 $\sqrt{a^2(1 - \epsilon^2\rho^2)} \geq \hat{h}$. 
  \end{property}
  
  \begin{property} Suppose $P = \{ \mathbf{x}: \| \mathbf{p}_1 - \mathbf{x} \| = \mathbf{v}(\mathbf{x} - \mathbf{p}_2) \}$
  is a paraboloid centered at $\mathbf{c}$ with axis vector $\mathbf{v}$ and parameter $c$, 
  and suppose that $HP = \{ \mathbf{x} : \mathbf{h}(\mathbf{x} - \mathbf{c}) = \hat{h} \}$, is a hyperplane with $\|\mathbf{h} \| = 1$, and 
  $P \cap HP \ne \emptyset$.   Let $\rho = \sqrt{1 - (\mathbf{h}\mathbf{v})^2}$.  Then $P \cap HP$ is a paraboloid of dimension $n-1$ if $\rho = 1$ or an ellipsoid of dimension $n-1$ if $\rho < 1$.
  
  \end{property}
  
  \begin{figure}
\begin{center}
\includegraphics[scale=.8]{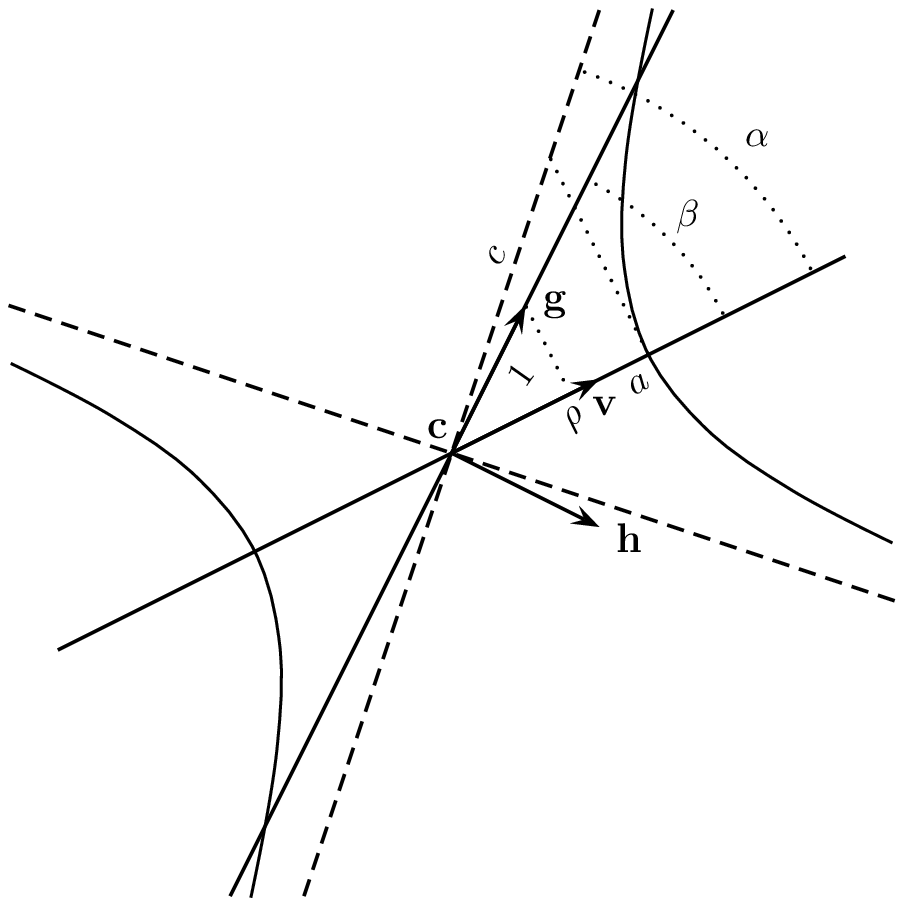}
\caption{Hyperbola $H$ projected onto the affine plane aff$(\mathbf{v}, \mathbf{h}, \mathbf{c})$}
\end{center}
\end{figure}

The relationship between the condition $\epsilon \rho > 1$ and the intersection of the hyperboloid $H$ and the hyperplane $H\!P$ are illustrated by Figure 1,
which shows the vectors $\mathbf{v}$ and $\mathbf{h}$, and the point $\mathbf{c}$, in the affine plane aff$(\mathbf{v}, \mathbf{h}, \mathbf{c})$.
  The projection of the hyperboloid $H$ onto aff$(\mathbf{v}, \mathbf{h}, \mathbf{c})$ is shown by the  hyperbola of two sheets, with parameters $a$ and $c$.
   The projection of the hyperplane $H\!P$ onto  aff$(\mathbf{v}, \mathbf{h}, \mathbf{c})$  is some line parallel to the ray 
   determined by the vector $\mathbf{g}$ that is orthogonal to $\mathbf{h}$.
  Let  $\beta$ be the angle between $\mathbf{g}$ and $\mathbf{v}$, which is the same angle between  the projection of $H\!P$ and $\mathbf{v}$.
  Then $\rho = \mathbf{g}\mathbf{v} = \cos(\beta)$.
  
  Let $C$ be the cone that is the asymptotic approximation to $H$, so that 
   the projection of $C$ onto aff$(\mathbf{v}, \mathbf{h}, \mathbf{c})$ is illustrated by the dashed lines through $\mathbf{c}$, that is, the 
  asymptotes of the hyperbola. 
  The angle between either asymptote and $\mathbf{v}$ is $\alpha$.
  
  Then, $\epsilon \rho > 1$ if and only if $\rho > a/c$ if and only if $\cos(\beta) > \cos(\alpha)$ if and only if $\beta < \alpha$ for $0 \leq \beta \leq \pi/2$,
  if and only if $H\!P \cap H \neq \emptyset$.  That is, $H\!P$ intersects both sheets of $H$.  
  Alternatively,  $\epsilon \rho< 1$ if and  only if $H\!P$ intersects one sheet of $C$, or the point $\mathbf{c}$, 
  If $\epsilon \rho< 1$, the distance squared between $H\!P$ and $\mathbf{c}$ is $\hat{\mathbf{h}}^2$, and 
  the minimum distance squared between $\mathbf{c}$ and  $H\!P \cap Q$  is $a^2(1 - \epsilon^2\rho^2)$,
  so that the condition $a^2(1 - \epsilon^2\rho^2) \leq \hat{h}^2$ implies $H\!P$ intersects one sheet of $Q$.

The following corollary gives expressions that are used to compute the vectors and parameters of the $n-1$ dimensional  conic section resulting from the intersection of a hyperplane and an $n$-dimensional hyperboloid.

\begin{corollary} If $H\!B \cap H\!P$ is an hyperboloid or an ellipsoid, then its vectors are given as follows: axis vector   $\mathbf{g}_1 = (\mathbf{v}-(\mathbf{v}\mathbf{h})\mathbf{h})/\| \mathbf{v}-(\mathbf{v}\mathbf{h})\mathbf{h} \|$ 
(the component of $\mathbf{v}$ that is orthogonal to the projection of $\mathbf{v}$  onto $\mathbf{h}$), 
center $\hat{\mathbf{c}} = \mathbf{c} +\hat{h}\mathbf{h} +\tilde{c}\mathbf{g}_1$, verex $\hat{\mathbf{a}} = \hat{\mathbf{c}} + \hat{a}\mathbf{g}_1$ (for a hyperboloid) or  $\hat{\mathbf{a}} = \hat{\mathbf{c}} \pm \hat{a}\mathbf{g}_1$ (for an ellipse) , focal points $\hat{\mathbf{c}} \pm \hat{a}\epsilon \rho \mathbf{g}_1$, where 
$\hat{a} = \frac{(1 - \epsilon^2)(a^2(1-\epsilon^2\rho^2)-\hat{h}^2)}{(1-\epsilon^2\rho^2)^2}$,  $\tilde{c} = \frac{\epsilon^2\rho\sigma\hat{h}}{(1-\epsilon^2\rho^2)^2}$.
If $H\!B \cap H\!P$ is a paraboloid, then its axis vector and center is the same as above.  Its focal points are $\hat{\mathbf{c}} \pm \bar{c} \mathbf{g}_1$, with parameters $\tilde{c} = \frac{\epsilon\sigma \hat{h}}{2}$,  $\hat{c} = \frac{(\epsilon^2\rho^2-1)\hat{h}^2 - (a^2-c^2)}{4 \tilde{c}}$.
\end{corollary}
   
The class of hyperboloids $H\!B_{j,k}$ that occur in the problem $M(P)$ have additional properties.
Let  $T = \{ \mathbf{p}_{j}, \mathbf{p}_{k}, \mathbf{p}_{l} \}$ be a set of three affinely independent
points from $P$ ordered so that $r_{j} \geq r_{k} \geq r_{l}$, and let $\mathcal{B} = \{ B_{j,k}, B_{j,l} , B_{k,l} \}$
denote the bisectors corresponding to the respective pairs of points from $T$.
 If the radii are unequal, that is  $r_{j} > r_{l} $,
then at least two pairs of points from $T$ have unequal radii, 
and the bisector corresponding to each of these pairs is one sheet of a hyperboloid. 

The following theorem constructs a unique hyperplane $H_T$  that contains the intersection of  any two 
bisectors in $\cal{B}$.  The theorem also shows that for any two bisectors in $\cal{B}$, say $B_{j,k}$ and $B_{j,l}$, then
$ B_{j,k} \cap B_{j,l} =  B_{j,k}\cap H_T = B_{j,l}\cap H_T$.
Also, if one of the bisectors in $\cal{B}$ is a hyperplane, then it is identical to $H_T$.
This result allows the intersection of two bisectors to be determined by the intersection of either one of the bisectors with the hyperplane $H_T$.

In the subsequent development, a bisector corresponding to two points with unequal radii will be referred to as a hyperboloid, even though the bisector is only one sheet of the hyperboloid generated by the two points and their radii.

\begin{theorem}
Suppose that   $T = \{ \mathbf{p}_{j}, \mathbf{p}_{k}, \mathbf{p}_{l} \}$ is a subset of three affinely independent points from $P$,
   ordered so that $r_{j} \geq r_{k} \geq r_{l}$, with $r_{j} > r_{l}$,
   and suppose that the intersection of the three bisectors from $\mathcal{B} = \{ B_{j,k}, B_{j,l} , B_{k,l} \}$ is nonempty.  
  Let   \begin{equation}
H_T = \{ \mathbf{x}: \mathbf{h}_T \mathbf{x}  =  \mathbf{h}_T \mathbf{d}_T \} =
 \{ \mathbf{x}: (\epsilon_{j,k} \mathbf{v}_{j,k} - \epsilon_{j,l} \mathbf{v}_{j,l} ) \mathbf{x}  = (\epsilon_{j,k} \mathbf{v}_{j,k} - \epsilon_{j,l} \mathbf{v}_{j,l}) \mathbf{d}_T \},
\label{Ajkjl} \end{equation}
   where  $\mathbf{d}_T$ is  the intersection of aff$(T)$,  the directrix $\mathbf{v}_{j,l}\mathbf{x} =  \mathbf{v}_{j,l}\mathbf{d}_{j,l}$ of $B_{j,l}$, and the directrix  $ \mathbf{v}_{j,k}\mathbf{x} =  \mathbf{v}_{j,k}\mathbf{d}_{j,k}$ of $B_{j,k}$.
If any one of the bisectors in $\cal{B}$ is a hyperplane, it is identical to $H_T$.
For each pair of bisectors in $\cal{B}$, their intersection equals the intersection of $H_T$ with either bisector in the pair.
\end{theorem}
 
 The next property shows how to compute the intersection point $\mathbf{d}_T$ under the assumptions of Theorem A.1.
 \begin{property}
Given a set $T =  \{ \mathbf{p}_{j}, \mathbf{p}_{k}, \mathbf{p}_{l} \}$ 
of affinely independent points from $P$,
so that $r_{j} \geq r_{k} \geq r_{l}$, with $r_{j} > r_{l}$,
and the hyperplane
$H_T = \{ \mathbf{x}: \mathbf{h}_T \mathbf{x}  = \mathbf{h}_T \mathbf{d}_T \},$
  containing the intersections of $ B_{j,k}, B_{j,l} $, and $ B_{k,l} $, then
  $$ \mathbf{d}_T = \mathbf{d}_{jl} + \frac{\mathbf{v}_{jk}(\mathbf{d}_{jk}-\mathbf{d}_{jl})}{\mathbf{v}_{jk}\mathbf{u}_T} \mathbf{u}_T, 
\;\;\;\text{where}\;\;\; \mathbf{u}_T = \frac{\mathbf{h}_T - (\mathbf{v}_{jl} \mathbf{h}_T) \mathbf{v}_{jl}}
{\parallel \mathbf{h}_T - (\mathbf{v}_{jl} \mathbf{h}_T) \mathbf{v}_{jl}\parallel}$$ 
is a unit vector  in aff$(T)$ that is orthogonal to $\mathbf{v}_{jk}$.  
\end{property}  


The following development applies the above results to problem $M(P)$
for the case where at least two points in the active set $S$ have unequal radii.
Order the points in $S$ by non-increasing radii, so that $S = \{ \mathbf{p}_{i_1}, \ldots , \mathbf{p}_{i_{s}} \}$, with 
$r_{i_1} \geq \ldots \geq r_{i_s}$.   
By the assumption of unequal radii, $r_{i_1} > r_{i_s}$, so that $B_{i_1,i_s}$ is a hyperboloid.. 
 
 In this case, determining the search path requires computing the parameters and vectors of a two-dimensional conic section that is a subset
  of $B_S = \cap_{\{\mathbf{p}_{i_j}, \mathbf{p}_{i_k}\} \subset S} B_{i_j,i_k}$,
the intersection of all bisectors of pairs of points in $S$.
  The following property gives an equivalent expression for $B_{S}$, as the intersection of the $s-1$ bisectors $B_{i_1,i_j}$ for $j = 2,\ldots,s$.
  \begin{property}
   $B_{S} = \cap_{j=2}^{s} B_{i_1,i_j}$.
   \end{property}
   
   The next property gives an equivalent expression for $B_S$ as the intersection of the hyperboloid $B_{i_1,i_s}$ with a sequence of $s-2$
   hyperplanes $H_j$ for $j = 2, \ldots,s-1$.  
   \begin{property}
   $B_S = B_{i_1,i_s} \cap_{j=2}^{s-1}H_j$, where $H_j = \{ \mathbf{x}: \mathbf{h}_j\mathbf{x} = \mathbf{h}_j\mathbf{d}_j \}$ and is constructed from each triple 
   of points $T_j = \{ \mathbf{p}_{i_1},  \mathbf{p}_{i_j},  \mathbf{p}_{i_s} \}$, for $j = 2, \ldots,s-1$,  as specified in Theorem A.1 and Property A.4.
   \end{property}
   
  The vectors and parameters of $B_S = B_{i_1,i_{s}} \cap_{j=2}^{{s}-1} H_j$ are computed sequentially. 
  Initially, the vectors and parameters of $B_{i_1,i_s}$ are computed and designated as: $\mathbf{v}_{1}:=\mathbf{v}_{i_1,i_s}$, 
$\mathbf{c}_{1}:=\mathbf{c}_{i_1,i_s}$, $\mathbf{d}_{1}:=\mathbf{d}_{i_1,i_s}$, $\mathbf{a}_{1}:=\mathbf{a}_{i_1,i_s}$, 
$\epsilon_{1}:=\epsilon_{i_1,i_s}$, $a_{1}:=a_{i_1,i_s}$,  $b_{1}:=b_{i_1,i_s}$,  $c_{1}:=c_{i_1,i_s}$.
Also, define $\mathbf{hp}_1 = \mathbf{0}$.

Then for $k = 2,\ldots,s-1$, given the vectors and parameters $\mathbf{v}_{k-1}$, 
$\mathbf{c}_{k-1}$, $\mathbf{d}_{k-1}$, $\mathbf{a}_{k-1}$, 
$\epsilon_{k-1}$, $a_{k-1}$,  $b_{k-1}$,  $c_{k-1}$,
of $B_{i_1,i_s}  \cap_{j=2}^{k-1} H_j$,
the expressions \eqref{HT12} through \eqref{pa} compute the vectors  
$\mathbf{v}_{k}$, 
$\mathbf{c}_{k}$, $\mathbf{d}_{k}$, $\mathbf{a}_{k}$, and parameters
$\epsilon_{k}$, $a_{k}$,  $b_{k}$,  $c_{k}$,
of  $B_{i_1,i_s} \cap_{j=2}^{k-1} H_j \cap H_k$.
For the interations $k=3,\ldots, s-1$, $H_k$ must be projected onto $\cap_{j=2}^{k-1}H_j$.
    
\begin{alignat}{3}
\text{If\;} &r_{i_1} > r_{i_k} \text{\;\;\;\;\;\;} \mathbf{h}_{k} =   (\epsilon_{i_1,i_k}\mathbf{v}_{i_1,i_k} -  \epsilon_{i_1,i_{s}}\mathbf{v}_{i_1,i_{s}})/
\parallel \epsilon_{i_1,i_k}\mathbf{v}_{i_1,i_k} -  \epsilon_{i_1,i_{s}}\mathbf{v}_{i_1,i_{s}}\parallel    \label{HT12} \\
\text{If\;} &r_{i_1} = r_{i_k} \text{\;\;\;\;\;\;} \mathbf{h}_{k}  =  (\mathbf{p}_{i_1} -  \mathbf{p}_{i_k})/
\parallel \mathbf{p}_{i_1} -  \mathbf{p}_{i_k} \parallel   \label{HT2}\\
\mathbf{hp}_{k} &  = ( \mathbf{h}_{k} - \sum_{j=2}^{k-1}(\mathbf{h}_{k}  \mathbf{hp}_j) \mathbf{hp}_j )/
\parallel  \mathbf{h}_{k} - \sum_{j=2}^{k-1}(\mathbf{h}_{k}  \mathbf{hp}_j) \mathbf{hp}_j \parallel & \label{ph}\\
\mathbf{u}_{k-1} & = (\mathbf{hp}_{k} - (\mathbf{hp}_{k} \mathbf{v}_{k-1})\mathbf{v}_{k-1}) / 
\parallel \mathbf{hp}_{k} - (\mathbf{hp}_{k} \mathbf{v}_{k-1})\mathbf{v}_{k-1}\parallel   \label{pu}\\
\mathbf{v}_{k} & =  (\mathbf{v}_{k-1} - (\mathbf{v}_{k-1} \mathbf{hp}_{k})\mathbf{hp}_{k}) / 
\parallel \mathbf{v}_{k-1} - (\mathbf{v}_{k-1} \mathbf{hp}_{k})\mathbf{hp}_{k}  \parallel   \label{pvv}\\
\mathbf{d}_{k} & =  \mathbf{d}_{k-1} + \frac{\mathbf{v}_{i_1,i_k} (\mathbf{d}_{i_1,i_k} - \mathbf{d}_{k-1})}{\mathbf{v}_{i_1,i_k} \mathbf{u}_{k-1}} \mathbf{u}_{k-1}   \label{pd1}\\
\hat{h}_{k} & = \mathbf{hp}_{k} (\mathbf{d}_{k} - \mathbf{c}_{k-1} ) \label{phh} \\
\rho_{k} &= \mathbf{v}_{k-1}   \mathbf{v}_{k} \\
\sigma_{k} &= \mathbf{v}_{k-1}   \mathbf{hp}_{k} \\
\epsilon_{k} & =  \epsilon_{k-1}  \rho_{k} \\
\text{If\;}  &\epsilon_{k} \neq 1, \notag \\
& \tilde{c}_{k}  =  \epsilon_{k-1}^2 \rho_{k} \sigma_{k} \hat{h}_{k} / (1 - \epsilon_{k}^2) \\
& a_{k}^2 =  [(1 - \epsilon_{k-1}^2)(a_{k-1}^2(1 - \epsilon_{k}^2) - \hat{h}_{k}^2)] / (1 - \epsilon_{k}^2)^2 \label{lca} \\
& \mathbf{c}_{k}  =  \mathbf{c}_{k-1} + \hat{h}_{k}\mathbf{hp}_{k} + \tilde{c}_{k} \mathbf{v}_{k}  \label{pc}\\
& \mathbf{a}_{k}  =   \mathbf{c}_{k} + a_{k}\mathbf{v}_{k}   \label{pa} \\
\text{If\;}  &\epsilon_{k} > 1, \text{\;\;\;\;\;\;} b_{k}^2  =  - a_{k}^2(1 - \epsilon_{k}^2 )  = c_k^2-a_k^2 \label{lcb}\\
\text{If\;} &\epsilon_{k} < 1, \text{\;\;\;\;\;\;} b_{k}^2  =   a_{k}^2(1 - \epsilon_{k}^2 ) = a_k^2 - c_k^2  \label{lcba}\\
\text{If\;} &\epsilon_k = 1, \notag \\
&\tilde{c}_{k}  =  \epsilon_{k-1} \sigma_{k} \hat{h}_{k} / 2  \label{pch}\\
 &\hat{c}_{k} =  [(1 - \epsilon_{k-1}^2  \sigma_k^2) \hat{h}_{k}^2 + b_{k-1}^2] / (4 \tilde{c}_{k}) \label{pctld} \\
&\mathbf{c}_{k}  =  \mathbf{c}_{k-1} + \hat{h}_{k}\mathbf{hp}_{k} + \hat{c}_{k} \mathbf{v}_{k}  \label{pnc}
\end{alignat}
From expression (\ref{ph}), $\mathbf{hp}_2 = \mathbf{h}_2$,  and for each $k > 2$,  
$\mathbf{hp}_k$ is the orthogonal complement of the projection of $\mathbf{h}_{k}$ onto
the intersection of the hyperplanes $H_{2} \cap \ldots \cap H_{k-1}$.
Furthermore, the vectors $\mathbf{hp}_k$  are mutually orthogonal.
Given $\mathbf{hp}_k$, expression (\ref{pu}) computes the  vector $\mathbf{u}_{k-1}$ to be orthogonal to  $\mathbf{v}_{k-1}$,
and in the plane determined by  $\mathbf{hp}_k$ and  $\mathbf{v}_{k-1}$.
Expression (\ref{pvv}) computes the principal axis vector $\mathbf{v}_k$ for $B_{i_1,i_s}  \cap_{j=2}^{k-1}H_j \cap H_k$.
Expression (\ref{pd1}) computes the point $\mathbf{d}_k$ which lies on the principal axis vector $\mathbf{v}_k$ 
through the center of  $B_{i_1,i_s}  \cap_{j=2}^{k-1}H_j \cap H_k$.
Expression (\ref{phh}) computes the distance $\hat{h}$ from the center $\mathbf{c}_{k-1}$ of 
 $B_{i_1,i_s}  \cap_{j=2}^{k-1}H_j $ to the
point $\mathbf{d}_k$ on the major axis of  $B_{i_1,i_s}  \cap_{j=2}^{k-1}H_j \cap H_k$.

If $\epsilon_{k-1}\rho_k \neq 1$, then the intersection is either a hyperboloid or an ellipsoid, and expression (\ref{pc})
computes the center $\mathbf{c}_k$  and expression (\ref{pa}) computes the vertex (or vertices)  $\mathbf{a}_k$ based on Corollary A.1.

If $\epsilon_{k-1}\rho_k = 1$, then the intersection is a paraboloid, and expression (\ref{pnc})
computes the center $\mathbf{c}_k$  (which is also the vertex) based on (\ref{pch}) and (\ref{pctld}) from Corollary A.1.

Thus $B_S$ is determined by intersecting the hyperboloid $B_{i_1,i_{s}}$  with ${s}-2$ hyperplanes, 
resulting in the conic section  of dimension $n-{s}+2$,  with center  $\mathbf{c}_S$,  axis vector $\mathbf{v}_S$, vertex $\mathbf{a}_S$, 
eccentricity $\epsilon_S$, and parameters $a_S$ and $b_S$. 

Although $B_{i_1,i_s}$ is a hyperboloid, the conic section $B_S$ may be a hyperboloid or an ellipsoid.  
At each intersection of $B_{i_1,i_{s}}\cap_{j=2}^{k-1}H_j$ with $H_k$ for $k = 2,\ldots,s-1$,   
 the eccentricity $\epsilon_k$ is reduced, so that for some $k$, it may be that $\epsilon_k < 1$, in which case $B_S$ is an ellipsoid.

 Given the conic section $B_S$, the following property  shows that for any vector $\mathbf{u}_S$ orthogonal to $\mathbf{v}_S$, 
there is a two-dimensional conic section in  aff$(\mathbf{c}_S, \mathbf{v}_S ,\mathbf{u}_S)$, with the same parameters and vectors as $B_S$.  

\begin{property}
Given the conic section $B_S$,
let $\mathbf{u}_S$ be any unit vector orthogonal to the axis vector $\mathbf{v}_S$.
If $B_S$ is a hyperboloid,
then $\{ \mathbf{y}(\alpha) = \mathbf{c}_S+ a_S\sec(\alpha) \mathbf{v}_S + b_S \tan (\alpha) \mathbf{u}_S, - \pi/2 < \alpha < \pi/2 \}$,
 with $b_S^2 = \sqrt{c_S^2 - a_S^2}$,
 is one sheet of a two-dimensional hyperbola in the affine plane aff$( \mathbf{c}_S, \mathbf{v}_S, \mathbf{u}_S )$,  with the same vectors and parameters as $B_S$. 
 If $B_S$ is an ellipsoid,
then $\{ \mathbf{y}(\alpha) = \mathbf{c}_S+ a_S\cos(\alpha) \mathbf{v}_S + b_S \sin (\alpha) \mathbf{u}_S, - \pi/2 < \alpha < \pi/2 \}$,
 with $b_S^2 = \sqrt{a_S^2 - c_S^2}$,
 is  a two-dimensional ellipse in the affine plane aff$( \mathbf{c}_S, \mathbf{v}_S, \mathbf{u}_S )$,  with the same vectors and parameters as $B_S$. 
 If $B_S$ is a paraboloid, then $\{ \mathbf{x}(\alpha) = \mathbf{c} + c\alpha^2\mathbf{v}_S + 2c\alpha\mathbf{u}_S, -\infty < \alpha < \infty \}$,
 is a two-dimensional parabola in the affine plane aff$(\mathbf{c}_S,\mathbf{v}_S, \mathbf{u}_S)$, with the same vectors and parameters as $B_S$.
 \end{property}
 


\end{document}